\documentclass[3p]{elsarticle}
\biboptions{numbers,sort&compress}
\usepackage{amssymb}
\usepackage{amsmath}
\usepackage{booktabs}
\usepackage{algorithm}
\usepackage{algorithmic}
\usepackage{subcaption}
\usepackage{multirow}
\usepackage{makecell}
\usepackage{adjustbox}
\usepackage{hyperref}

\usepackage{float}

\usepackage{xcolor}

\journal{Applied Mathematical Modelling}
\begin{document}

\begin{frontmatter}



\title{Physics-Informed Kolmogorov-Arnold Networks for multi-material elasticity problems in electronic packaging}


\author[address1]{Yanpeng Gong}
\ead{yanpeng.gong@bjut.edu.cn}
\author[address1]{Yida He}
\author[address2]{Yue Mei\corref{mycorrespondingauthor}}
\ead{meiyue@dlut.edu.cn}
 \author[address3,address4]{Xiaoying Zhuang\corref{mycorrespondingauthor}}
 \ead{zhuang@iop.uni-hannover.de}
\author[address1]{Fei Qin}
 \author[address5]{Timon Rabczuk}
 % addresses
\address[address1]{Institute of Electronics Packaging Technology and Reliability, Department of Mechanics, Beijing University of Technology, Beijing, 100124, China}
\address[address2]{State Key Laboratory of Structural Analysis for Industrial Equipment, Department of Engineering Mechanics, International Research Center for Computational Mechanics, Dalian University of Technology,Dalian 116023, China}
\address[address3]{Chair of Computational Science and Simulation Technology, Institute of Photonics, Department of Mathematics and Physics, Leibniz University Hannover, 30167 Hannover, Germany}
\address[address4]{Department of Geotechnical Engineering, College of Civil Engineering, Tongji University, Shanghai, 200092, China}
\address[address5]{Institute of Structural Mechanics, Bauhaus University Weimar, 99423 Weimar, Germany}

\cortext[mycorrespondingauthor]{Corresponding author}

% \affiliation{organization={},%Department and Organization

%% Abstract
\begin{abstract}

This paper proposes a Physics-Informed Kolmogorov-Arnold Network for analyzing elasticity problems in multi-material electronic packaging structures. The method replaces traditional Multi-Layer Perceptrons with Kolmogorov-Arnold Networks within an energy-based Physics-Informed Neural Network framework. By constructing admissible displacement fields satisfying essential boundary conditions and optimizing network parameters through numerical integration, the proposed method effectively handles material property discontinuities. Unlike traditional methods that require domain decomposition and interface constraints for multi-material problems, Kolmogorov-Arnold Networks’ trainable B-spline activation functions provide inherent piecewise characteristics. This capability stems from B-splines’ local support, which enables effective approximation of discontinuities despite their individual smoothness. Consequently, this approach enables  accurate approximation across the entire domain using a single network and simplifying the computational framework. Numerical experiments demonstrate that the proposed method achieves excellent accuracy and robustness in multi-material elasticity problems, validating its practical potential for electronic packaging analysis. Source codes are available at https://github.com/yanpeng-gong/PIKAN-MultiMaterial.
\end{abstract}


% \begin{highlights}
% \end{highlights}

%% Keywords
\begin{keyword}
Physics-Informed Neural Network \sep Kolmogorov-Arnold networks \sep Multi-material structures \sep Electronic packaging structures \sep Deep energy method


\end{keyword}

\end{frontmatter}



\section{Introduction}

Electronic packaging provides critical mechanical support, electrical interconnection, and thermal management for high-density integrated electronic devices~\cite{moreau2012semiconductor}. Its reliability analysis involves challenges including multi-material systems, multi-scale phenomena, and interface discontinuities. Numerical simulation techniques, benefiting from easy implementation, high accuracy, and environment-independent operation, have become mainstream methods for electronic packaging reliability analysis, including finite element method (FEM)~\cite{kim2020analyses,feng2023engineered}, boundary element method (BEM)~\cite{gong2022isogeometric,yu2021thermal}, coupled algorithms~\cite{gong2025coupled,qin2022application}, finite difference method (FDM)~\cite{zhang2004finite}, and phase-field models~\cite{gong2024application,GONG2025111039}, among others, with each method advancing through specific applications. However, facing the growing design complexity and computational efficiency demands of modern packaging structures, traditional methods still have inherent limitations in computational cost, local field resolution, and physical modeling, creating urgent demand for more efficient and flexible computational frameworks.

The development of artificial intelligence (AI) technology and enhanced computational resources have provided opportunities for deep learning as a new approach to electronic packaging reliability problems. Researchers have successfully applied machine learning to construct surrogate models for rapid thermal-mechanical coupling behavior assessment~\cite{kang2023machine}, geometric parameter co-design~\cite{zaghari2025co}, and process prediction~\cite{yu2025development}, demonstrating significant computational efficiency advantages. However, these data-driven methods typically require substantial high-confidence training data; in scenarios with limited engineering data, their predictions may deviate from physical laws, affecting model generalization and physical consistency~\cite{yuan2022towards}. Therefore, combining machine learning with traditional computational methods offers tremendous potential for transforming modeling and simulation in engineering and scientific domains, opening new pathways for next-generation computer-aided analysis and design~\cite{rabczuk2024machine}.

In recent years, deep learning methods have provided new numerical paradigms for solving solid mechanics problems by avoiding traditional mesh generation and improving computational efficiency~\cite{WOS:001223472000001,WOS:001068952500001}. To reduce neural networks' dependence on high-quality data, researchers embed physical laws into loss functions to enhance models through physics-informed constraints. Physics-Informed Neural Networks (PINNs), proposed by Raissi et al.~\cite{raissi2019physics} in 2019, became the pioneering work in this direction. The implementation of PINNs benefits from the mature application of automatic differentiation functionality in machine learning frameworks such as PyTorch and TensorFlow. Compared to traditional numerical algorithms, PINNs leverage the powerful approximation capabilities of neural networks to provide a new paradigm for solving complex problems that conventional numerical methods struggle to address~\cite{cuomo2022scientific}. Furthermore, PINNs demonstrate significant potential in handling high-dimensional problems, as their loss functions are typically constructed through Monte Carlo sampling, which is inherently an effective tool for high-dimensional integration. PINNs encompass two forms: strong-form original PINNs and energy-form PINNs. Samaniego et al.~\cite{samaniego2020energy} introduced the Deep Ritz Method proposed by E and Yu~\cite{yu2018deep} into computational mechanics, forming energy-form PINNs called the Deep Energy Method (DEM). Strong-form PINNs construct loss functions by directly combining PDEs through weighted residual methods, while DEM constructs loss functions based on variational principles such as the minimum potential energy principle in mechanics. The flexibility and potential of the PINN framework continue to drive its development and application in various cutting-edge engineering fields~\cite{WOS:001223472000001,MA2025107559,WOS:001523847800001,CHEN2025109735,AMININIAKI2021,XIONG2025,WOS:001571017400001}. In electronic packaging structure analysis, although there are limited studies using traditional data-driven surrogate models based on physical simulation data~\cite{yao2022physics}, research directly applying PINNs to such simulations has rarely been reported.

Within the standard PINN framework, multi-layer perceptrons have been widely adopted due to their structural simplicity and universal approximation properties~\cite{hornik1989multilayer}. However, the framework is inherently flexible, and the approximation function is not limited to MLPs. Research has successfully extended to various architectures such as RBFs~\cite{GAO2025127517}, Wavelets \cite{ZHU2026105174}, and others. Recently, Kolmogorov-Arnold Networks (KANs), proposed by Liu et al.~\cite{liu2024kan}, employ trainable one-dimensional activation functions that simultaneously replace linear weights and provide nonlinear activation behavior, with nodes simplified to perform only input summation, have also been introduced into the PINN framework, rapidly spawning multiple improved studies focusing on fair comparisons with MLPs~\cite{SHUKLA2024117290}, architectural innovations~\cite{WANG2025117518,WOS:001575943800001}, and other aspects. When dealing with PDE problems involving high-dimensional or complex solution spaces, MLPs often require extremely wide or deep networks to achieve the required accuracy, leading to an enormous number of parameters. This increases optimization difficulty and computational cost~\cite{han2018solving}, while KANs provide an alternative parameterization approach through their spline-based learnable activation functions. Meanwhile, the spectral bias problem in MLP training may lead to insufficient learning capability for high-frequency components of PDE solutions, whereas the flexible function representation provided by KANs has been indicated to be more conducive to capturing multi-scale features of solutions and alleviating spectral bias phenomena~\cite{wang2025express,meshir2025study}. This provides a new potential avenue for solving PDE problems with multi-scale characteristics or solutions with sharp transitions. It should be noted that the completeness of the universal approximation theory for classical KANs is still being refined, with only a few variants rigorously proven to date~\cite{TOSCANO2025107831}. Nevertheless, the excellent computational performance demonstrated by this architecture in a wide range of scientific computing problems motivates us to further explore its application potential in PDE solving. 

Electronic packaging structures are typically composed of heterogeneous materials with different properties, exhibiting multi-material, multi-interface structural characteristics. For multi-material structures with sharp interfaces, the meshfree nature of PINNs enables them to avoid the mesh distortion issues that may arise in FEM. When applying PINNs to solve multi-material problems, the current mainstream approach uses different neural networks to describe physical fields in different material domains. Diao et al.~\cite{diao2023solving} proposed a strong-form PINNs framework for solving multi-material solid mechanics problems by partitioning the computational domain through domain decomposition, utilizing independent sub-networks to characterize mechanical behavior of different materials, and introducing interface regularization terms to handle contact relationships and discontinuities. Wang et al.~\cite{wang2022cenn} focused on variational problems and proposed a conservative energy method based on subdomain neural networks for solving heterogeneous material and complex geometry problems. However, PINNs based on subdomain partitioning have inherent limitations. To accurately describe interface contact relationships, additional hyperparameters must be introduced, which typically rely on empirical parameter tuning. Although research has explored hyperparameter optimization strategies~\cite{sheng2021pfnn,wang2022and}, obtaining optimal parameter combinations remains challenging when solving specific problems. Furthermore, the domain decomposition method requires independent optimization of neural networks for each subdomain, which increases training complexity and often causes the optimization process to converge to local rather than global optima.

Although Physics-Informed Kolmogorov-Arnold Networks (PIKANs) have received increasing attention~\cite{liu2024kan,SHUKLA2024117290,WANG2025117518,WOS:001575943800001, toscano2025pinns, mostajeran2025, faroughi2025}, the distinctive contribution of this study lies in developing a straightforward PIKAN solution framework specifically tailored for multi-material solid mechanics problems with various interface configurations commonly encountered in electronic packaging structures, integrating KANs with the DEM to replace traditional MLPs. To our knowledge, we present and validate for the first time this integrated approach for electronic packaging multi-material structure analysis, focusing on the persistent challenge of modeling multi-material structures with inherent interface discontinuities. Within the PINN framework, the strong form has broader generality, but its application typically requires extensive hyperparameter tuning. In contrast, the Deep Energy Method, based on variational principles, transforms the governing equations into an extremum problem of energy functionals. This not only reduces the number of required hyperparameters, but also typically achieves higher computational accuracy and efficiency due to the reduced order of functional derivatives~\cite{li2021physics}, though its generality is usually weaker than the strong form. For elastic analysis problems of multi-material structures, this study aims to develop a more concise solution paradigm, and therefore adopts the DEM with fewer hyperparameters. It should be noted that theoretical analysis of numerical convergence and stability for PINNs and their variational (energy-based) forms has made certain progress~\cite{xu2025convergence,WOS:001601418300001}. For instance, rigorous convergence proofs have been established for specific higher-order multi-scale PINN methods in composite mechanics~\cite{WOS:001305621200001}. However, for the emerging KAN architecture and its applications in physics-informed learning, rigorous convergence theory for KAN/PIKAN architectures remains an active research frontier. 

Based on the above research background, this paper employs the PIKAN method to solve multi-material problems in solid mechanics and extends its application to electronic packaging multi-material structure analysis. This work constructing admissible displacement fields~\cite{wang2022cenn,sheng2021pfnn} that strictly satisfy essential boundary conditions. The main contributions of our work include: 
\begin{itemize}
    \item \textbf{Development of PIKAN methodology} that replaces traditional MLPs with KANs in deep energy method, leveraging KANs’ inherent piecewise function characteristics to naturally accommodate material property discontinuities.
    \item \textbf{Elimination of domain decomposition requirements} by using a single KAN to approximate displacement fields across the entire computational domain without subdomain partitioning or interface continuity constraints.
    \item \textbf{Comprehensive validation through diverse examples} including cantilever beams with different interface geometries and other electronic packaging models, demonstrating superior accuracy compared to traditional domain decomposition approaches.
    \item \textbf{Systematic investigation of numerical integration schemes and hyperparameter sensitivity analysis}, with complete algorithm implementation details and source code to be provided in our GitHub repository, facilitating practical application, comparison, and further development of PIKANs for multi-material problems.
\end{itemize}

This paper is organized as follows: Section \ref{sec:DEM overview} presents the governing equations for elasticity problems in solid mechanics, the DEM, and the construction method for admissible displacement fields. Section~\ref{sec:PIKAN-methodology} provides a detailed description of the proposed PIKAN method. Section~\ref{sec:multi_material_elasticity_problems} presents specific strategies for applying PIKAN to solve multi-material problems. Section~\ref{sec:numerical_examples} validates the effectiveness of the proposed method through multi-material model examples and applies it to electronic packaging multi-material structure analysis. Finally, Section~\ref{sec:conclusion} summarizes the research conclusions.

\section{Deep energy method-based computational solid mechanics}
\label{sec:DEM overview}
%% Labels are used to cross-reference an item using \ref command.
This section introduces the Deep Energy Method (DEM) for solving elasticity problems. Then we present the construction of admissible displacement fields that naturally satisfy essential displacement boundary conditions.

\subsection{Deep energy method for solving elasticity problems}

For elastic systems, based on the principle of minimum potential energy in variational elasticity theory, the Deep Energy Method (DEM) constructs the total potential energy functional of the system as a loss function for optimization \cite{samaniego2020energy}. For problems without body forces, the specific mathematical expression of this loss function is
\begin{equation}
\label{eq:energy_loss}
\mathcal{L} = \Psi_{\text{in}} - \Psi_{\text{ex}} = \int_{\Omega} \frac{1}{2} \sigma_{ij} \varepsilon_{ij} \mathrm{d}\Omega - \int_{\Gamma_\sigma} \bar{p}_i u_i \mathrm{d}\Gamma
\end{equation}
where $\Psi_{\text{in}}$ and $\Psi_{\text{ex}}$ represent the internal strain energy and external force potential energy, respectively;  $\sigma_{ij}$ and $\varepsilon_{ij}$ denote the stress and strain tensors, respectively; $u_i$ represents the $i$-th displacement component; and $\bar{p}_i$ is the traction applied on the natural boundary. DEM shares similarities with FEM in adopting energy-based variational weak forms, but its key distinction lies in the fact that DEM directly optimizes the total potential energy functional, thereby eliminating the steps of preliminary derivation and discretization of the weak form of governing equations required by the FEM. For the electronic packaging multi-material structures of interest, we focus on mechanical responses dominated by external mechanical loads and material interface effects. The body force term is omitted in the energy functional construction, based on engineering simplification appropriate for the scale and loading characteristics of the studied problems \cite{BaiAn2023,WOS:000793133700005}. The proposed framework can straightforwardly incorporate body force terms to handle more general scenarios when needed.

Among all kinematically admissible displacement fields $\mathbf{u}(\mathbf{x})$ that satisfy the essential boundary conditions, we seek the solution that minimizes the total potential energy functional
\begin{equation}
\label{eq:variational_problem}
\small
\mathbf{u}(\mathbf{x}) = \underset{\mathbf{u}(\mathbf{x})}{\arg\min} \mathcal{L}(\mathbf{u}(\mathbf{x})) = \underset{\mathbf{u}(\mathbf{x})}{\arg\min} \left\{\int_{\Omega} \psi(\boldsymbol{\varepsilon}(\mathbf{u}(\mathbf{x}))) \mathrm{d}\Omega - \int_{\Gamma_\sigma} \overline{\mathbf{p}} \cdot \mathbf{u}(\mathbf{x}) \mathrm{d}\Gamma\right\}
\end{equation}
where $\mathbf{u}(\mathbf{x})$ is approximated through a neural network as $\mathbf{u}(\mathbf{x}) \approx \mathbf{F}(\mathbf{x}; \boldsymbol{\theta})$, $\mathbf{F}(\cdot)$ denotes the neural network mapping, and $\psi(\boldsymbol{\varepsilon}(\mathbf{u}(\mathbf{x})))$ is the strain energy density function, which is a function of the strain tensor $\boldsymbol{\varepsilon}(\mathbf{u}(\mathbf{x}))$.

\subsection{Construction of admissible displacement fields}

The DEM selects the displacement field that minimizes the potential energy among all admissible displacement fields. The key to realizing admissible displacement fields lies in constructing displacement trial functions that naturally satisfy essential boundary conditions and embedding this constraint directly into the network architecture to ensure that the displacement field consistently satisfies essential boundary conditions throughout the optimization process. 

The general structure of an admissible displacement field can be expressed as
\begin{equation}
\label{eq:admissible_displacement}
\mathbf{u}^{\text{pred}}(\mathbf{x}) = \mathbf{P}(\mathbf{x}) + \mathbf{D}(\mathbf{x}) \odot \mathbf{F}(\mathbf{x}; \boldsymbol{\theta})
\end{equation}
where $\mathbf{P}(\mathbf{x})$ is the boundary function that satisfies essential boundary conditions. Note that for fixed boundaries, the prescribed essential boundary condition value is \textbf{0}, i.e., $\mathbf{P}(\mathbf{x})=\mathbf{0}$. $\mathbf{F}(\mathbf{x}; \boldsymbol{\theta})$ is the generalized network obtained through training by minimizing the potential energy loss function. $\odot$ denotes element-wise multiplication. $\mathbf{D}(\mathbf{x})$ is the distance function, representing the shortest distance from the current coordinate point to the essential boundary $\Gamma_u$ 

\begin{equation}
    \label{eq:distance_function}
    \mathbf{D}(\mathbf{x}) = \min_{\mathbf{y} \in \Gamma_u} \|\mathbf{x} - \mathbf{y}\|
    \end{equation}

The boundary and distance functions have flexible implementation forms; in addition to neural networks, fitting methods such as radial basis functions can also be employed~\cite{wang2022cenn}. The key advantage of this construction is that when the coordinate point is located on the boundary $\Gamma_u$, the distance function $\mathbf{D}(\mathbf{x})=\mathbf{0}$, and the admissible displacement field obtained from Eq.~\eqref{eq:admissible_displacement} naturally satisfies the enforced essential boundary conditions. In the validation examples presented in Section~\ref{sec:numerical_examples} of this study, given that the essential boundary geometries are relatively simple, $\mathbf{D}(\mathbf{x})$ consistently adopts the analytical distance function form~\cite{SUKUMAR2022114333}, thereby avoiding additional errors introduced by neural network approximation while accurately enforcing essential boundary conditions. This method enables the admissible displacement field to naturally satisfy boundary constraints during generalized network training, without requiring the introduction of boundary condition penalty terms and associated hyperparameter tuning in its loss function. For detailed discussion on admissible displacement fields, please refer to~\cite{SUKUMAR2022114333}.

The neural network parameters of the generalized network $\mathbf{F}(\mathbf{x}; \boldsymbol{\theta})$ are trained and optimized by minimizing the potential energy loss function
\begin{equation}
    \label{eq:generalized_network_optimization}
    \begin{aligned}
    \mathbf{F}(\mathbf{x}; \boldsymbol{\theta}) & = \underset{\boldsymbol{\theta}}{\arg\min} \mathcal{L} \left(\mathbf{u}^{\text{pred}}(\mathbf{x}; \boldsymbol{\theta})\right) \\
    & = \underset{\boldsymbol{\theta}}{\arg\min} \left\{\int_{\Omega} \psi\left(\boldsymbol{\varepsilon}(\mathbf{u}^{\text{pred}}(\mathbf{x}; \boldsymbol{\theta}))\right) \mathrm{d} \Omega - \int_{\Gamma_\sigma} \overline{\mathbf{p}} \cdot \mathbf{u}^{\text{pred}}(\mathbf{x}; \boldsymbol{\theta}) \mathrm{d}\Gamma\right\}
    \end{aligned}
  \end{equation}
Note that $\mathbf{F}(\mathbf{x}; \boldsymbol{\theta})$ represents the generalized network output, while the actual displacement field is given by Eq.~\eqref{eq:admissible_displacement}. 

\section{Methodology of PIKAN}
\label{sec:PIKAN-methodology}

The previous section provided a detailed introduction to using DEM for solving elasticity problems. In this work, instead of using the MLP traditionally employed in DEM, we adopt KAN as the neural network architecture. This section first introduces the fundamental principles of KAN, then presents the PIKAN method and finally describes the sample point distribution strategies employed in the analysis.
\begin{figure}[htbp]
  \centering
  \includegraphics[width=0.75\textwidth]{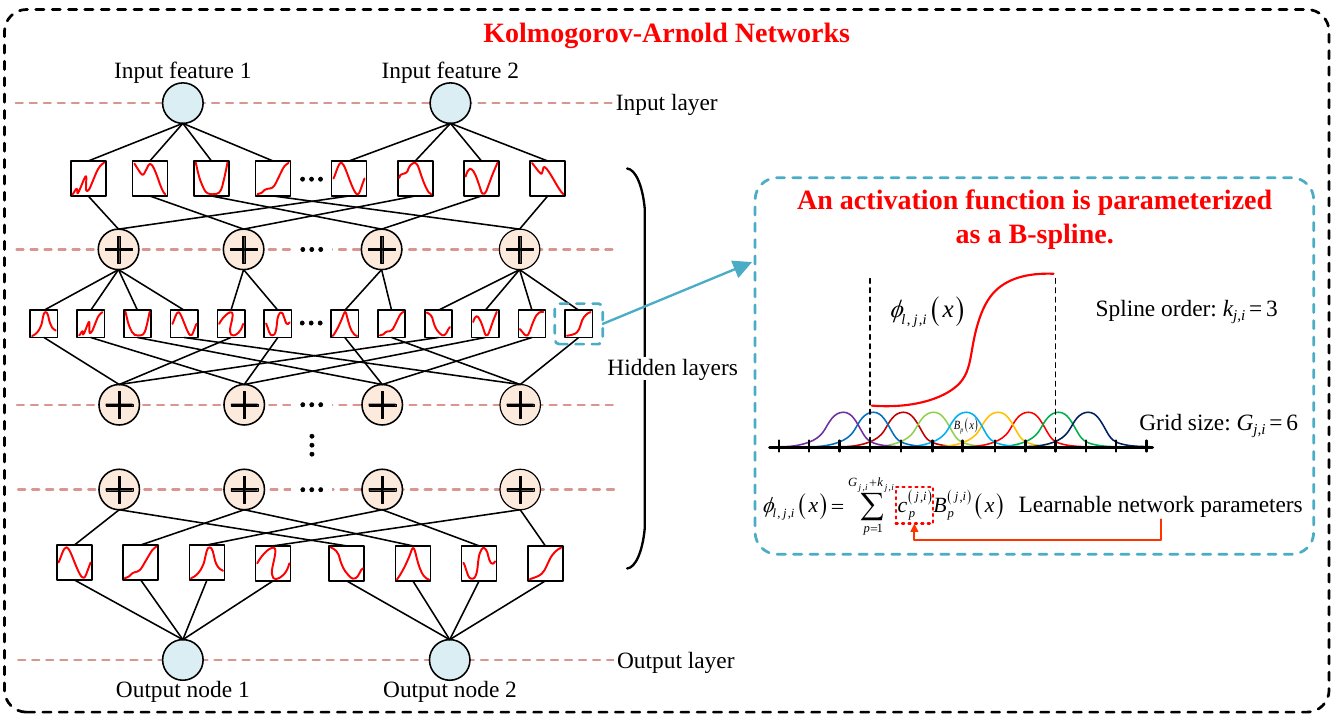}
  \caption{Kolmogorov-Arnold Network schematic (adapted and redrawn from the original KAN architecture~\cite{liu2024kan}).}
  \label{fig:kan_schematic}
\end{figure}

\subsection{Kolmogorov-Arnold Networks (KANs)}
\label{sec:kan_overview}

KANs are based on the Kolmogorov-Arnold representation theorem and feature a fully connected structure. Traditional MLPs achieve layer-wise transformations by applying fixed nonlinear activation functions to linear combinations of inputs at each node, where weights and biases are learned parameters. During backpropagation training, MLPs update these parameters by computing gradients of the loss function with respect to weights and biases.

In contrast, KANs move activation functions from nodes to edges, where each edge contains  learnable univariate functions parameterized as B-spline functions. Instead of applying fixed activation functions at nodes after linear transformations (as in MLPs), KANs apply learnable activation functions on each edge connection. Nodes in KANs perform only summation operations on the outputs from incoming edges. Their core advantage lies in adaptive learning capability for complex functional relationships and the resulting potentially higher parameter efficiency (i.e., achieving comparable accuracy with fewer parameters)~\cite{liu2024kan}. It must be objectively noted that choosing KANs over MLPs involves trade-offs. In many benchmark tests, their accuracies are often comparable~\cite{SHUKLA2024117290}. The primary advantage of KANs lies in potentially mitigating spectral bias and more efficiently learning high-frequency features~\cite{meshir2025study,wang2025express}, but the more complex gradient computations of B-spline functions during training increase computational cost~\cite{liu2024kan}. The piecewise function approximation characteristics of KANs are better suited for handling solution field discontinuities at multi-material interfaces, which traditional MLPs struggle to naturally capture.  The KAN architecture is illustrated in Fig.~\ref{fig:kan_schematic}.

The core idea of the Kolmogorov-Arnold Representation Theorem (KART) is that any continuous multivariate function f defined on a bounded domain can be expressed as a composition of finite univariate continuous functions and addition operations
\begin{equation}
\label{eq:kolmogorov_arnold_theorem}
f(\mathbf{x}) = f(x_1, x_2, \ldots, x_n) = \sum_{h=1}^{2n+1} \Phi_h\left(\sum_{g=1}^n \phi_{h,g}(x_g)\right)
\end{equation}
This provides theoretical inspiration for expressing multivariate functions through univariate function compositions. The design of KANs parameterizes the univariate functions $\phi$ and $\Phi$ in the theorem as learnable B-spline functions and extends them into nested, multi-layer network structures to enhance expressiveness. However, this deep KAN architecture goes beyond the scope of the original KART. As a new network paradigm, its theoretical validity is established and explored by the approximation theorem proposed and developed in the original paper itself~\cite{liu2024kan}. This work focuses on the application of this architecture, and the following will introduce its network structure and implementation details in detail.

% Define a KAN layer with $N_{\text{in}}$ input neurons and $N_{\text{out}}$ output neurons. For a general KAN layer $l$, the activation functions are denoted as
% \begin{equation} 
% \label{eq:kan_layer} 
% \Phi_l = \{\phi_{l,j,i}\}, \quad i = 1,2,\ldots,N_l; \quad j = 1,2,\ldots,N_{l+1}
% \end{equation}
% where each $\phi_{l,j,i}$ represents a trainable univariate function.

Define a KAN layer with $N_{\text{in}}$ input neurons and $N_{\text{out}}$ output neurons
\begin{equation}
\label{eq:kan_layer}
\mathbf{\Phi} = \{\phi_{h,g}\}, \quad g = 1,2,3,\ldots,N_{\text{in}}; \quad h = 1,2,3,\ldots,N_{\text{out}}
\end{equation}
where $\phi_{h,g}$ has trainable parameters. The original Eq.~\eqref{eq:kolmogorov_arnold_theorem} represents a combination of two KAN layers: the inner function with $N_{\text{in}} = n, N_{\text{out}} = 2n+1$, and the outer function with $N_{\text{in}} = 2n+1, N_{\text{out}} = 1$. By stacking more KAN layers, deeper KANs can be obtained. The shape of a KAN is represented by an integer array $[N_0, N_1, N_2, \ldots, N_L]$, where $N_i$ is the number of neurons in the $i$-th layer. Between the $l$-th and $(l+1)$-th layers, there are $N_l N_{l+1}$ activation functions. The activation function connecting the $i$-th neuron $(l,i)$ in the $l$-th layer and the $j$-th neuron $(l+1,j)$ in the $(l+1)$-th layer is denoted as
\begin{equation}
\label{eq:activation_function_notation}
\phi_{l,j,i}, \quad l = 0,1,\ldots,L-1; \quad i = 1,2,\ldots,N_l; \quad j = 1,2,\ldots,N_{l+1}
\end{equation}

The activation value of the $(l+1)$-th layer is computed as
\begin{equation} 
\label{eq:layer_computation} 
x_{l+1,j} = \sum_{i=1}^{N_l} \phi_{l,j,i}(x_{l,i}), \quad j = 1,2,\ldots,N_{l+1}
\end{equation}
In vector form, this can be written as
\begin{equation}
\label{eq:vector_form}
\mathbf{x}_{l+1} = \begin{bmatrix} 
\sum_{i=1}^{N_l} \phi_{l,1,i}(x_{l,i}) \\ 
\sum_{i=1}^{N_l} \phi_{l,2,i}(x_{l,i}) \\ 
\vdots \\ 
\sum_{i=1}^{N_l} \phi_{l,N_{l+1},i}(x_{l,i}) 
\end{bmatrix} = \mathbf{\Phi}_l(\mathbf{x}_l)   
\end{equation}
where $\mathbf{\Phi}_l$ represents the function matrix of the $l$-th KAN layer.

Since the activation functions are parameterized as B-spline functions, each activation function $\phi_{l,j,i}$ can be expressed as
\begin{equation}
\label{eq:bspline_form}
\phi_{l,j,i}(x) = \sum_{p=1}^{G_{j,i}+k_{j,i}} c_p^{(j,i)} B_p^{(j,i)}(x)
\end{equation}
where $G_{j,i}$ and $k_{j,i}$ are the grid size and spline order for the activation function connecting neuron $(l,i)$ to neuron $(l+1,j)$, $c_p^{(j,i)}$ are the trainable B-spline coefficients, and $B_p^{(j,i)}(x)$ are the corresponding B-spline basis functions.

The total number of trainable parameters in layer $l$ is
\begin{equation}
N_{\text{params},l} = \sum_{i=1}^{N_l} \sum_{j=1}^{N_{l+1}} (G_{j,i} + k_{j,i})
\end{equation}

To enhance the function approximation capability, scaling coefficients are introduced for each activation function~\cite{WANG2025117518}. For each activation function $\phi_{l,j,i}$, a corresponding scaling coefficient $m_{j,i}$ is defined. The scaling coefficient matrix $\mathbf{M}_l$ is
\begin{equation}
    \label{eq:scaling_matrix}
    \mathbf{M}_l = \begin{bmatrix}
    m_{1,1} & m_{1,2} & \cdots & m_{1, N_l} \\
    m_{2,1} & m_{2,2} & \cdots & m_{2, N_l} \\
    \vdots & \vdots & \ddots & \vdots \\
    m_{N_{l+1}, 1} & m_{N_{l+1}, 2} & \cdots & m_{N_{l+1}, N_l}
    \end{bmatrix}
\end{equation}
The scaled activation functions become
\begin{equation}
\label{eq:scaled_activation}
\phi_{l,j,i}^*(x) = m_{j,i} \cdot \phi_{l,j,i}(x)
\end{equation}

Additionally, a linear transformation matrix $\mathbf{W}_l$ and nonlinear activation function $\sigma(\cdot)$ are introduced to further enhance the network's expressive and fitting capabilities. The final output of the KAN layer is computed as
\begin{equation}
\label{eq:kan_final_output}
x_{l+1,j} = \sum_{i=1}^{N_l} m_{j,i} \cdot \phi_{l,j,i}(x_{l,i}) + \sum_{i=1}^{N_l} w_{j,i} \cdot \sigma(x_{l,i}), \quad j = 1,2,\ldots,N_{l+1}
\end{equation}
where $w_{j,i}$ are elements of the weight matrix $\mathbf{W}_l$. The introduction of $\sigma(\cdot)$ enables the network to capture smoother solutions, complementing the B-spline functions. The scaling coefficients $m_{j,i}$ and weights $w_{j,i}$ are both trainable parameters optimized during training.

Following common practice to simplify implementation and training, this work adopts a uniform configuration where all activation functions within each KAN layer share the same grid size $G_{j,i}$ and B-spline order $k_{j,i}$. The trainable parameters in KAN layers are given in Table~\ref{tab:kan_hyperparameters}.
\begin{table}[htbp]
\tiny
  \centering
  \caption{Trainable parameters of KAN layers}
  \label{tab:kan_hyperparameters}
  \begin{tabular}{llll}
  \toprule
  Parameter & Type & Quantity per layer & Description \\
  \midrule
  $c_p^{(j,i)}$ & B-spline coefficients & $N_l \times N_{l+1} \times (G_{j,i}+k_{j,i})$ & Coefficients for B-spline basis functions \\
  $m_{j,i}$ & Scaling factors & $N_l \times N_{l+1}$ & Element-wise scaling factors for activation functions \\
  $w_{j,i}$ & Linear weights & $N_l \times N_{l+1}$ & weights for the nonlinear activation function $\sigma(\cdot)$\\
  \bottomrule
  \end{tabular}
\end{table}
For a KAN network with $L$ layers, given an input vector $\mathbf{x}^{(0)} \in \mathbb{R}^{N_0}$, the forward propagation through the network is computed sequentially as
\begin{align}
\mathbf{x}^{(L)}=\mathbf{\Phi}_{L-1}^*\left(\mathbf{\Phi}_{L-2}^*\left(\cdots\mathbf{\Phi}_0^*\left(\mathbf{x}^{(0)}\right)\cdots\right)\right)
\end{align}
where $\mathbf{\Phi}_l^*(\mathbf{x}^{(l)})=\mathbf{\Phi}_l(\mathbf{x}^{(l)}) \odot \mathbf{M}_l+\mathbf{W}_l \cdot \sigma(\mathbf{x}^{(l)})$, $0\leq l \leq (L-1)$ and $\mathbf{x}^{(L)}$ is the network's final output.

\subsection{Physics-Informed Kolmogorov-Arnold Networks (PIKANs)}
\label{sec:pikan_overview}

\begin{figure}[htbp]
    \centering
    \includegraphics[width=0.75\textwidth]{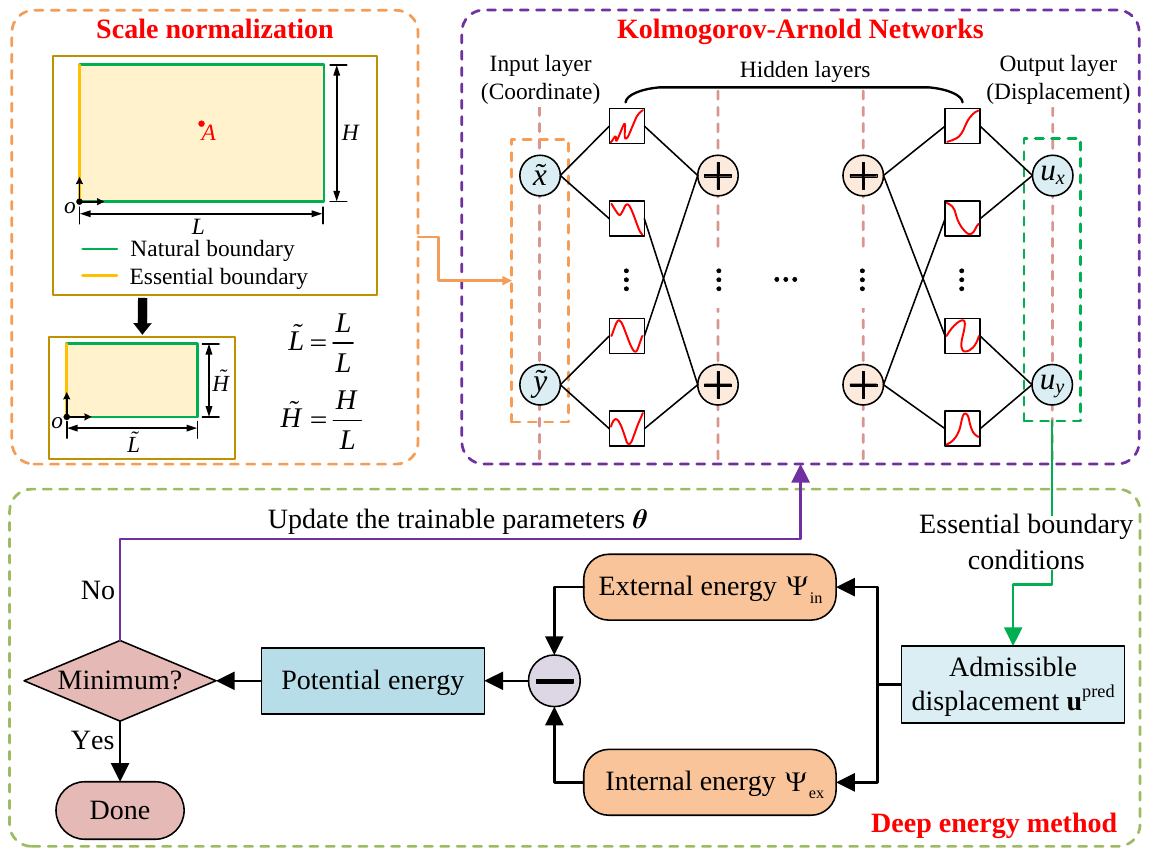}
    \caption{PIKAN schematic: Deep Energy Method using KAN instead of MLP.}
    \label{fig:pikan_schematic}
\end{figure}

As discussed in the Introduction, PIKAN has emerged as a widely explored new paradigm in the field of physics-informed machine learning~\cite{liu2024kan,SHUKLA2024117290,WANG2025117518,WOS:001575943800001, toscano2025pinns, mostajeran2025, faroughi2025}. This section presents the specific PIKAN framework adopted in this study. As illustrated in Fig.~\ref{fig:pikan_schematic}, KAN replaces the MLP traditionally used in the Deep Energy Method for computational solid mechanics.

Since KANs employ B-spline activation functions defined on preset grids ([0,1] or [-1,1]), and the algorithm takes sample point coordinates within the computational domain as neural network inputs, while actual computational models have arbitrary geometric dimensions that may exceed the preset grid range, normalization of network inputs is necessary. Geometric normalization is crucial for ensuring training stability and convergence speed~\cite{toscano2025pinns,mostajeran2025,wang2023expert}, and its theoretical necessity has been thoroughly explored through rigorous analyses based on frameworks such as Neural Tangent Kernel (NTK), particularly in studies of Chebyshev KANs~\cite{mostajeran2025,faroughi2025}. This work implements this standardization strategy in specific problems based on these established theoretical foundations.

As shown in the upper left corner of Fig.~\ref{fig:pikan_schematic}, we consider a 2D geometric domain with characteristic length \textit{$L$} and characteristic height $H$, where $L \geq H$, and point \textit{o} is the coordinate origin. For any point $A(x_A, y_A)$ within the domain, we perform coordinate normalization
\begin{equation}
\label{eq:coordinate_normalization}
\tilde{x}_A = \frac{x_A}{L}, \quad \tilde{y}_A = \frac{y_A}{L}
\end{equation}
The normalized coordinates serve as actual inputs to the KAN, ensuring B-spline activation functions operate effectively within their definition domain.

To mitigate training instability caused by KAN layer outputs potentially exceeding the definition domain of subsequent spline functions, this work adopts a validated approach of applying tanh activation functions to constrain intermediate layer outputs~\cite{SHUKLA2024117290,mostajeran2025}
\begin{equation}
\label{eq:kan_output_with_tanh}
x_{l+1,j}^{\text{new}} = \tanh\left(x_{l+1,j}\right) = \tanh\left(\sum_{i=1}^{N_l} m_{j,i} \cdot \phi_{l,j,i}(x_{l,i}) + \sum_{i=1}^{N_l} w_{j,i} \cdot \sigma(x_{l,i})\right)
\end{equation}
for $j = 1, 2, \ldots, N_{l+1}$ and layers $l = 0, 1, \ldots, L-2$. Note that the final layer ($L-1$)  does not require tanh activation, as the output displacement field represents physical quantities with inherently unbounded ranges that should not be artificially constrained.

\subsection{Sample point distribution strategies}
\label{sec:sample_point_distribution}

\begin{figure}[htbp]
    \centering
    \begin{minipage}{0.49\textwidth}
        \centering
        \includegraphics[width=\textwidth]{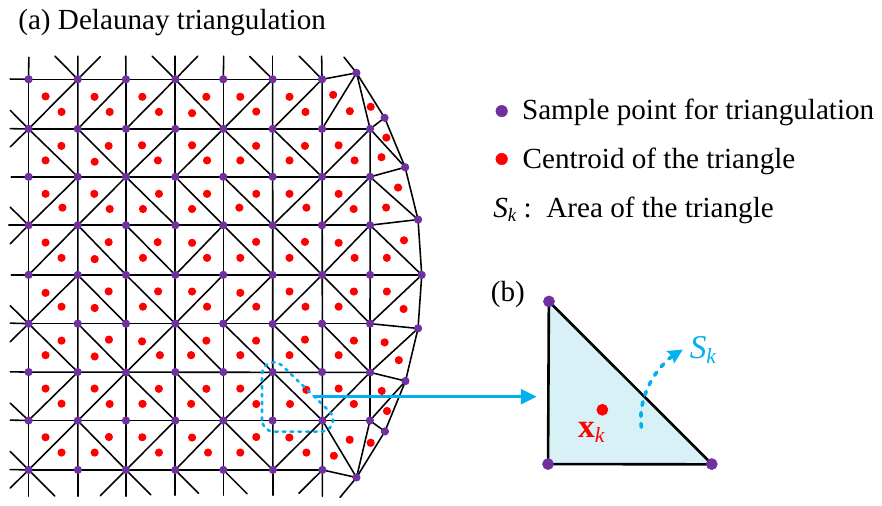}
        \caption{Triangular integration scheme: sample points and control areas.}
        \label{fig:triangular_integration}
    \end{minipage}
    \hfill
    \begin{minipage}{0.49\textwidth}
        \centering
        \includegraphics[width=0.88\textwidth]{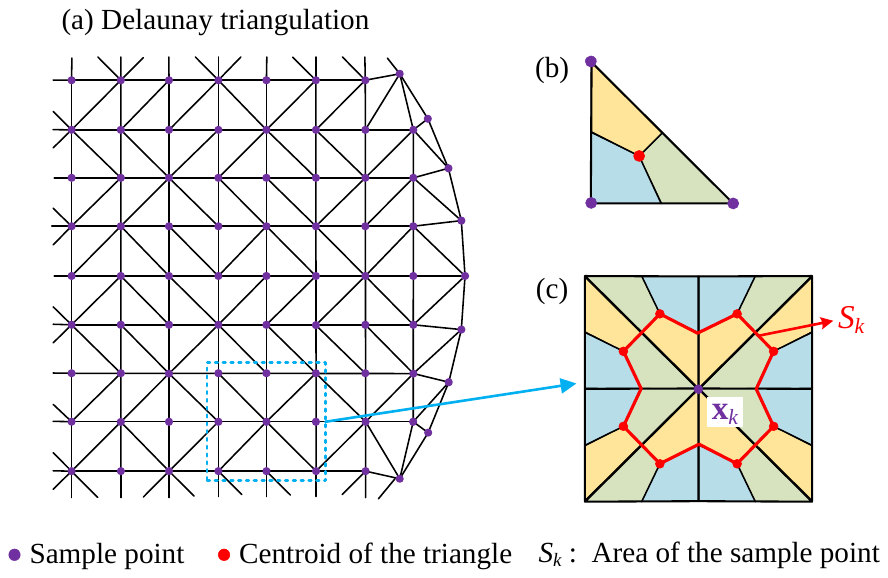}
        \caption{Delaunay integration scheme: sample points and control areas.}
        \label{fig:delaunay_integration}
    \end{minipage}
\end{figure}

We employ various numerical integration schemes to approximate the energy function, including Monte Carlo, trapezoidal, and Simpson integration. For these schemes, sample points follow uniform distribution throughout the computational domain, maintaining consistent spacing to form a regular grid. Additional boundary points can be incorporated along curved boundaries to better capture domain geometry, as illustrated by the purple points in Figs.~\ref{fig:triangular_integration} and \ref{fig:delaunay_integration}.

Additionally, we implement two triangular mesh-based integration methods for computing strain energy in the loss function. The first method employs triangular integration rules by subdividing the computational domain into triangular meshes and using triangle centroids as sample points for neural network training. Each triangle's area serves as the control area for its corresponding centroid. The triangular integration formula is
\begin{equation}
\label{eq:triangular_integration}
\int_{\Omega} F(\mathbf{x}) \mathrm{d}\Omega = \sum_{k=1}^N F(\mathbf{x}_k) \cdot S_k
\end{equation}
where $N$ is the number of triangular elements, $\mathbf{x}_k$ is the centroid of the $k$-th triangle, $S_k$ is its area, and $F$ is the integrand function.

Fig.~\ref{fig:triangular_integration} illustrates the triangular integration scheme, where purple points represent the original uniform grid, red points are triangle centroids serving as neural network training points, and the blue triangle area in Fig.~\ref{fig:triangular_integration}(b) shows the control area $S_k$ for sample point $\mathbf{x}_k$. The mesh density should be chosen based on problem complexity.

The second method, shown in Fig.~\ref{fig:delaunay_integration}, uses the same integral formula as Eq.~\eqref{eq:triangular_integration} but employs the uniform grid points (purple) as training points rather than centroids. Each triangle is subdivided by connecting its centroid to the midpoints of its edges, creating three equal areas distributed to adjacent training points. For instance, in Fig.~\ref{fig:delaunay_integration}(b), the green area is assigned to the bottom-right sample point. Each training point's control area comprises one-third of all surrounding triangles. This approach is termed Delaunay integration, where the red-outlined regions in Fig.~\ref{fig:delaunay_integration}(c) represent the control area $S_k$ for point $\mathbf{x}_k$.

The choice of sampling strategy and numerical integration scheme directly affects loss function computation. For classical numerical integration schemes based on regular sample point distributions (Monte Carlo, trapezoidal, Simpson), this study employs their standard-defined weight coefficients to compute integral approximations. For Delaunay triangulation-based numerical integration, depending on different sample points (triangle centroids or uniformly distributed points), corresponding area coefficients are selected to approximate the energy integral. Section~\ref{sec:numerical_examples} systematically compares computational results across different schemes.

\section{Multi-material elasticity problems}
\label{sec:multi_material_elasticity_problems}

Recent work has demonstrated KAN's effectiveness in solving PDEs~\cite{liu2024kan} and heterogeneous problems~\cite{WANG2025117518}.
 Based on spline interpolation theory, B-splines naturally exhibit piecewise polynomial characteristics with controllable smoothness properties. In multi-material elasticity problems, displacement fields remain continuous across material interfaces while stress fields (and consequently displacement gradients) may exhibit discontinuities due to material property changes. The piecewise nature of B-spline functions makes KAN theoretically well-suited for capturing such solution behaviors without explicit domain decomposition. 

Therefore, when solving multi-material problems, this work abandons the traditional domain decomposition approach that assigns separate neural networks to different material regions. Instead, we employ a single KAN network to predict displacement fields across the entire computational domain without requiring subdomain division. We construct a unified energy-based loss function for the complete multi-material model and optimize all trainable KAN parameters simultaneously, obtaining a neural network capable of predicting physically consistent displacement fields that naturally accommodate material property variations.

\begin{figure}[htbp]
  \centering
  \includegraphics[width=0.75\textwidth]{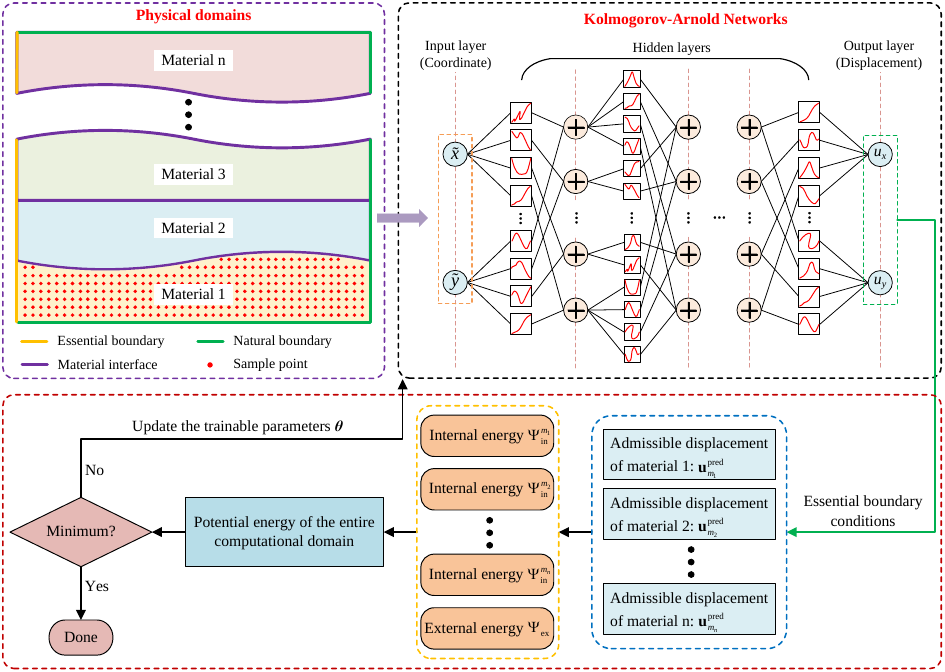}
  \caption{PIKAN flowchart for multi-material elasticity problems.}
  \label{fig:pikan_multimaterial_flowchart}
\end{figure}

Fig.~\ref{fig:pikan_multimaterial_flowchart} illustrates the PIKAN framework for solving multi-material elasticity problems. As shown in the upper left, sample point coordinates are generated throughout different material domains according to geometric features and material distribution. The training point arrangement follows the integration schemes introduced in Section 3.3, with arrangements varying by integration method. Simultaneously, boundary condition points are generated, and all normalized coordinates serve as neural network inputs. It should be noted that since we employ a single KAN network to predict the entire computational domain, no interface continuity penalty terms are required in the loss function, eliminating the need for explicit interface point generation.

A Kolmogorov-Arnold network $\mathbf{F}(\mathbf{x}; \boldsymbol{\theta})$ is constructed to predict displacement fields across all material domains simultaneously. Material properties are incorporated during energy calculation, where strain energy density is computed using the appropriate material constants for each domain based on point locations. This approach constructs a unified admissible displacement field $\mathbf{u}^{\text{pred}}(\mathbf{x})$ that naturally accommodates material discontinuities through KAN's piecewise function characteristics.

The strain energies from all material domains are summed to obtain total strain energy, and subtracting external work yields the potential energy of the entire computational domain. Using the principle of minimum potential energy, an optimizer minimizes the loss function with respect to KAN parameters $\boldsymbol{\theta}$
\begin{equation}
\label{eq:pikan_optimization_multimaterial}
\mathbf{F}(\mathbf{x}; \boldsymbol{\theta}) = \underset{\boldsymbol{\theta}}{\arg\min} \mathcal{L}_{\text{PIKAN}} = \underset{\boldsymbol{\theta}}{\arg\min} \left\{\sum_{i=1}^{n} \Psi_{\text{in}}^{m_i} - \Psi_{\text{ex}}\right\}
\end{equation}
until convergence is achieved. $\Psi_{\text{in}}^{m_i}$ represents the strain energy in the $i$-th material domain and $\Psi_{\text{ex}}$ is the total external work. The complete PIKAN solution procedure is summarized in Algorithm \ref{alg:pikan_multimaterial} in the \hyperlink{AppendixA}{Appendix A}.

\section{Numerical examples}
\label{sec:numerical_examples}

This section demonstrates the accuracy of PIKAN through several numerical examples. The first two examples explore PIKAN's application to multi-material problems, while the latter two analyze electronic packaging structures. Unless otherwise specified, all geometric dimensions in the following examples are given in millimeters (mm). To assess accuracy and convergence, we consider three error metrics: absolute error, relative error and relative $L_2$ norm error. The absolute error is defined as
\begin{equation}
\label{eq:absolute_error}
\text{Abs} = |\mathbf{u}_{\text{pred}} - \mathbf{u}_{\text{ref}}|
\end{equation}
where $\mathbf{u}_{\text{pred}}$ represents the admissible displacement solution constructed from the neural network output and $\mathbf{u}_{\text{ref}}$ is the reference solution. This metric evaluates prediction accuracy at each sample point. The relative error is computed by 
\begin{equation}
    \label{eq:relative_error}
    \text{Re} = \left\lvert \frac{\mathbf{u}_{\text{pred}} - \mathbf{u}_{\text{ref}}}{\mathbf{u}_{\text{ref}}}\right\rvert \times 100\%
\end{equation}
The relative $L_2$ norm error is defined as
\begin{equation}
\label{eq:relative_l2_error}
L_2 = \frac{\sqrt{\sum_{i=1}^{N} |u^i_{\text{pred}} - u^i_{\text{ref}} |^2}}{\sqrt{\sum_{i=1}^{N} |u^i_{\text{ref}}|^2}}
\end{equation}
where $N$ represents the total number of sample points. Since analytical solutions are unavailable for these examples, we use FEM results as reference solutions to validate the proposed method's accuracy. In the FEM setup, four-node bilinear quadrilateral elements with reduced integration (CPS4R) were employed, and mesh convergence studies were conducted for all analyses to ensure the validity of the reference solutions.

All examples employ a consistent computational framework. We use the L-BFGS optimizer with initial learning rate 0.001  and constrained to a maximum of 20 iterations per optimization step. The KAN network architectural hyperparameters and optimizer settings employed in this study primarily follow the experience ranges widely adopted in the original KAN paper and recent physics-informed learning applications \cite{liu2024kan,WANG2025117518}. These were manually fine-tuned through iterative adjustments to achieve accurate solution predictions. Complete hyperparameter specifications are provided for each example to ensure reproducibility. The total trainable parameter counts reported for all examples in this section refer exclusively to the KAN network parameters as defined in Table~\ref{tab:kan_hyperparameters}, and do not include any parameters from boundary condition treatment functions. Neural network training was performed on an NVIDIA GeForce RTX 4060 Laptop GPU.

\subsection{Cantilever beam analysis across material interfaces}
\label{sec:cantilever_beam_analysis}

\begin{figure}[htbp]
  \centering
  \begin{subfigure}[b]{0.48\textwidth}
      \centering
      \includegraphics[width=0.9\textwidth]{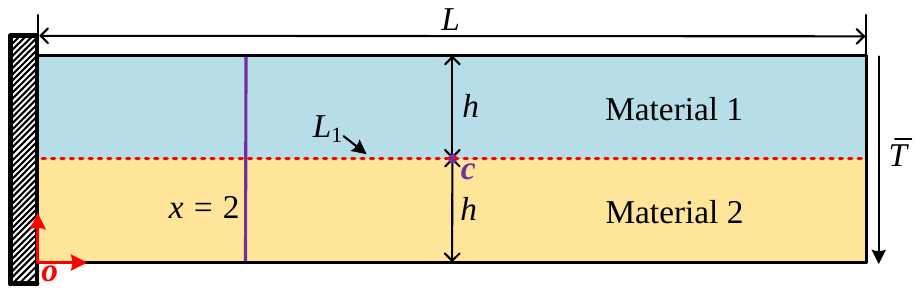}
      \caption{Straight interface}
      \label{fig:cantilever_straight}
  \end{subfigure}
  %\hfill
  \begin{subfigure}[b]{0.48\textwidth}
      \centering
      \includegraphics[width=0.9\textwidth]{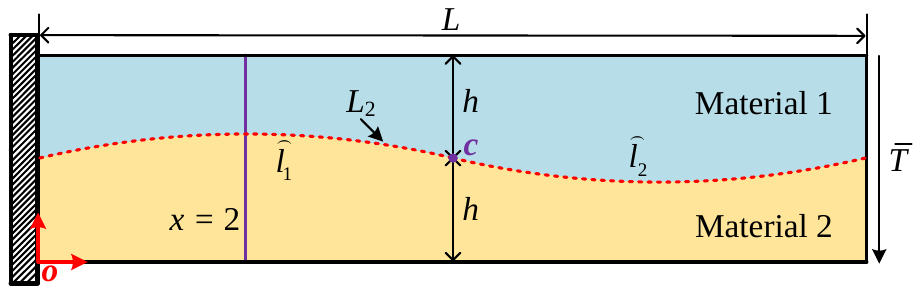}
      \caption{Wavy interface}
      \label{fig:cantilever_wavy}
  \end{subfigure}
  %\hfill
  \begin{subfigure}[b]{0.48\textwidth}
      \centering
      \includegraphics[width=0.9\textwidth]{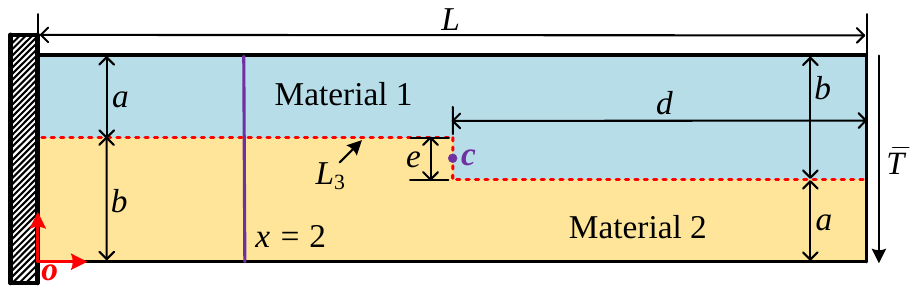}
      \caption{Stepped interface}
      \label{fig:cantilever_stepped}
  \end{subfigure}
   %\hfill
  \begin{subfigure}[b]{0.48\textwidth}
      \centering
      \includegraphics[width=0.9\textwidth]{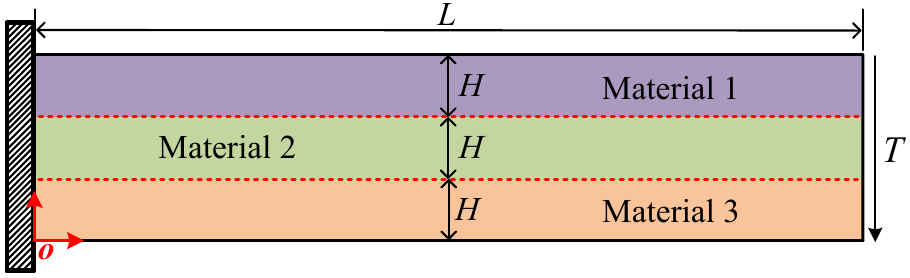}
      \caption{Three-layer composite cantilever beam }
      \label{fig:three layer beam}
  \end{subfigure}
  \caption{Cantilever beam models with different material interface geometries.}
  \label{fig:cantilever_beam_models}
\end{figure}

We investigate PIKAN's performance through three cantilever beam bending problems with different material interface geometries, as shown in Fig.~\ref{fig:cantilever_beam_models}. Each beam has length $L = 8$ and height $2h = 2$, with the coordinate origin at the bottom-left corner. The straight interface (Fig.~\ref{fig:cantilever_straight}) is located at distance $h = 1$ from the bottom edge. The wavy interface (Fig.~\ref{fig:cantilever_wavy}) consists of two arc segments $\overset{\frown}{l}_1$ and $\overset{\frown}{l}_2$ that are antisymmetric about the beam's center point c $(4,1)$. Arc $\overset{\frown}{l}_1$ belongs to a circle centered at $(2,-8)$ and arc $\overset{\frown}{l}_2$ to a circle centered at $(6,10)$, both with radius $r = \sqrt{85}$. The stepped interface (Fig.~\ref{fig:cantilever_stepped}) has geometric parameters $a = 1.2$, $b = 0.8$, $d = 4$, and $e = 0.4$. As shown in Fig.~\ref{fig:cantilever_beam_models}, all models are subjected to a shear traction $\bar{T} = 6$ N/mm applied at the right boundary. The left end is fully constrained. 

The material parameters are: Material 1 with Young's modulus $E_1 = 8,500$ MPa and Material 2 with $E_2 = 43,000$ MPa, both having Poisson's ratio $\nu = 0.3$.
 The method uses a KAN structure of $[2,5,5,5,2]$ with grid size 10, B-spline order 3, and grid range $[0,1]$, totaling 1050 trainable parameters.
The material interfaces are depicted as red dashed lines $L_1$, $L_2$, and $L_3$ in Fig.~\ref{fig:cantilever_beam_models}, representing straight, wavy, and stepped configurations for the three cantilever beam models.

% The material interfaces are depicted as red dashed lines and mathematically defined as $L_1$, $L_2$, and $L_3$ for the straight, wavy, and stepped configurations, respectively.
% \begin{align}
%     \small
%     \label{eq:interface_definitions}
%     L_1 &= \{(x, y): 0 \leq x \leq L, y = h\} \\
%     L_2 &= \{(x, y): (x-2)^2+(y+8)^2 = r^2 \text{ for } 0 \leq x \leq L/2, \nonumber \\
%         &\phantom{= \{(x, y):}\,(x-6)^2+(y-10)^2 = r^2 \text{ for } L/2 \leq x \leq L\} \\
%     L_3 &= \{(x, y): \begin{cases}
%            0 \leq x \leq L/2, & y = b \\
%            x = L/2, & a \leq y \leq b \\
%            L/2 \leq x \leq L, & y = a
%            \end{cases}\}
%     \end{align}

% \begin{table}[htbp]
%   \scriptsize
%   \centering
%   \caption{Material parameters used in cantilever beam models}
%   \label{tab:material_parameters}
%   \begin{tabular}{ccc}
%   \toprule
%   Material & Young's Modulus $E$ (MPa) & Poisson's Ratio $\nu$ \\
%   \midrule
%   1 & 8500 & 0.3 \\
%   2 & 43000 & 0.3 \\
%   \bottomrule
%   \end{tabular}
%   \end{table}

\begin{figure}[htbp]
  \centering
  % ---------- 第一行 ----------
  \begin{subfigure}[b]{0.3\linewidth}
    \centering
    \includegraphics[width=\linewidth]{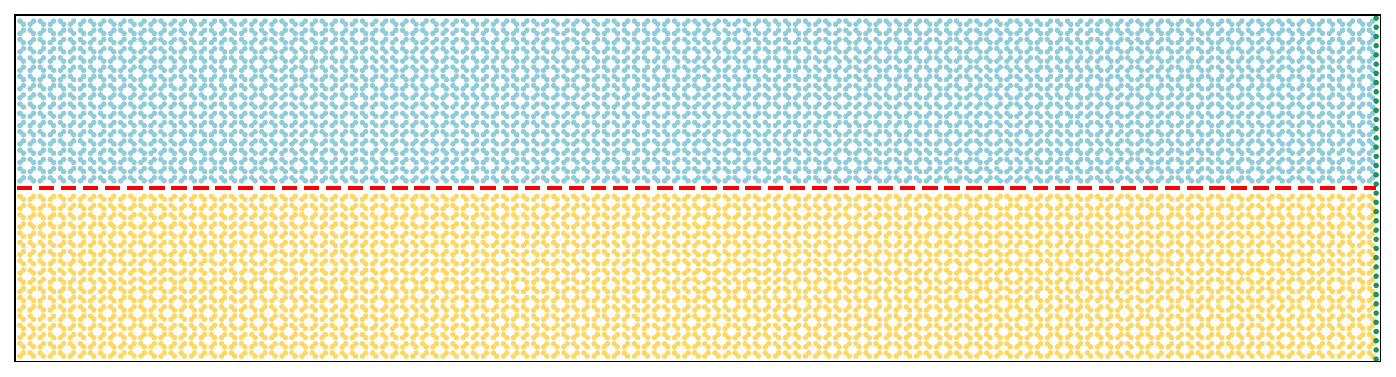}
    \caption{Straight interface - triangular}
    \label{fig:training_points_straight_triangular}
  \end{subfigure}\hspace{0.8em}
  \begin{subfigure}[b]{0.3\linewidth}
    \centering
    \includegraphics[width=\linewidth]{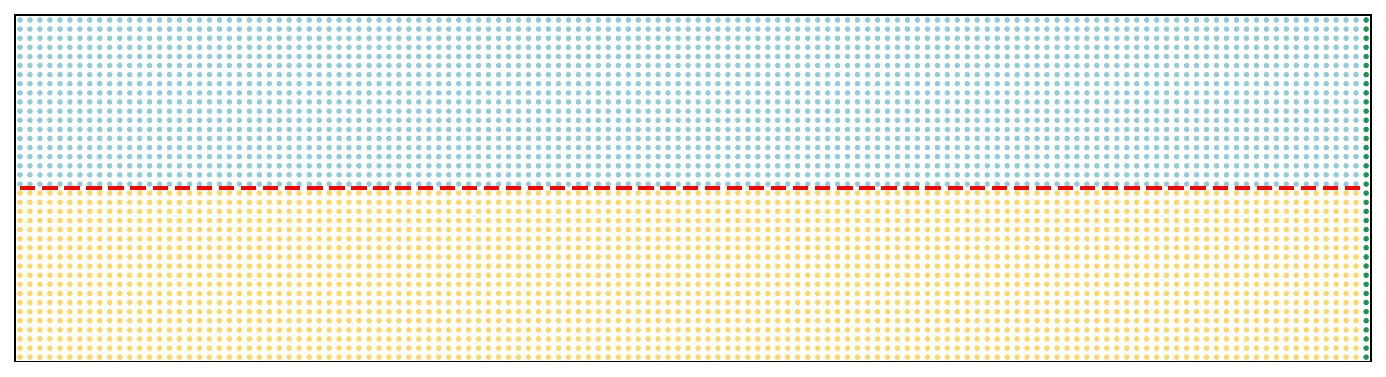}
    \caption{Straight interface - uniform}
    \label{fig:training_points_straight_uniform}
  \end{subfigure}\hspace{0.8em}
  \begin{subfigure}[b]{0.3\linewidth}
    \centering
    \includegraphics[width=\linewidth]{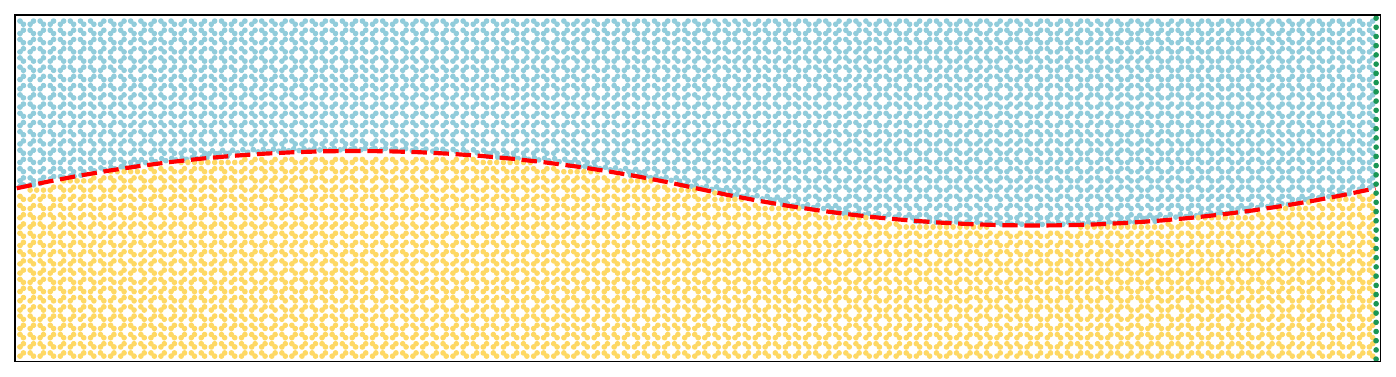}
    \caption{Wavy interface - triangular}
    \label{fig:training_points_wavy_triangular}
  \end{subfigure}

  \medskip
  % ---------- 第二行 ----------
  \begin{subfigure}[b]{0.3\linewidth}
    \centering
    \includegraphics[width=\linewidth]{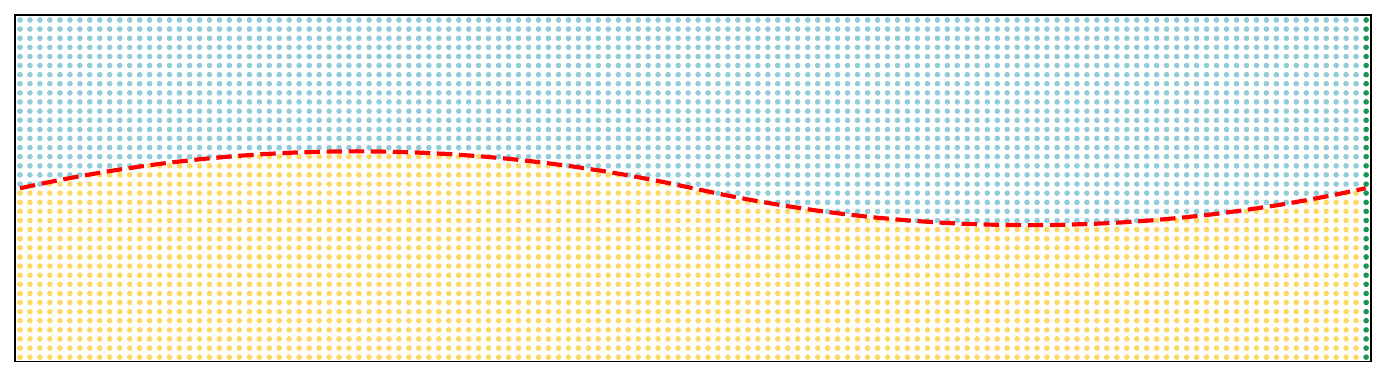}
    \caption{Wavy interface - uniform}
    \label{fig:training_points_wavy_uniform}
  \end{subfigure}\hspace{0.8em}
  \begin{subfigure}[b]{0.3\linewidth}
    \centering
    \includegraphics[width=\linewidth]{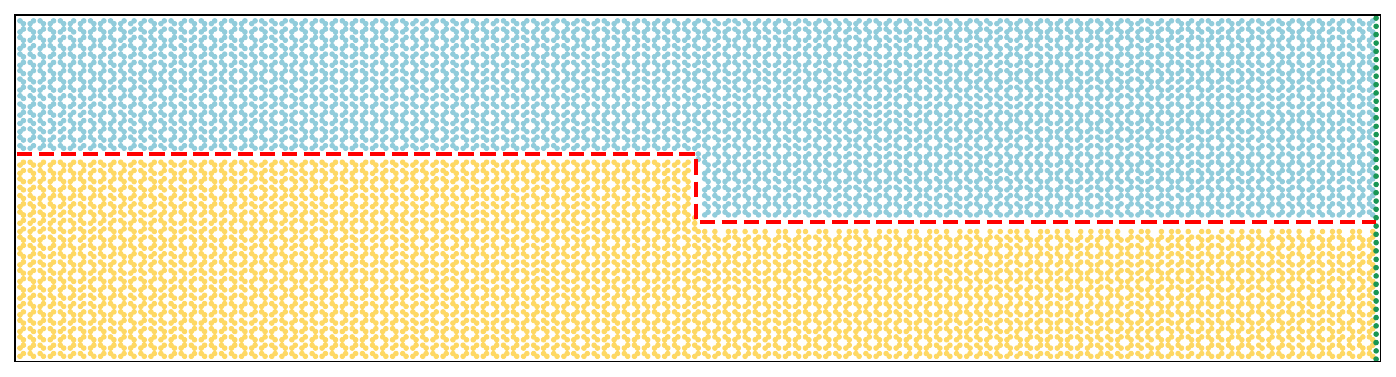}
    \caption{Stepped interface - triangular}
    \label{fig:training_points_stepped_triangular}
  \end{subfigure}\hspace{0.8em}
  \begin{subfigure}[b]{0.3\linewidth}
    \centering
    \includegraphics[width=\linewidth]{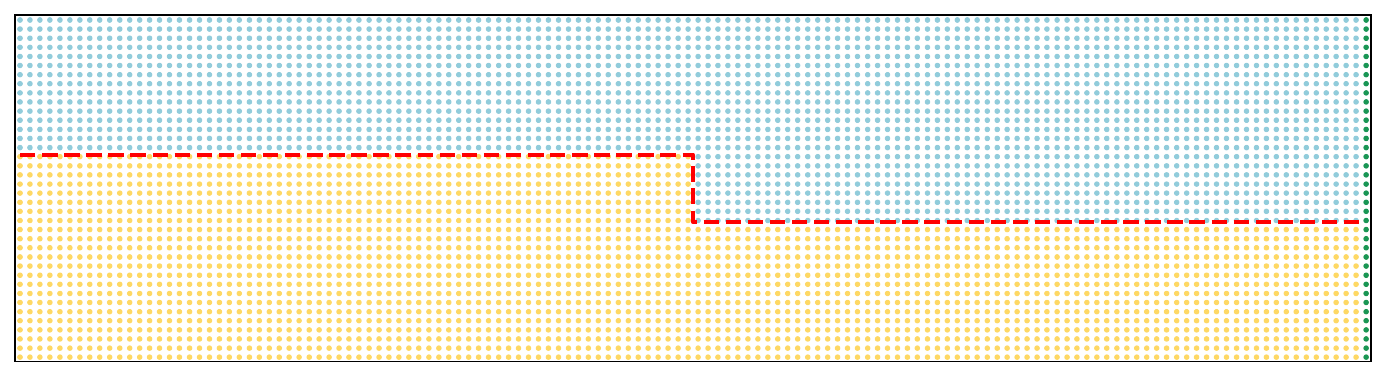}
    \caption{Stepped interface - uniform}
    \label{fig:training_points_stepped_uniform}
  \end{subfigure}

  \caption{Training point distributions: (a,c,e) triangular scheme; (b,d,f) uniform scheme.}
  \label{fig:training_points}
\end{figure}

Fig.~\ref{fig:training_points} illustrates the distribution of training points used in the analysis of different models. Green points on the right boundary represent natural boundary condition points, with 101 points uniformly distributed at spacing $dy = 0.02$ (this boundary point distribution is consistent across all cantilever beam models). For the centroid sampling strategy, the interior point counts are: 79,200 points for the straight interface beam (Fig.~\ref{fig:training_points_straight_triangular}), 80,278 points for the wavy interface beam (Fig.~\ref{fig:training_points_wavy_triangular}), and 79,163 points for the stepped interface beam (Fig.~\ref{fig:training_points_stepped_triangular}). In contrast, the uniform sampling shown in Figs.~\ref{fig:training_points_straight_uniform}, \ref{fig:training_points_wavy_uniform}, and \ref{fig:training_points_stepped_uniform} employ identical distribution patterns with 40,501 interior points for each beam configuration.

\begin{figure}[htbp]
    \centering
    \begin{subfigure}[b]{0.32\textwidth}
        \centering
        \includegraphics[width=\textwidth]{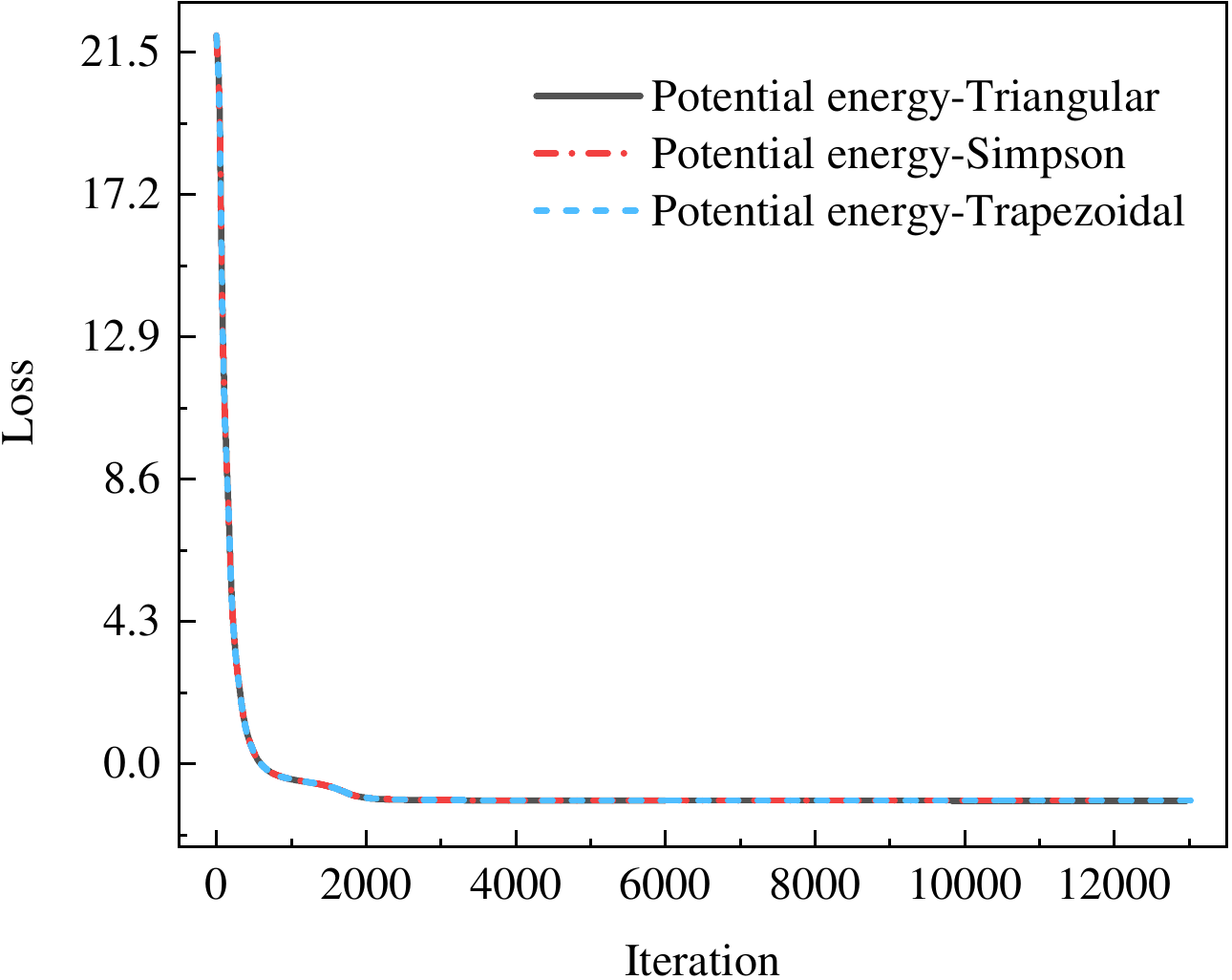}
        \caption{Straight interface}
        \label{fig:loss_straight}
    \end{subfigure}
    \hfill
    \begin{subfigure}[b]{0.32\textwidth}
        \centering
        \includegraphics[width=1.\textwidth]{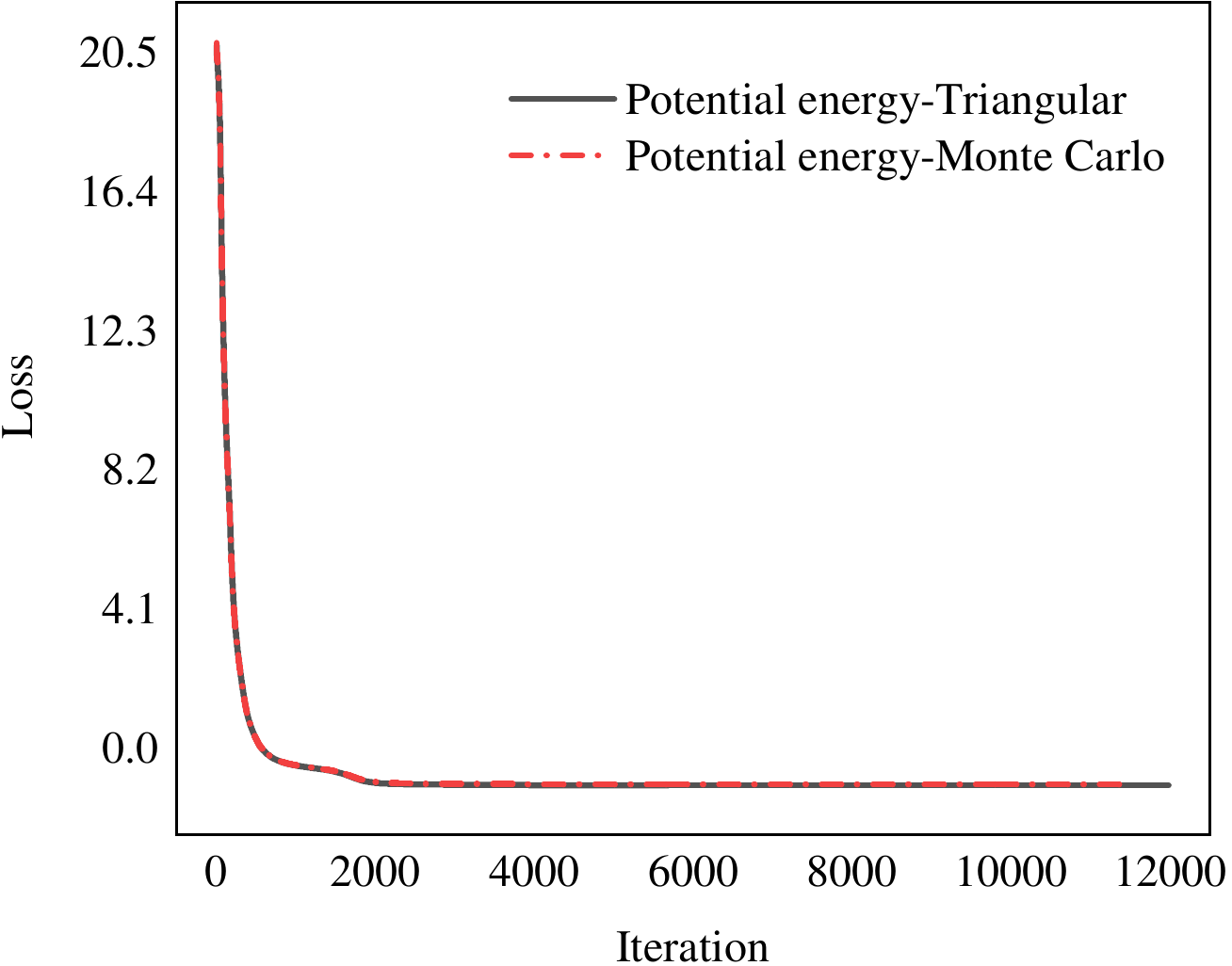}
        \caption{Wavy interface}
        \label{fig:loss_wavy}
    \end{subfigure}
    \hfill
    \begin{subfigure}[b]{0.32\textwidth}
        \centering
        \includegraphics[width=0.97\textwidth]{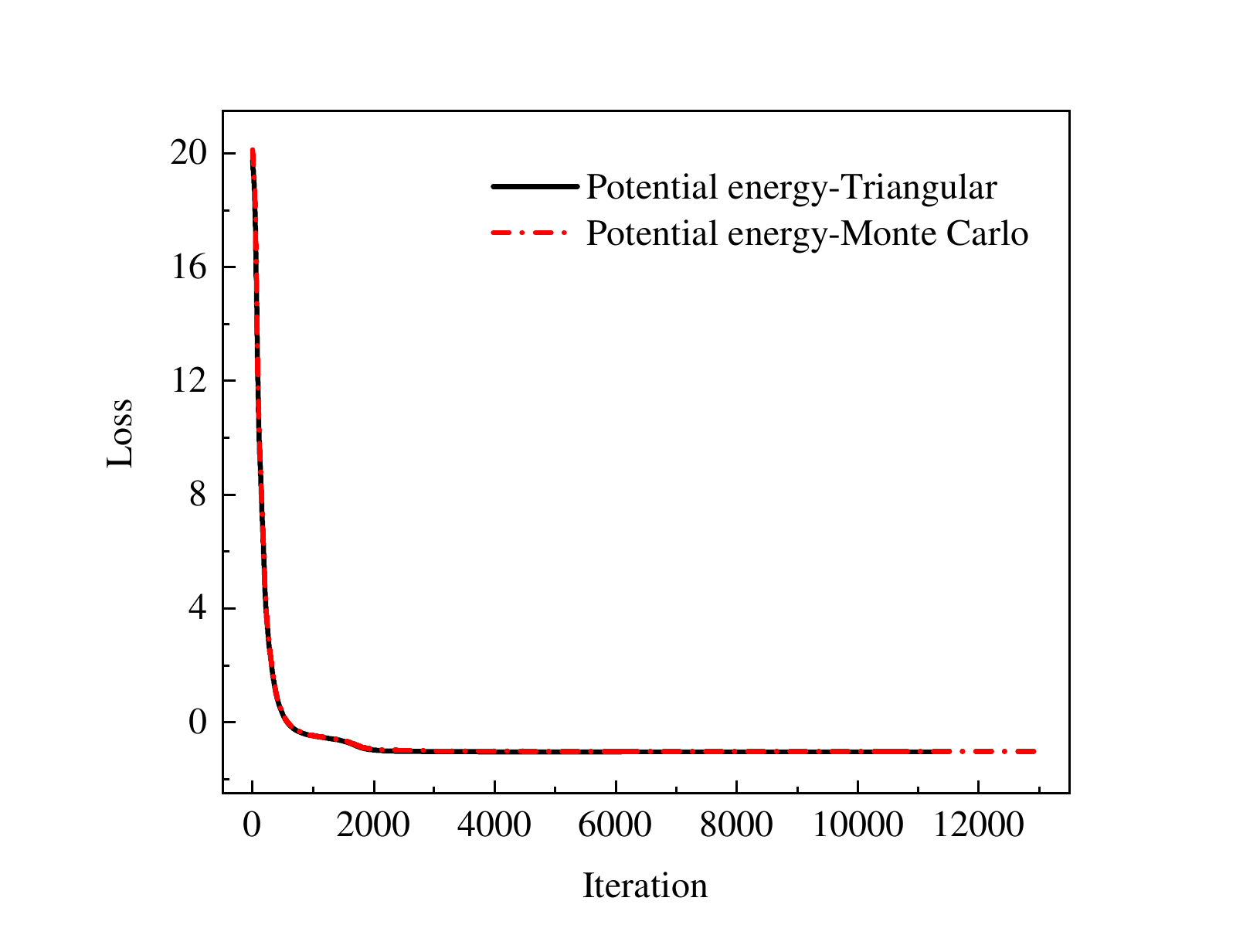}
        \caption{Stepped interface}
        \label{fig:loss_stepped}
    \end{subfigure}
    
    \caption{Evolution of energy-based loss functions for different interface geometries.}
    \label{fig:loss_evolution}
\end{figure}

Fig.~\ref{fig:loss_evolution} shows the evolution of loss functions during training for different numerical integration methods: triangular, Simpson, trapezoidal and Monte Carlo schemes. The triangular scheme uses the point distribution from Fig.~\ref{fig:training_points}a, c, e, while Simpson, trapezoidal and Monte Carlo schemes use the uniform distribution from Fig.~\ref{fig:training_points}b, d, f. All schemes exhibit convergent behavior, with the triangular integration method generally achieving faster convergence and lower final energy values. Despite different interface geometries, all models demonstrate similar convergence trends with increasing training iterations, indicating the robustness of the PIKAN approach across various geometric configurations.

\begin{figure}[H]
  \centering
  
  \begin{subfigure}[b]{0.85\textwidth}
      \centering
      \includegraphics[width=\textwidth]{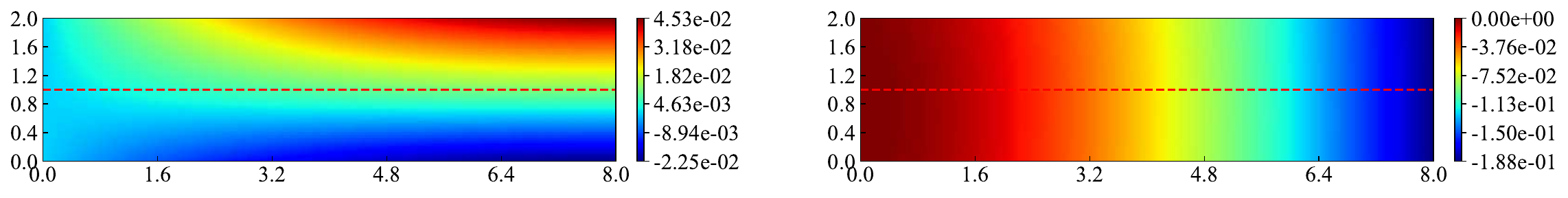}
      \caption{Straight interface - PIKAN}
      \label{fig:disp_straight_pikan}
  \end{subfigure}
  
  \begin{subfigure}[b]{0.85\textwidth}
      \centering
      \includegraphics[width=\textwidth]{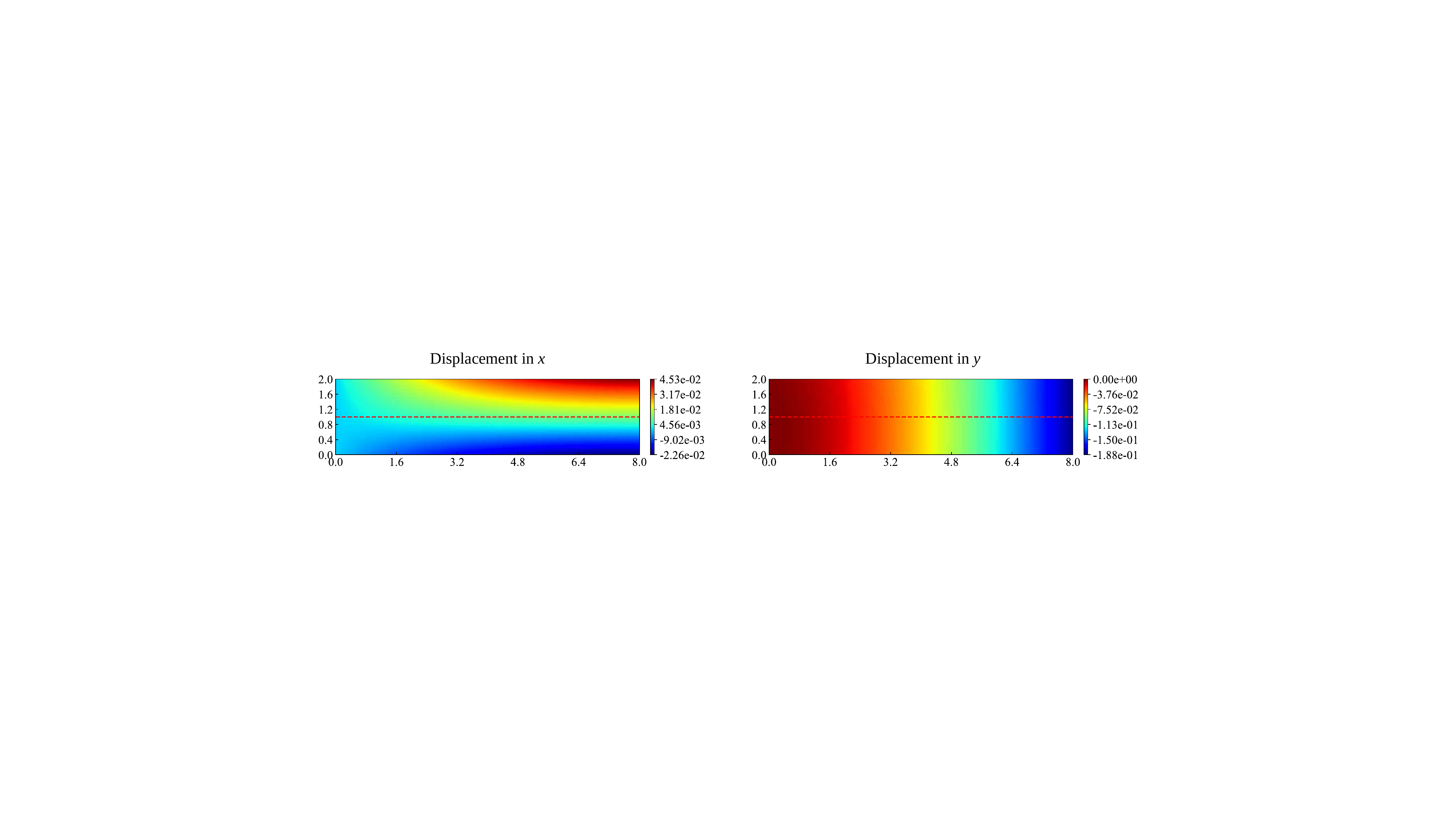}
      \caption{Straight interface - FEM}
      \label{fig:disp_straight_fem}
  \end{subfigure}
  
  \begin{subfigure}[b]{0.85\textwidth}
      \centering
      \includegraphics[width=\textwidth]{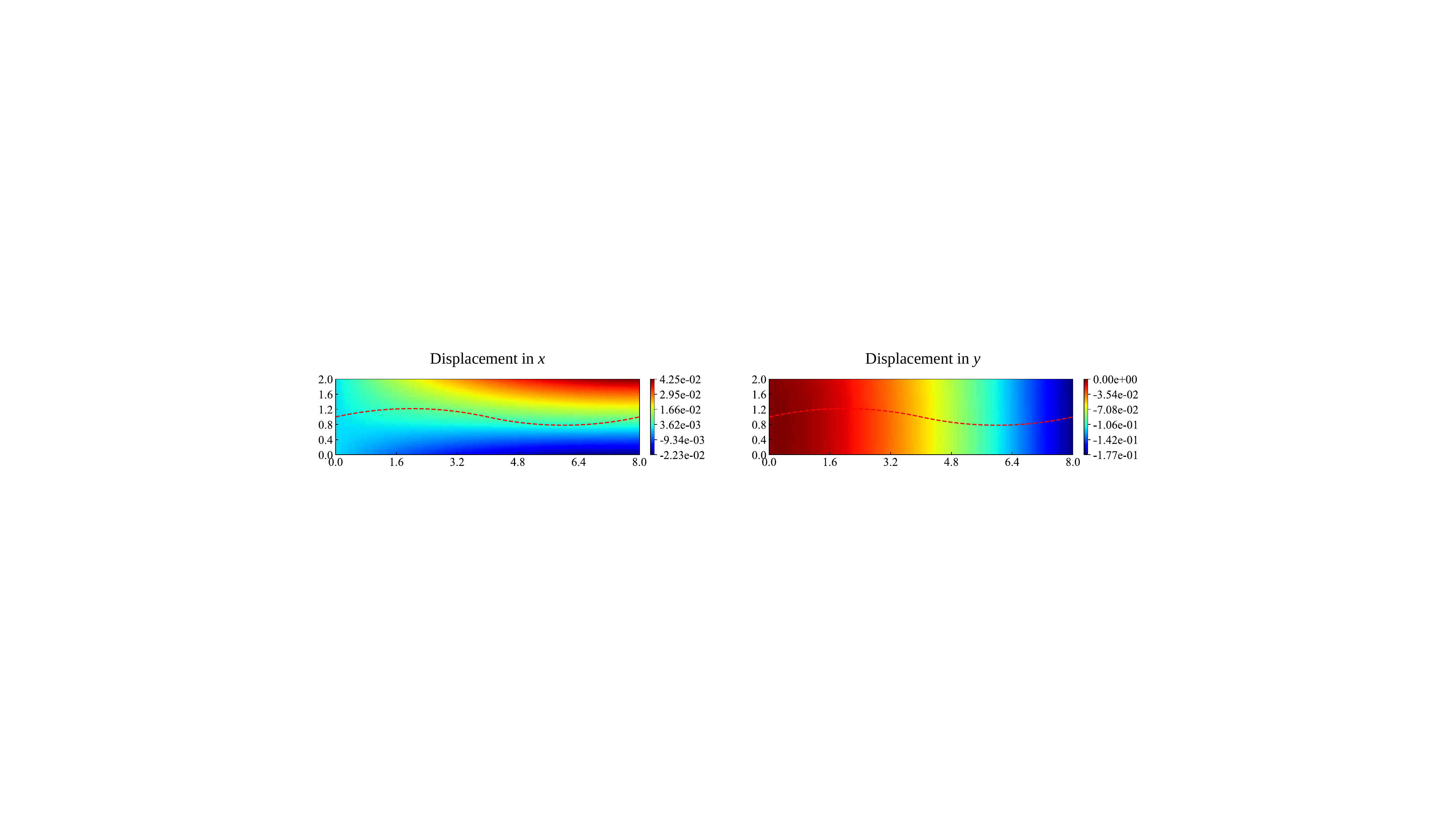}
      \caption{Wavy interface - PIKAN}
      \label{fig:disp_wavy_pikan}
  \end{subfigure}
  
  \begin{subfigure}[b]{0.85\textwidth}
      \centering
      \includegraphics[width=\textwidth]{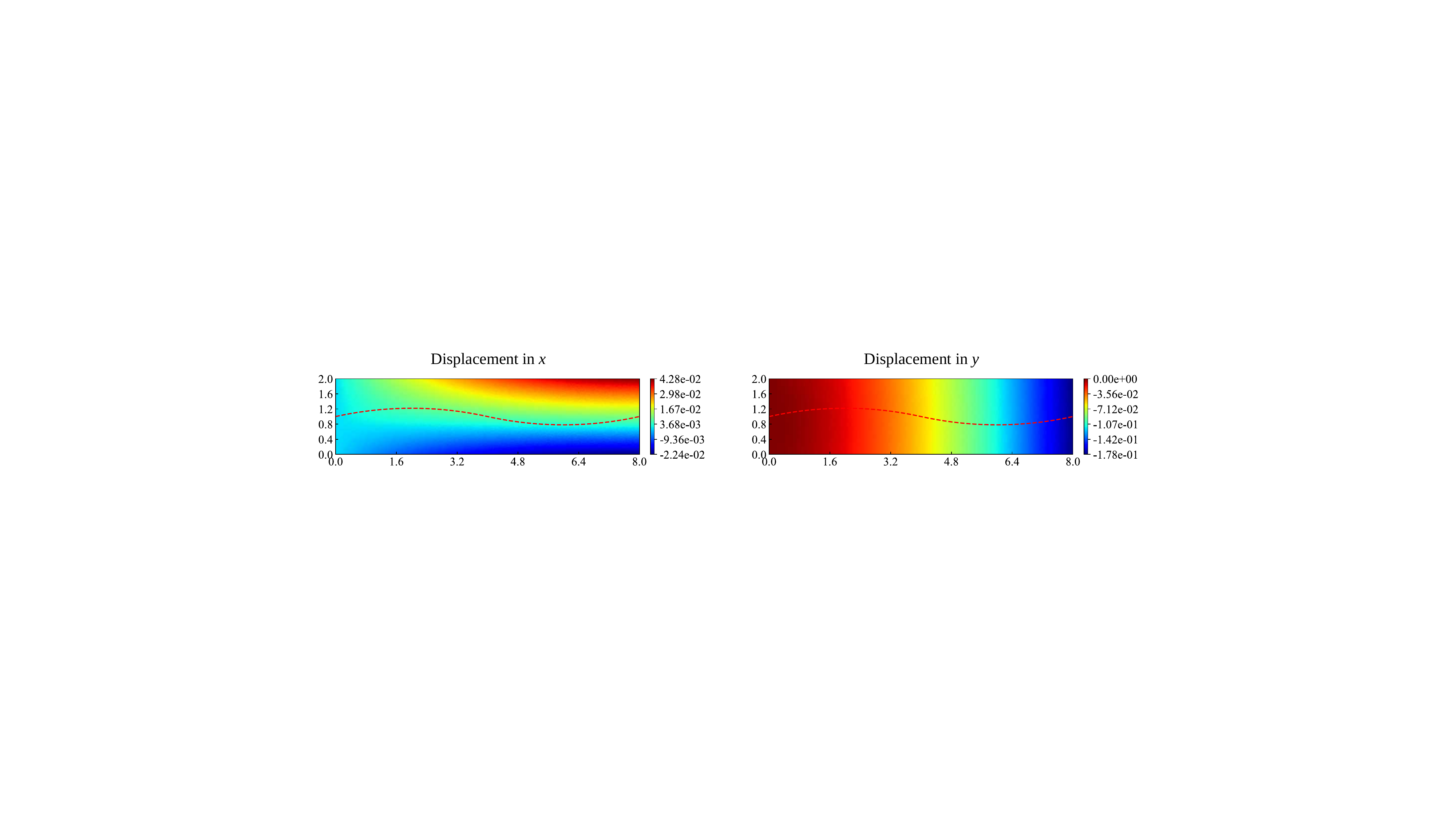}
      \caption{Wavy interface - FEM}
      \label{fig:disp_wavy_fem}
  \end{subfigure}
  
  \begin{subfigure}[b]{0.85\textwidth}
      \centering
      \includegraphics[width=\textwidth]{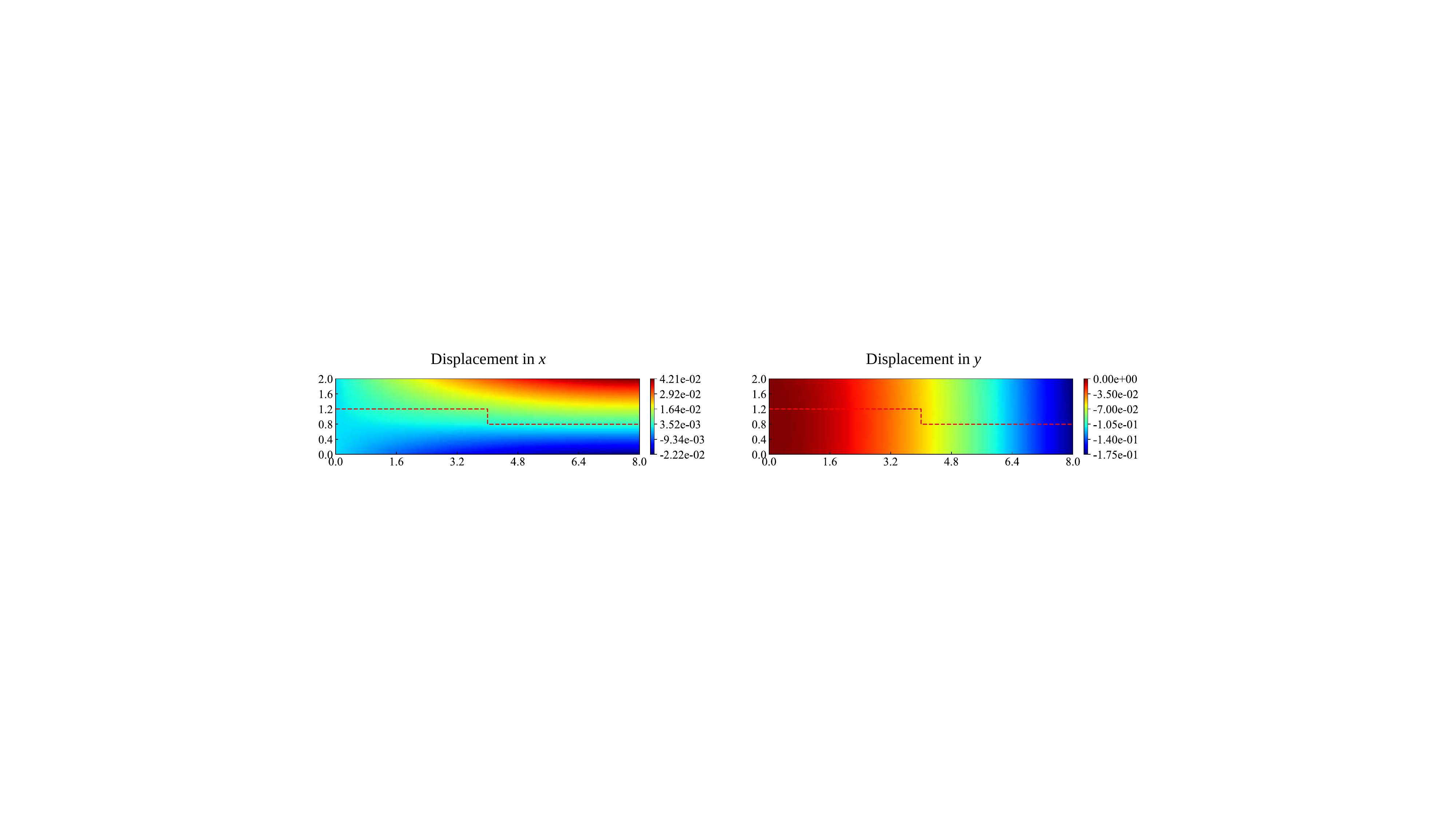}
      \caption{Stepped interface - PIKAN}
      \label{fig:disp_stepped_pikan}
  \end{subfigure}
  
  \begin{subfigure}[b]{0.85\textwidth}
      \centering
      \includegraphics[width=\textwidth]{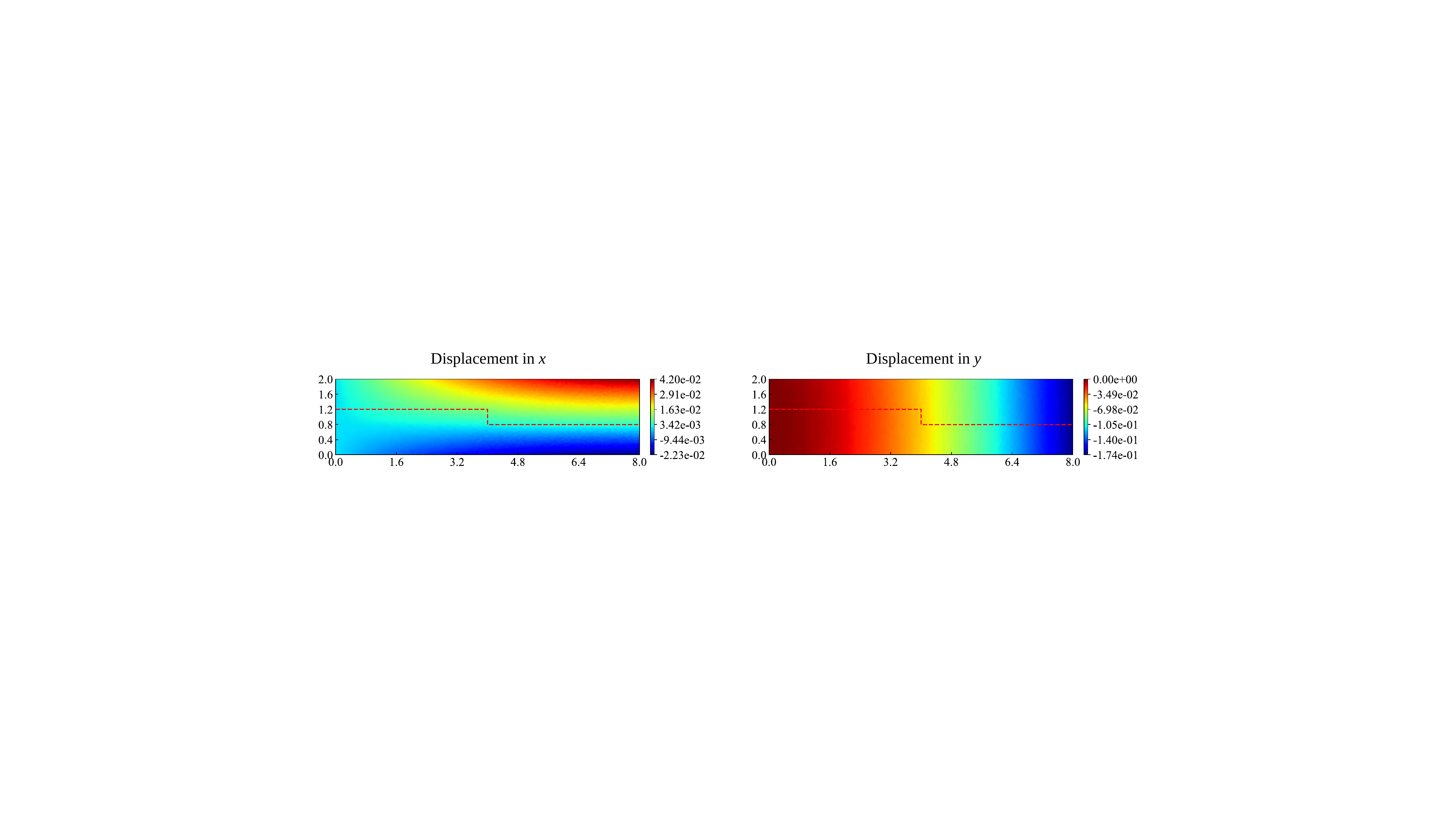}
      \caption{Stepped interface - FEM}
      \label{fig:disp_stepped_fem}
  \end{subfigure}
  
  \caption{Displacement fields $u_x$ and $u_y$ for different interface geometries: (a), (c), (e) PIKAN predictions; (b), (d), (f) FEM solutions.}
  \label{fig:displacement_comparison_cantilever}
\end{figure}

\begin{figure}[htbp]
  \centering
  
  \begin{subfigure}[b]{0.85\textwidth}
      \centering
      \includegraphics[width=\textwidth]{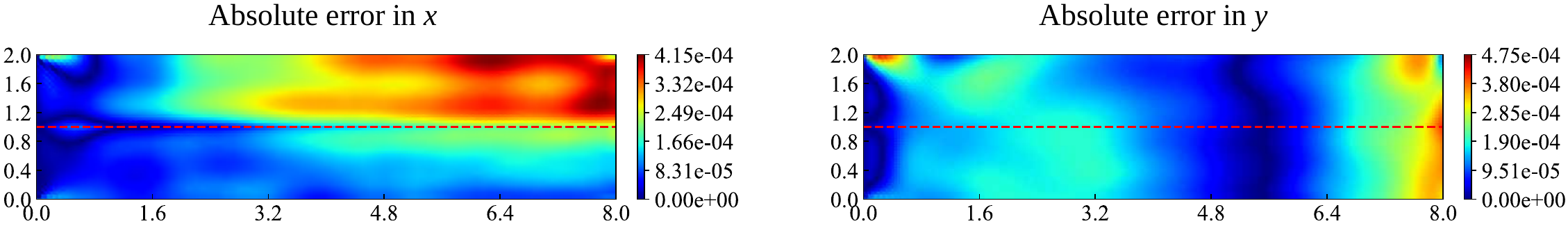}
      \caption{Straight interface}
      \label{fig:absolute_error_line_surface}
  \end{subfigure}
  
  \begin{subfigure}[b]{0.85\textwidth}
      \centering
      \includegraphics[width=\textwidth]{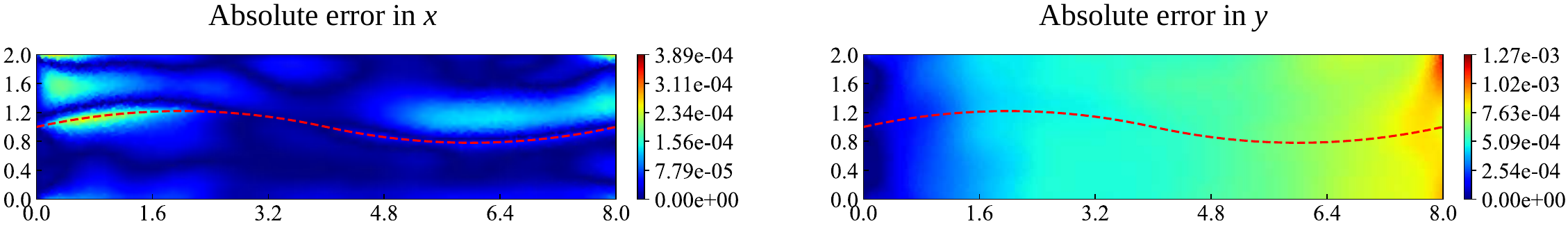}
      \caption{Wavy interface}
      \label{fig:absolute_error_s_surface}
  \end{subfigure}
  
  \begin{subfigure}[b]{0.85\textwidth}
      \centering
      \includegraphics[width=\textwidth]{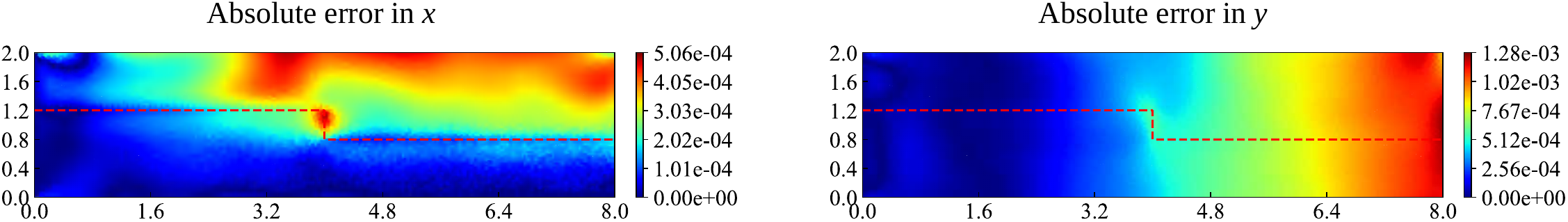}
      \caption{Stepped interface}
      \label{fig:absolute_error_z_surface}
  \end{subfigure}
  
  \caption{Absolute error contours between PIKAN and FEM solutions for $u_x$ and $u_y$.}
  \label{fig:Absolute_error_comparison_cantilever}
\end{figure}

Fig.~\ref{fig:displacement_comparison_cantilever} presents a detailed comparison of displacement fields between PIKAN results (using triangular integration) and FEM reference solutions. To provide a quantitative and intuitive assessment of PIKAN’s prediction accuracy, Fig.~\ref{fig:Absolute_error_comparison_cantilever} presents pointwise absolute error plots for cantilever beam models with different material interface configurations. PIKAN accurately reproduces both $u_x$ and $u_y$ displacement components across all interface configurations, demonstrating excellent agreement with FEM solutions. The method successfully captures complex displacement near material interfaces, particularly evident in the wavy and stepped configurations where interface geometry creates significant stress concentrations and discontinuities in material properties. 

To further validate PIKAN's computational accuracy, we compare results with the Conservative Energy Neural Network (CENN)~\cite{wang2022cenn} method using domain decomposition. CENN applies traditional DEM with separate MLPs for different subdomains, requiring interface continuity conditions that introduce additional penalty terms and hyperparameters in the loss function. The CENN loss function takes the form
\begin{equation}
\label{eq:cenn_loss}
\mathcal{L}_{\text{CENN}} = \sum_{i=1}^{n} \Psi_{\text{in}}^{i} - \Psi_{\text{ex}} + \beta \sum_{j=1}^{N_{\text{inter}}} \left\|\mathbf{u}_{\text{inter}}^{+}(\mathbf{x}_{\text{inter},j}) - \mathbf{u}_{\text{inter}}^{-}(\mathbf{x}_{\text{inter},j})\right\|^2
\end{equation}

Compared to single material domain energy-based loss functions, Eq.~\eqref{eq:cenn_loss} includes additional interface continuity penalty terms, where $\mathbf{u}_{\text{inter}}^{+}$ and $\mathbf{u}_{\text{inter}}^{-}$ represent the admissible displacement fields from different material subnetworks at the interface, and $\beta$ is the penalty parameter controlling interface condition enforcement. The penalty parameter is calculated using
\begin{equation}
\label{eq:beta_calculation}
\beta = -\lambda \cdot \ln\left[\tanh\left(N_{\text{inter}} / N_{\text{dom}}\right)\right]
\end{equation}
where $\lambda = 1000$, $N_{\text{inter}}$ is the number of interface sample points, and $N_{\text{dom}}$ is the total number of domain sample points. CENN employs identical MLP architectures for each material domain, each consisting of 5 hidden layers with 20 neurons per layer, yielding a total of 3,564 trainable parameters. The configuration uses tanh activation functions and the same triangular integration scheme for training sample point generation.

\begin{figure}[htbp]
    \centering
    
    \begin{subfigure}[b]{0.32\textwidth}
        \centering
        \includegraphics[width=\textwidth]{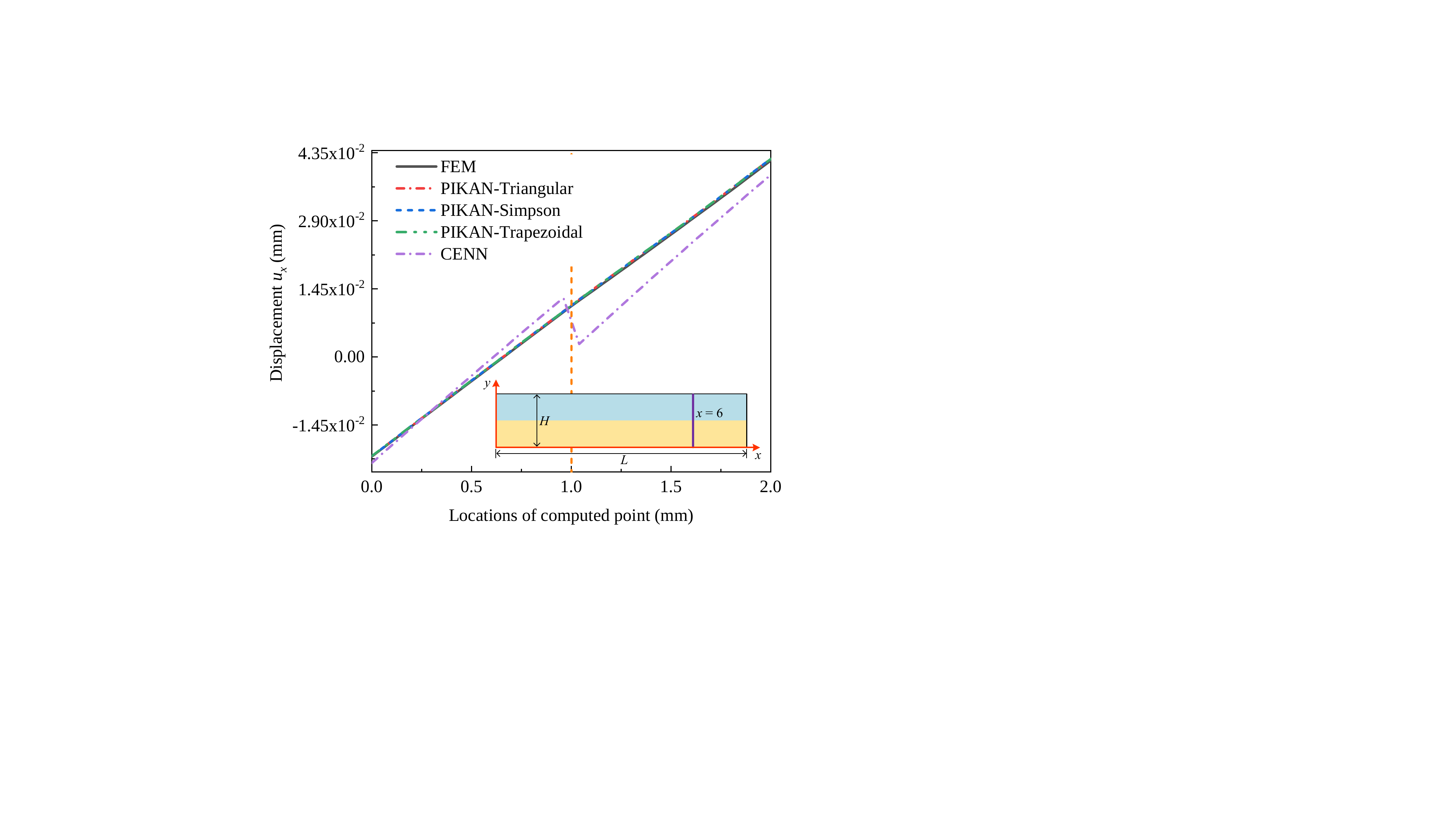}
        \caption{Straight interface}
        \label{fig:disp_x6_straight}
    \end{subfigure}
    \hfill
    \begin{subfigure}[b]{0.32\textwidth}
        \centering
        \includegraphics[width=\textwidth]{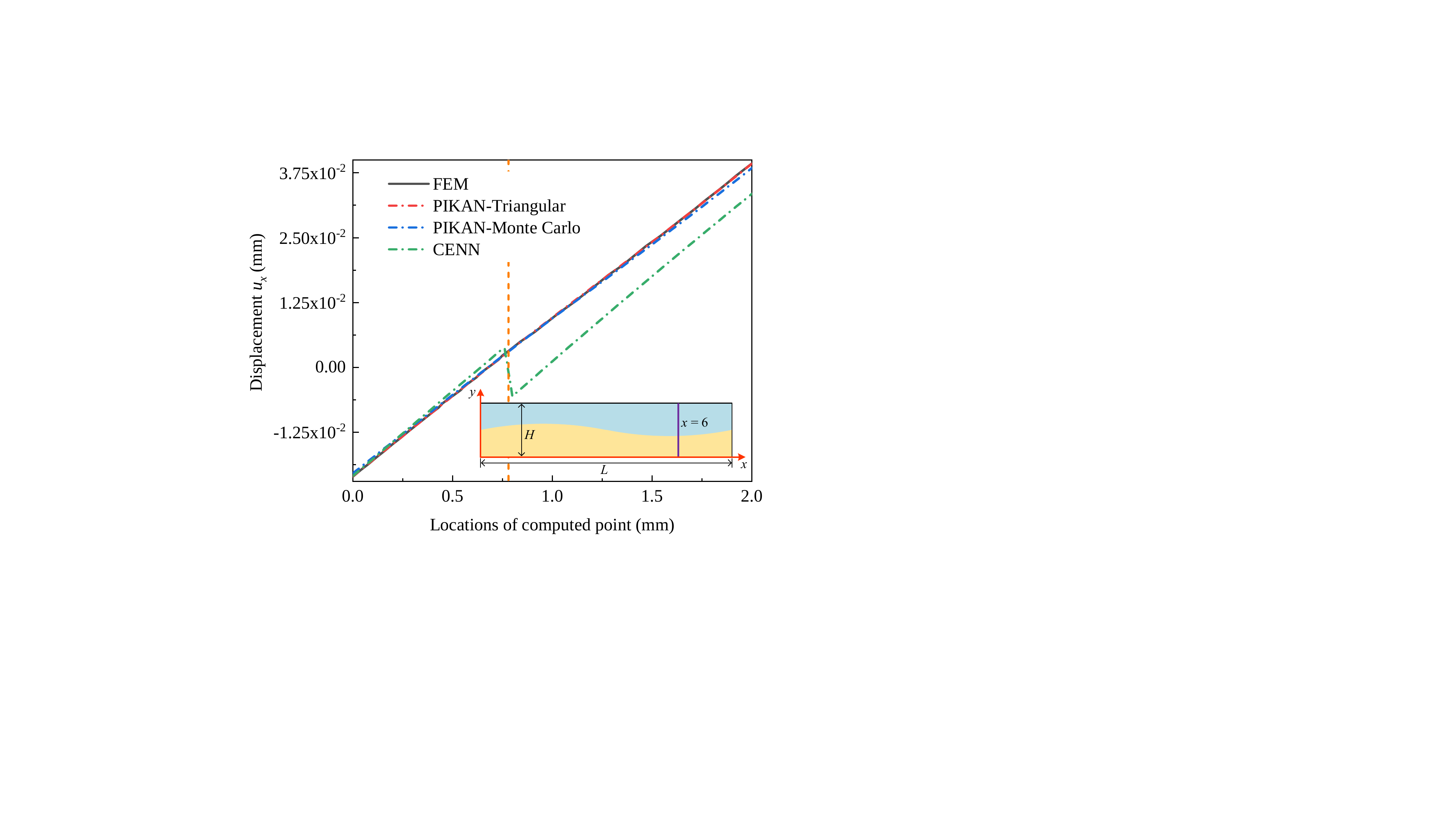}
        \caption{Wavy interface}
        \label{fig:disp_x6_wavy}
    \end{subfigure}
    \hfill
    \begin{subfigure}[b]{0.32\textwidth}
        \centering
        \includegraphics[width=\textwidth]{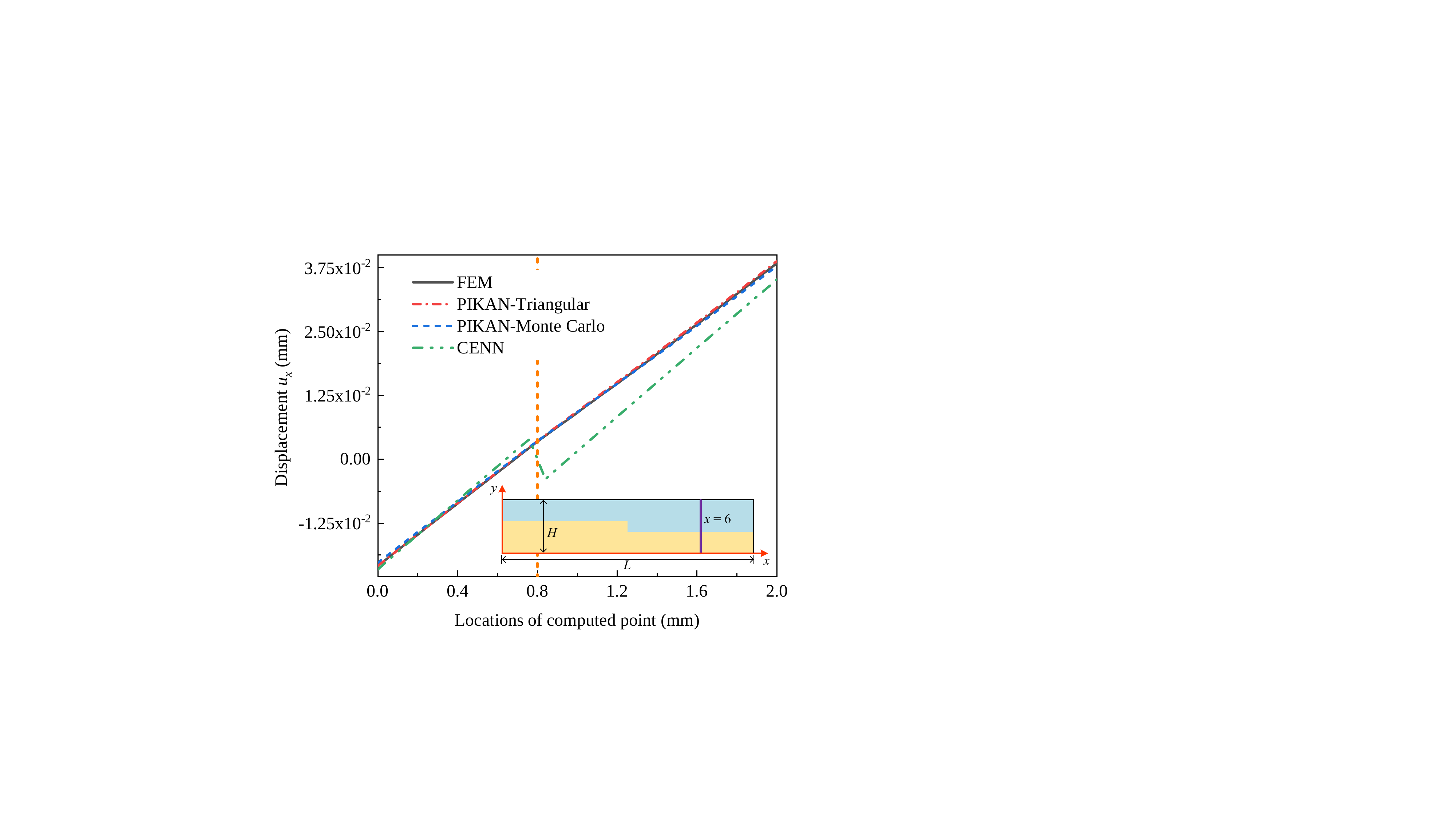}
        \caption{Stepped interface}
        \label{fig:disp_x6_stepped}
    \end{subfigure}
    
   \caption{Displacement ($u_x$) at $x = 6$ for different interface geometries: PIKAN vs. CENN and FEM solutions. The orange dashed line denotes the interface position.}
    \label{fig:displacement_comparison_x6}
\end{figure}

Fig.~\ref{fig:displacement_comparison_x6} compares $u_x$ displacement profiles along $x = 6$ for all three cantilever beam models. PIKAN results using different integration schemes are compared with CENN predictions and FEM reference solutions. The figures show PIKAN's superior performance across all interface geometries, closely matching FEM solutions with excellent accuracy and stability. Among the integration schemes, triangular, Simpson and trapezoidal methods show the best agreement with FEM solutions, while Monte Carlo integration exhibits slightly reduced precision due to its stochastic nature. 
In contrast, CENN produces noticeable displacement discontinuities at material interfaces due to its subdomain decomposition approach. Under complex interface geometries (particularly wavy and stepped configurations), it shows computational instability with significant localized errors, highlighting the inherent challenges of enforcing interface continuity conditions through penalty parameter tuning in domain decomposition methods.

\begin{figure}[htbp]
    \centering
    \begin{subfigure}[b]{0.32\textwidth}
        \centering
        \includegraphics[width=\textwidth]{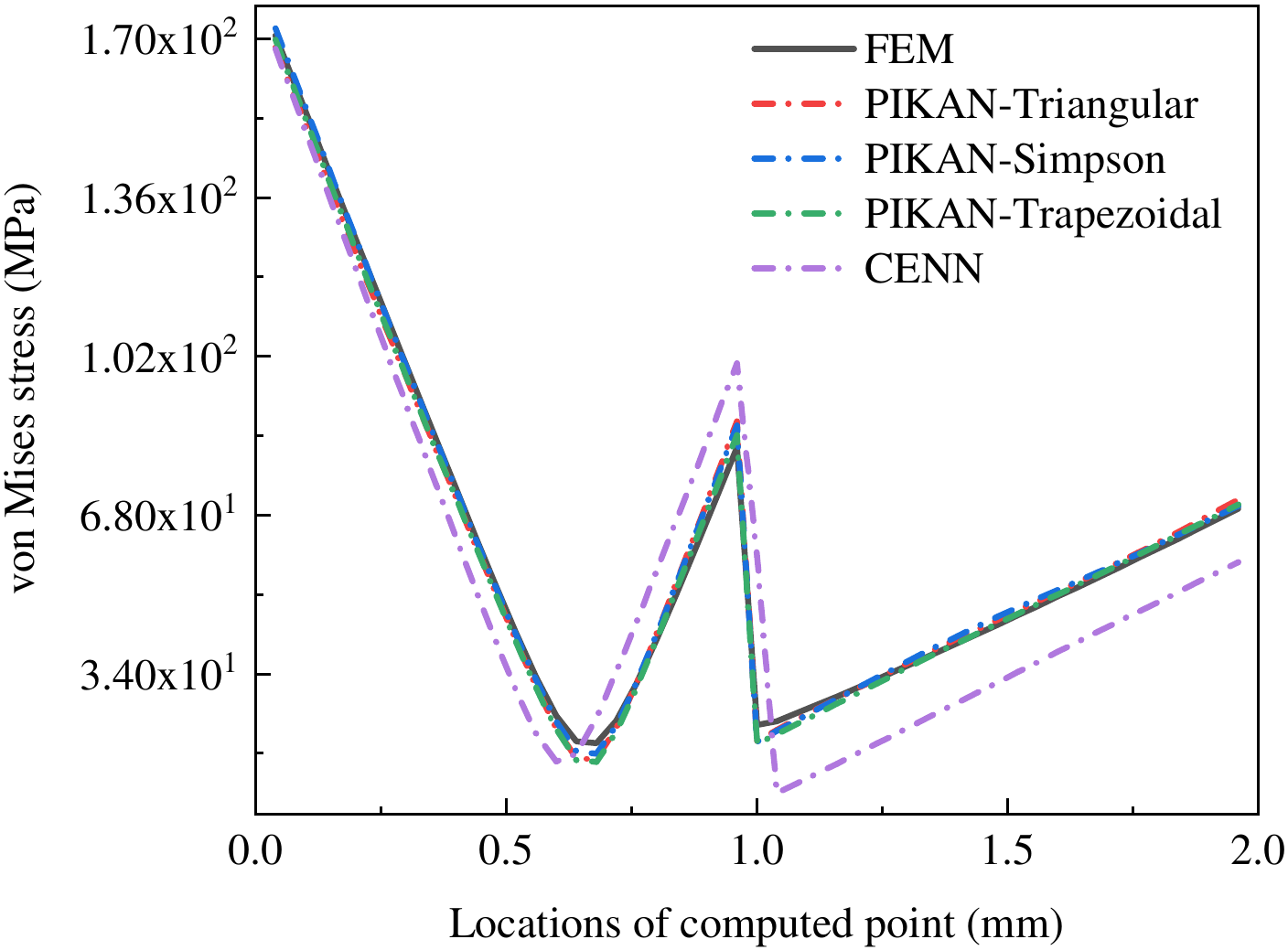}
        \caption{Straight interface}
        \label{fig:von_mises_x2_straight}
    \end{subfigure}
    \hfill
    \begin{subfigure}[b]{0.32\textwidth}
        \centering
        \includegraphics[width=\textwidth]{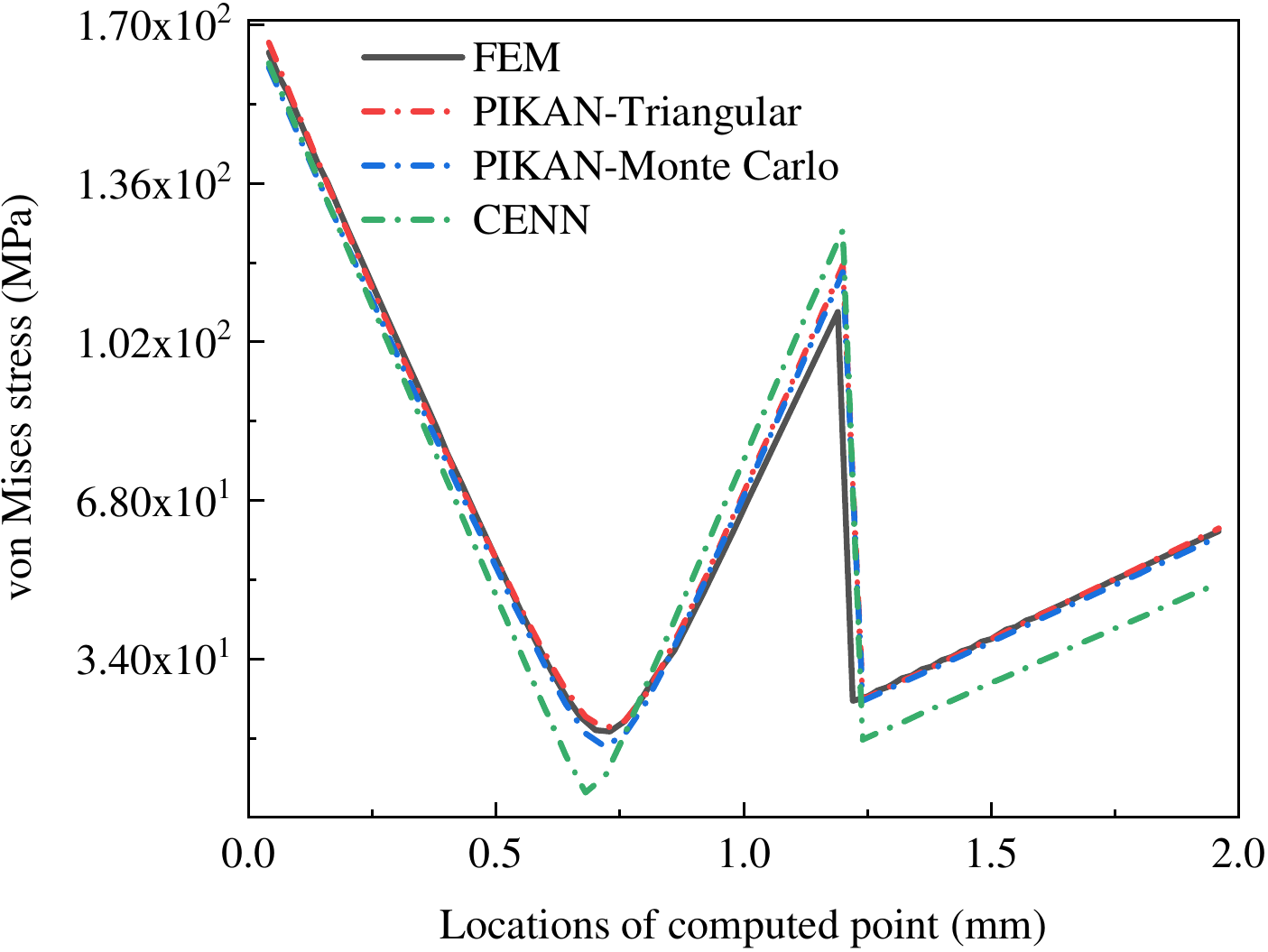}
        \caption{Wavy interface}
        \label{fig:von_mises_x2_wavy}
    \end{subfigure}
    \hfill
    \begin{subfigure}[b]{0.32\textwidth}
        \centering
        \includegraphics[width=\textwidth]{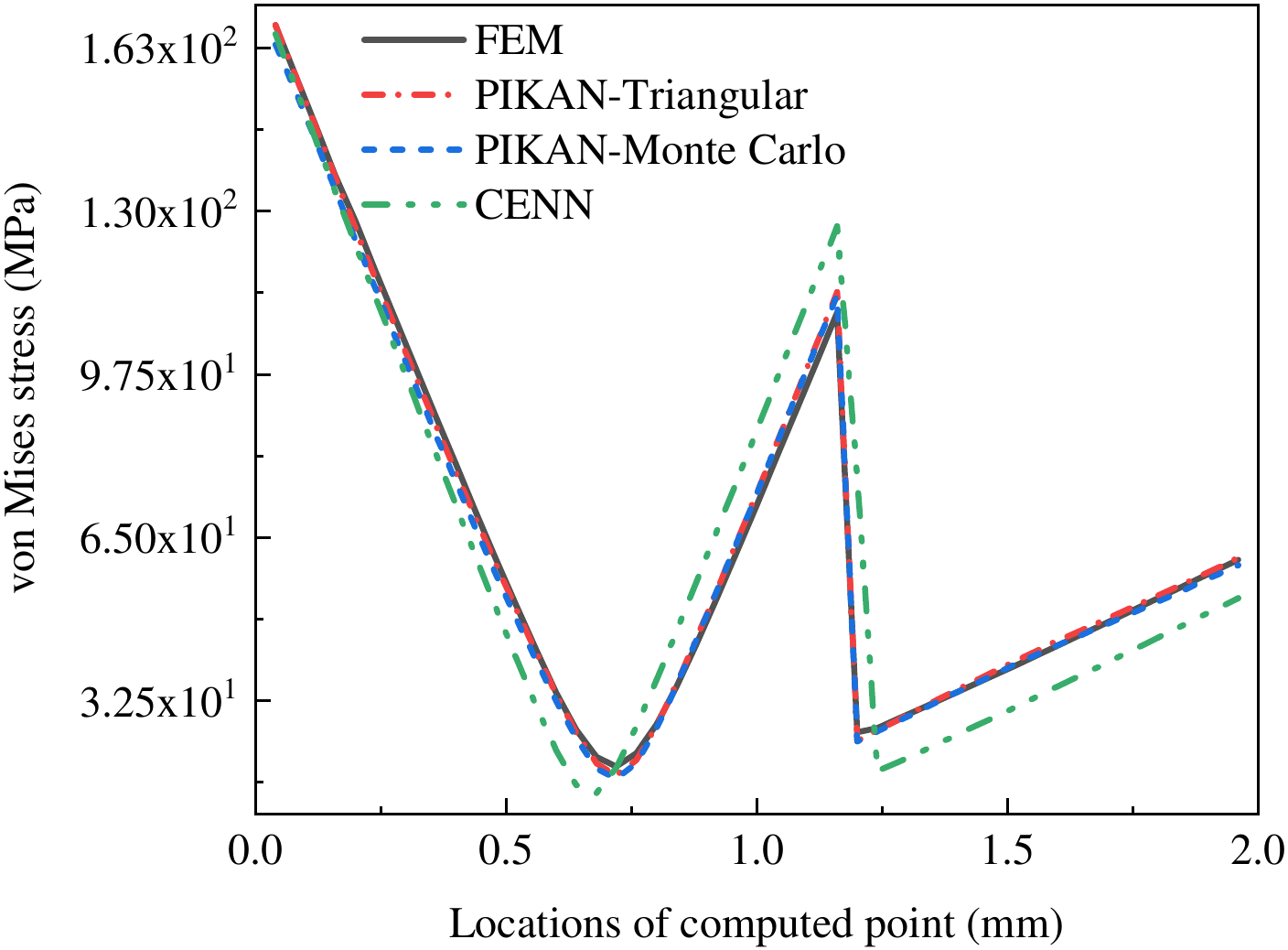}
        \caption{Stepped interface}
        \label{fig:von_mises_x2_stepped}
    \end{subfigure}
    
    \caption{von Mises stress at $x = 2$: PIKAN vs. CENN and FEM solutions.}
    \label{fig:von_mises_x2}
\end{figure}

We investigate PIKAN's performance in computing von Mises stress for cantilever beam models by analyzing stress profiles along the vertical line $x = 2$, $0 \leq y \leq 2$ (as shown in Fig. \ref{fig:cantilever_beam_models}). PIKAN results are compared with CENN predictions and FEM reference solutions. Fig.~\ref{fig:von_mises_x2} demonstrates that PIKAN predictions exhibit consistency with FEM solutions across all interface geometries. As shown in Fig. \ref{fig:von_mises_x2}, stress oscillations are observed, which can be attributed to two main factors. First, energy-based methods primarily minimize potential energy through displacement field optimization rather than directly enforcing pointwise stress equilibrium conditions. Second, stress solutions are computed as derivatives of displacement fields combined with constitutive relations, making them highly sensitive to small displacement errors, particularly near boundaries where geometric constraints and material property variations create stress concentrations. Consequently, the numerical accuracy of stress field is inherently lower than that of displacement~\cite{BaiAn2023}. Since PIKAN achieves superior displacement accuracy compared to CENN's subdomain approach in multi-material problems, it also yields more accurate stress results.
In contrast, CENN exhibits significant stress prediction errors, particularly at material interfaces, due to displacement discontinuities arising from its subdomain decomposition approach. The use of independent neural networks for different material regions results in pronounced stress fluctuations at material boundaries. Across all interface geometries (straight, wavy, and stepped configurations), CENN demonstrates computational instability in non-interface regions with substantial localized errors, confirming the limitations observed in the displacement analysis. Moreover, a comprehensive hyperparameter sensitivity analysis is presented in \hyperlink{AppendixB}{Appendix B}, revealing that optimal performance requires joint optimization of grid size, B-spline order, and network architecture rather than independent tuning.

To systematically evaluate PIKAN’s capability in handling multi-material interfaces and high parameter contrast problems, we further analyze a three-layer composite cantilever beam model beyond the previously discussed examples. As shown in Fig.~\ref{fig:cantilever_beam_models}d, the beam has length $L = 8$ and height $3H = 1.8$. The elastic moduli of the three materials are set with significant contrast: $E_1 = 3,000$ MPa, $E_2 = 20,000$ MPa, and $E_3 = 100,000$ MPa, with Poisson’s ratio $\nu = 0.3$ for all materials, the shear traction $T=6$ N/mm. This design creates a rigorous test scenario combining multiple heterogeneous interfaces with extreme parameter jumps to validate the method's performance.

\begin{figure}[htbp]
  \centering
  
  \begin{subfigure}[b]{0.85\textwidth}
      \centering
      \includegraphics[width=\textwidth]{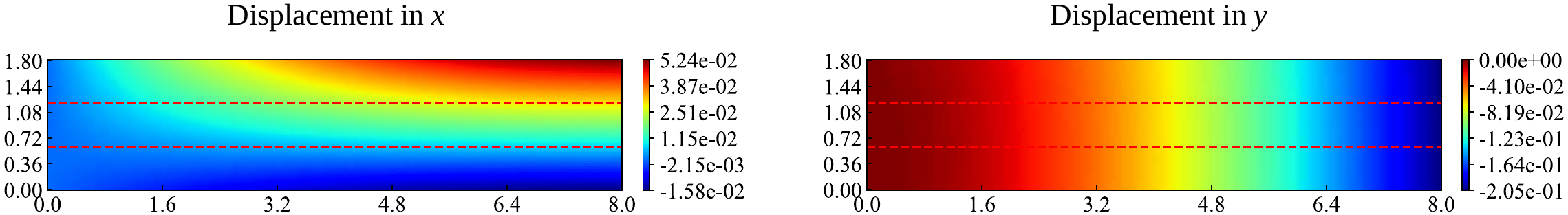}
      \caption{Displacement field - PIKAN}
      \label{fig:beam_three_PIKAN}
  \end{subfigure}
  
  \begin{subfigure}[b]{0.85\textwidth}
      \centering
      \includegraphics[width=\textwidth]{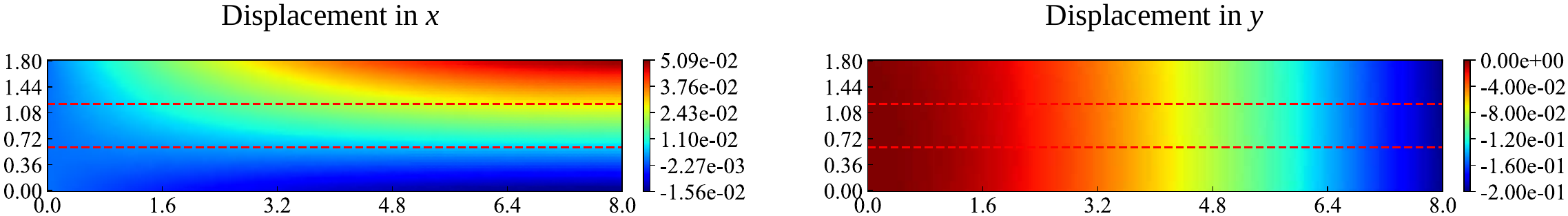}
      \caption{Displacement field - FEM}
      \label{fig:beam_three_FEM}
  \end{subfigure}
  
  \begin{subfigure}[b]{0.85\textwidth}
      \centering
      \includegraphics[width=\textwidth]{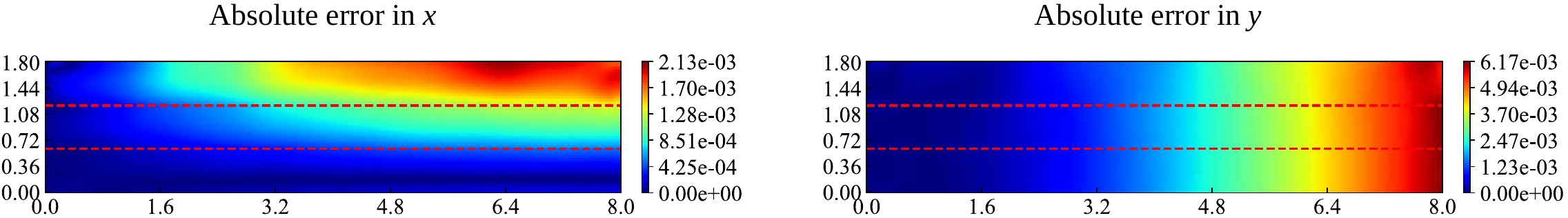}
      \caption{Absolute error}
      \label{fig:beam_three_absolute_error}
  \end{subfigure}
  
  \caption{Displacements and absolute rrrors (x- and y-directions) in three-layer composite cantilever beam model.}
  \label{fig:beam_three}
\end{figure}

Fig.~\ref{fig:beam_three}a, b compares PIKAN-predicted displacement fields $u_x$ and $u_y$ with reference FEM solutions, showing excellent agreement. More quantitatively, the pointwise absolute error distributions shown in Fig.~\ref{fig:beam_three}c demonstrate that even at material interfaces and in regions where elastic modulus spans two orders of magnitude, the maximum error remains well-controlled. This directly validates PIKAN’s accurate resolution capability for high-contrast parameters and interface discontinuities.

\

\subsection{Heterogeneous plate with central hole}

\begin{figure}[htbp]
\centering
\begin{subfigure}{0.45\textwidth}
    \centering
    \includegraphics[width=0.8\textwidth]{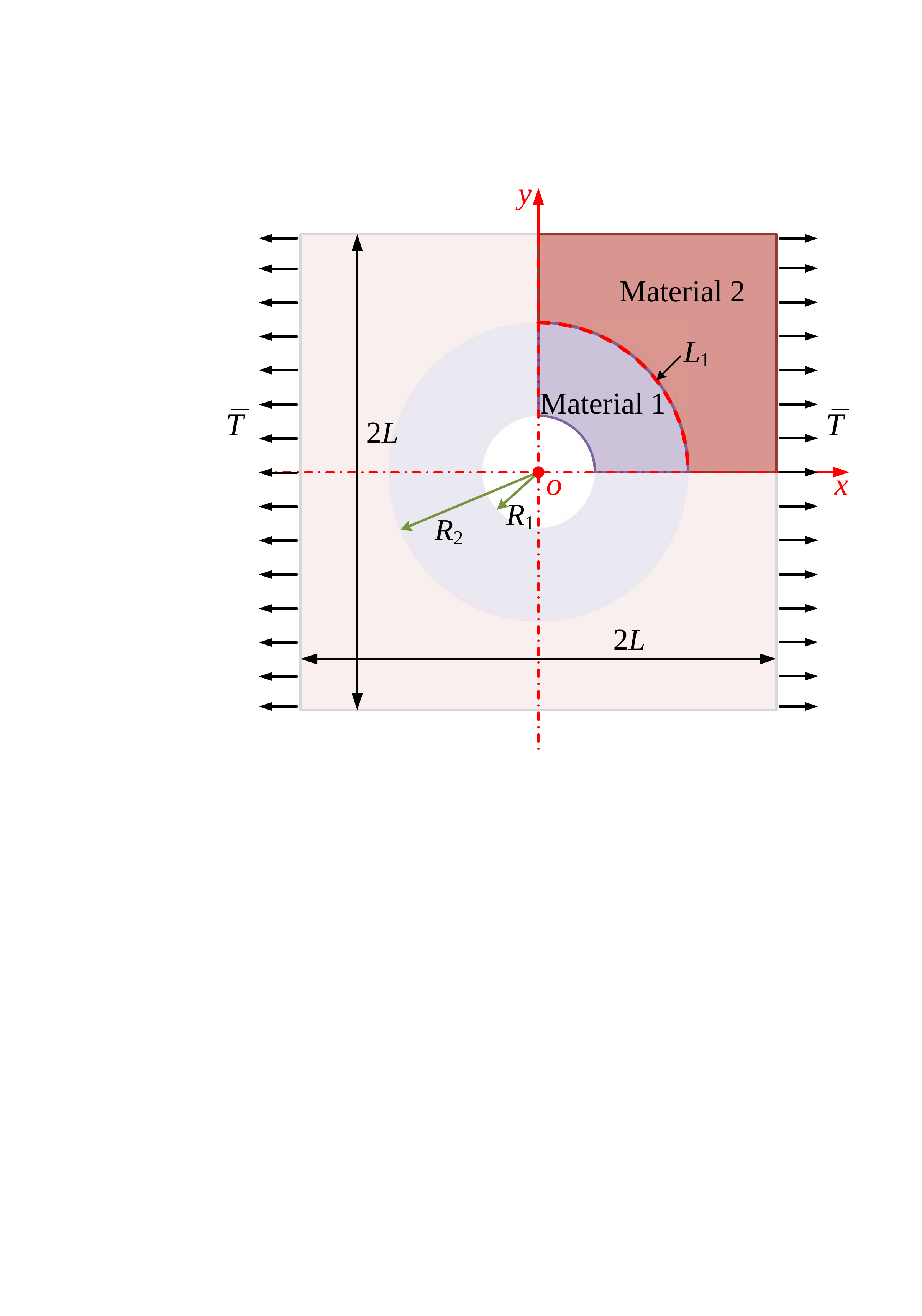}
    \caption{Full model}
    \label{fig:full_model}
\end{subfigure}
% \hfill
\begin{subfigure}{0.45\textwidth}
    \centering
    \includegraphics[width=0.8\textwidth]{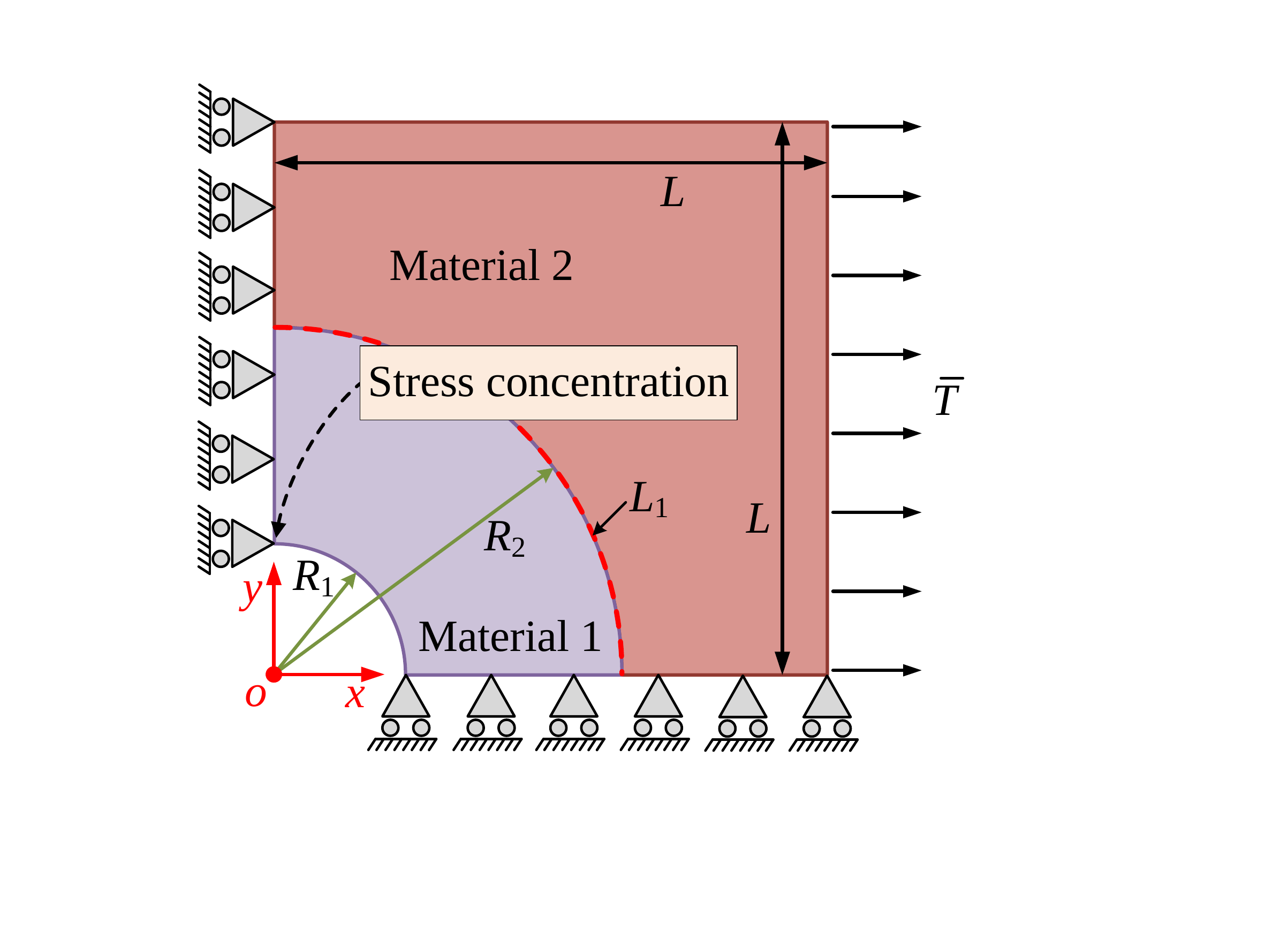}
    \caption{Quarter model}
    \label{fig:quarter_model}
\end{subfigure}
\caption{Heterogeneous plate with central hole model.}
\label{fig:heterogeneous_plate_model}
\end{figure}

We apply PIKAN to analyze a multi-material square plate with a central hole under tensile loading, as shown in Fig.~\ref{fig:heterogeneous_plate_model}a. Due to symmetry, only a quarter-plate model with dimensions $L \times L$, where $L = 21$, is analyzed. The coordinate origin $o$ is located at the bottom-left corner of the quarter model.
As illustrated in Fig.~\ref{fig:heterogeneous_plate_model}b, the central circular hole has radius $R_1 = 5$ with traction-free boundaries. The material interface (indicated by the red dashed line) is a circular boundary with radius $R_2 = 13$, where material 1 occupies the inner region and material 2 occupies the outer region. A uniform tensile load $\bar{T} = 120$ N/mm is applied to the right boundary of the quarter model. Symmetry boundary conditions are imposed along the left and bottom edges. The heterogeneous plate consists of two materials: Material 1 with Young's modulus $E_1 = 10,000$ MPa and Poisson's ratio $\nu_1 = 0.3$, and Material 2 with $E_2 = 1,000$ MPa and $\nu_2 = 0.4$.

\begin{figure}[htbp]
  \centering
  % ---------- 第一组 ----------
  \begin{minipage}{0.49\linewidth}
    \centering
    \begin{subfigure}[t]{0.49\linewidth}
      \includegraphics[width=\linewidth]{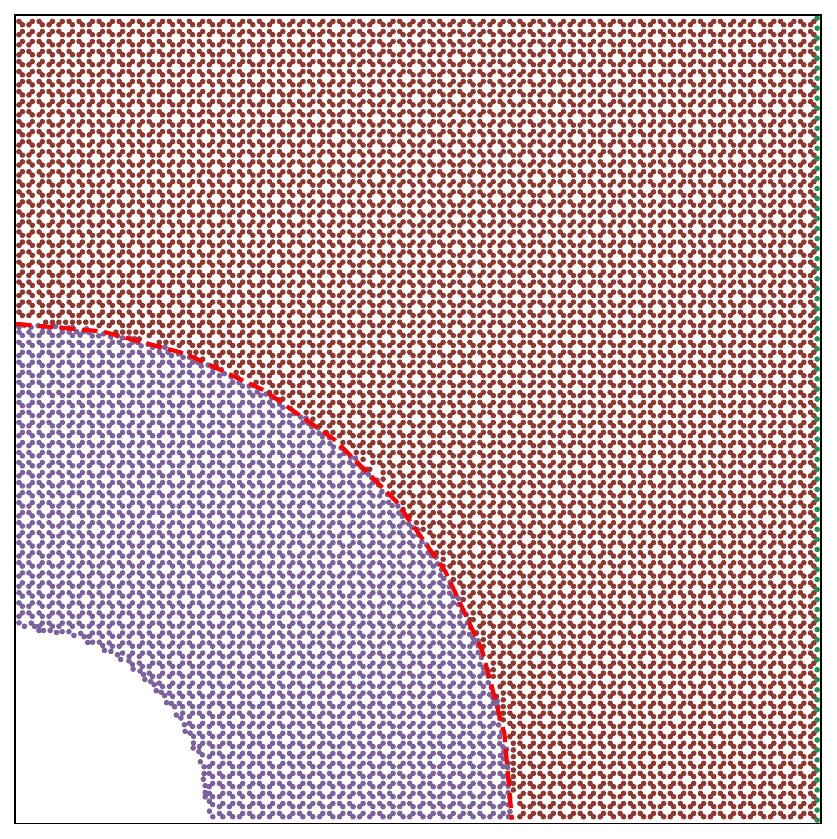}
      \caption{Triangular centroid distribution scheme}
    \end{subfigure}%
    \hspace{0.01em}
    \begin{subfigure}[t]{0.49\linewidth}
      \includegraphics[width=\linewidth]{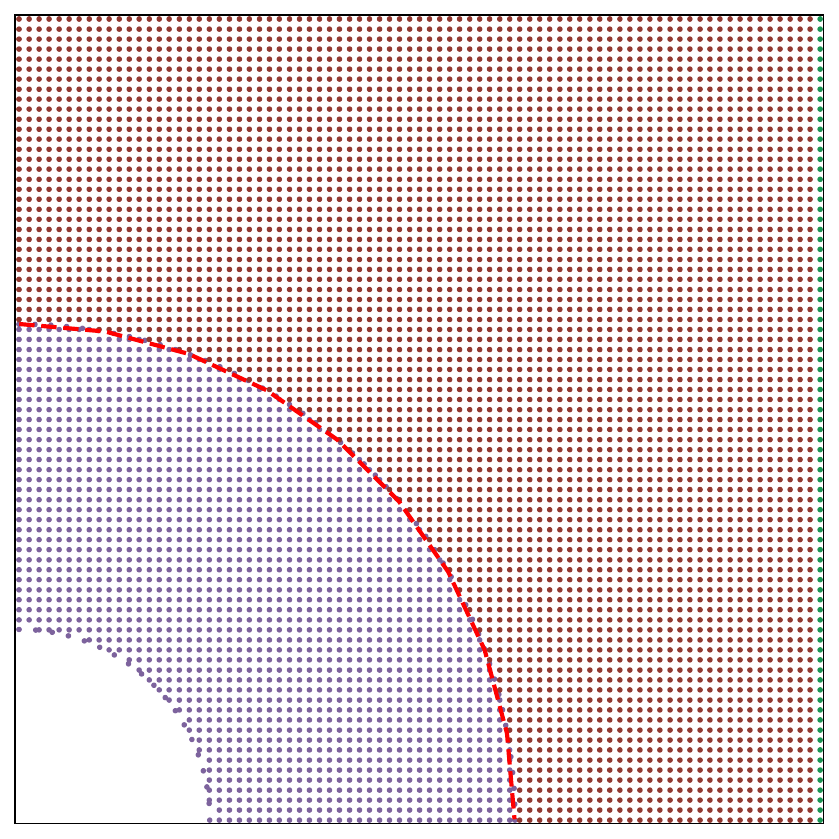}
      \caption{Uniform point distribution scheme}
    \end{subfigure}
    \caption{Sample point distribution for heterogeneous plate model.}
    \label{fig:sample_point_distribution}
  \end{minipage}%
  \hspace{0.5em} 
  % ---------- 第二组 ----------
  \begin{minipage}{0.48\linewidth}
    \centering
    \includegraphics[width=0.9\linewidth]{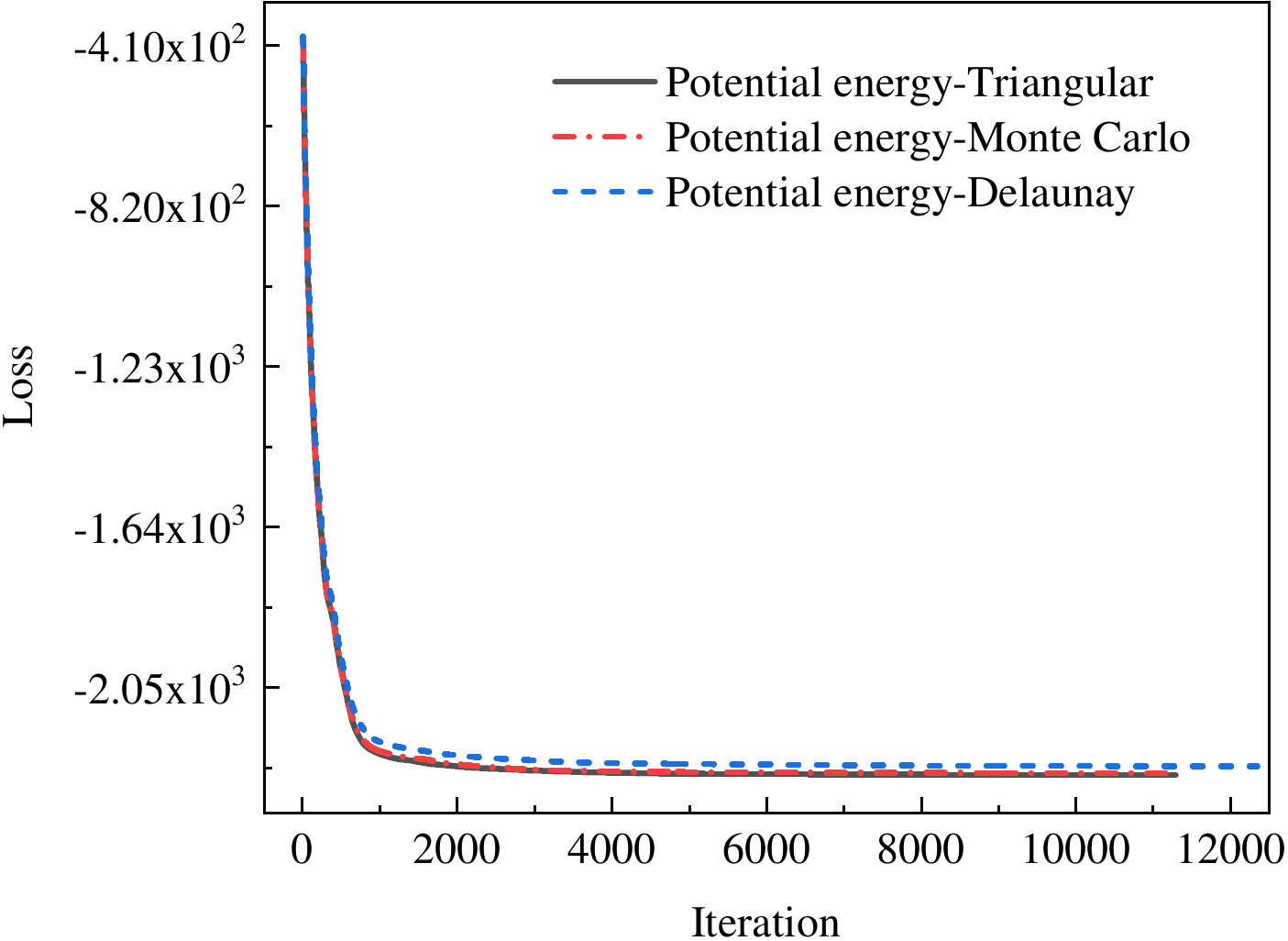}
    \caption{Loss function evolution during training for heterogeneous plate model.}
    \label{fig:energy_convergence}
  \end{minipage}
\end{figure}

The KAN architecture for this model is configured as $[2,5,5,5,2]$ with grid size 20, B-spline order $k=3$, and grid range $[0,1]$. The neural network contains 1,750 trainable parameters. Fig.~\ref{fig:sample_point_distribution} shows the sample point distribution patterns applied to the multi-material plate model. Purple points represent sample points in Material 1 (inner region), maroon points represent sample points in Material 2 (outer region), and green points along the loaded boundaries are used to apply natural boundary conditions, totaling 211 points.
Triangle centroid points comprise 22,786 interior points in Material 1 and 61,624 interior points in Material 2. For uniformly distributed sample points, Material 1 contains 11,574 interior points and Material 2 contains 31,121 interior points. This difference in point density allows for investigating the effect of sampling resolution on solution accuracy. For comparison, results from CENN are also presented. In CENN, each material domain employs a separate MLP with 6 hidden layers of 30 neurons each, totaling 4,802 trainable parameters. The training sample points follow the triangular distribution pattern shown in Fig.~\ref{fig:sample_point_distribution}a.

\begin{figure}[htbp]
  \centering
  % ---------- 第一组：PIKAN vs. FEM ----------
  \begin{minipage}{0.49\linewidth}
    \centering
    \begin{subfigure}{0.49\linewidth}
      \includegraphics[width=\linewidth]{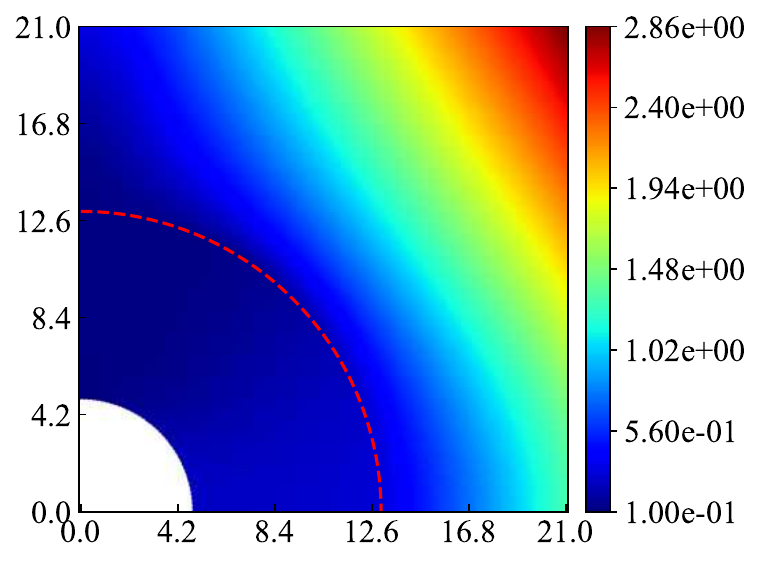}
      \caption{PIKAN results}
      \label{fig:pikan_solution}
    \end{subfigure}%
    \hspace{0.01em}
    \begin{subfigure}{0.49\linewidth}
      \includegraphics[width=\linewidth]{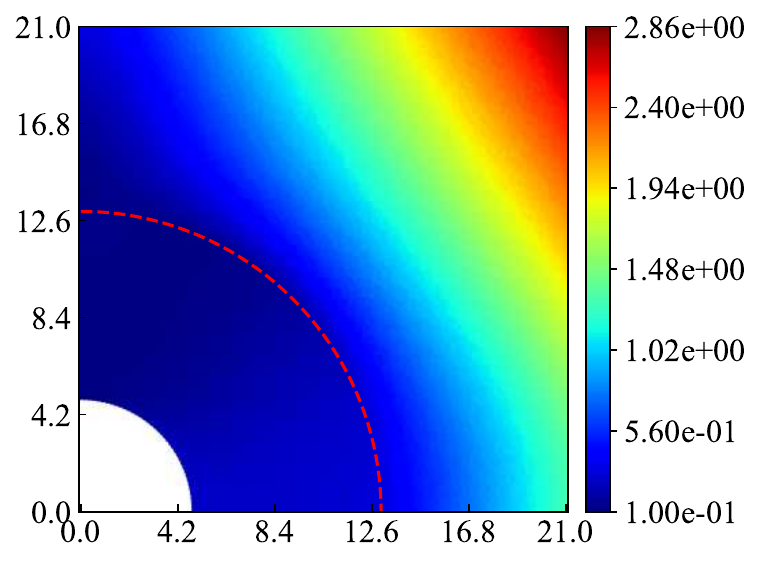}
      \caption{FEM solutions}
      \label{fig:fem_solution}
    \end{subfigure}
    \caption{Displacement $u_{\text{mag}}$ for heterogeneous plate model: PIKAN vs. FEM.}
    \label{fig:displacement_comparison}
  \end{minipage}%
  \hspace{0.01em}   % ← 两组之间距离，随意改
  % ---------- 第二组：误差 ----------
  \begin{minipage}{0.49\linewidth}
    \centering
    \begin{subfigure}[t]{0.49\linewidth}
      \includegraphics[width=\linewidth]{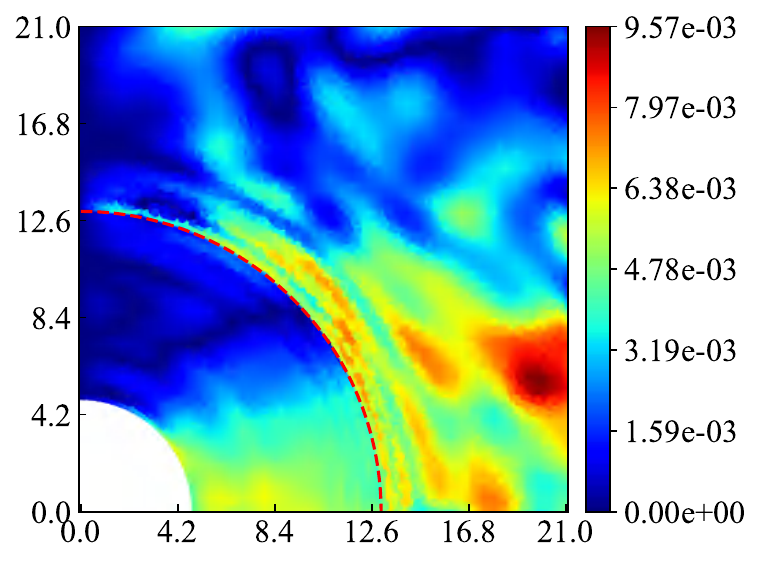}
      \caption{$u_x$ absolute error}
      \label{fig:ux_error_heterogeneous_plate}
    \end{subfigure}%
    \hspace{0.01em}
    \begin{subfigure}[t]{0.49\linewidth}
      \includegraphics[width=\linewidth]{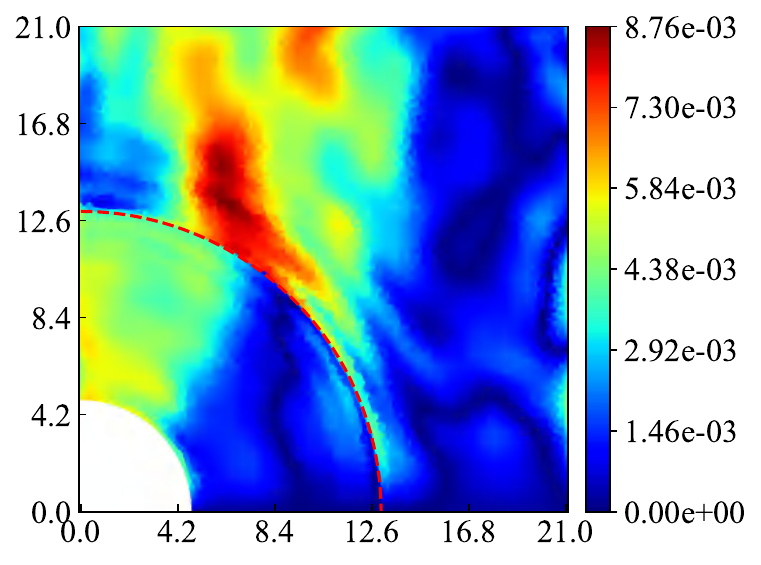}
      \caption{$u_y$ absolute error}
      \label{fig:uy_error_heterogeneous_plate}
    \end{subfigure}
    \caption{Absolute error contours for displacement fields in heterogeneous plate model.}
    \label{fig:error_contours_heterogeneous_plate}
  \end{minipage}
\end{figure}

Fig.~\ref{fig:energy_convergence} demonstrates the evolution of energy-based loss functions during training for the heterogeneous plate model. As the neural network undergoes iterative optimization, the total potential energy gradually converges to its minimum value, validating the principle of minimum potential energy. All numerical integration schemes exhibit similar convergence behavior, but the Delaunay scheme exhibits marginally weaker convergence in the loss compared to the other two schemes. Fig.~\ref{fig:displacement_comparison} presents a detailed comparison between PIKAN predictions (using triangular integration) and FEM reference solutions for displacement magnitude $u_{\text{mag}}$ (calculated according to Eq.~\eqref{eq:displacement_magnitude}).
\begin{equation}
\label{eq:displacement_magnitude}
u_{\text{mag}} = \sqrt{u_x^2 + u_y^2}
\end{equation} 
The results demonstrate excellent agreement between PIKAN and FEM solutions, with PIKAN successfully capturing the displacement distribution around the central hole and the smooth transition across the material interface. The displacement patterns clearly show the influence of material property differences, with higher displacements occurring in the softer outer material (Material 2).
Fig.~\ref{fig:error_contours_heterogeneous_plate} shows the pointwise absolute error contours for displacement components $u_x$ and $u_y$ compared to FEM reference solutions. The error distribution reveals that PIKAN maintains high accuracy throughout most of the domain, with maximum errors typically below $9.57\times 10^{-3}$.

\begin{figure}[htbp]  
    \centering
    \begin{minipage}{0.45\textwidth}
        \centering
        \includegraphics[width=0.99\textwidth]{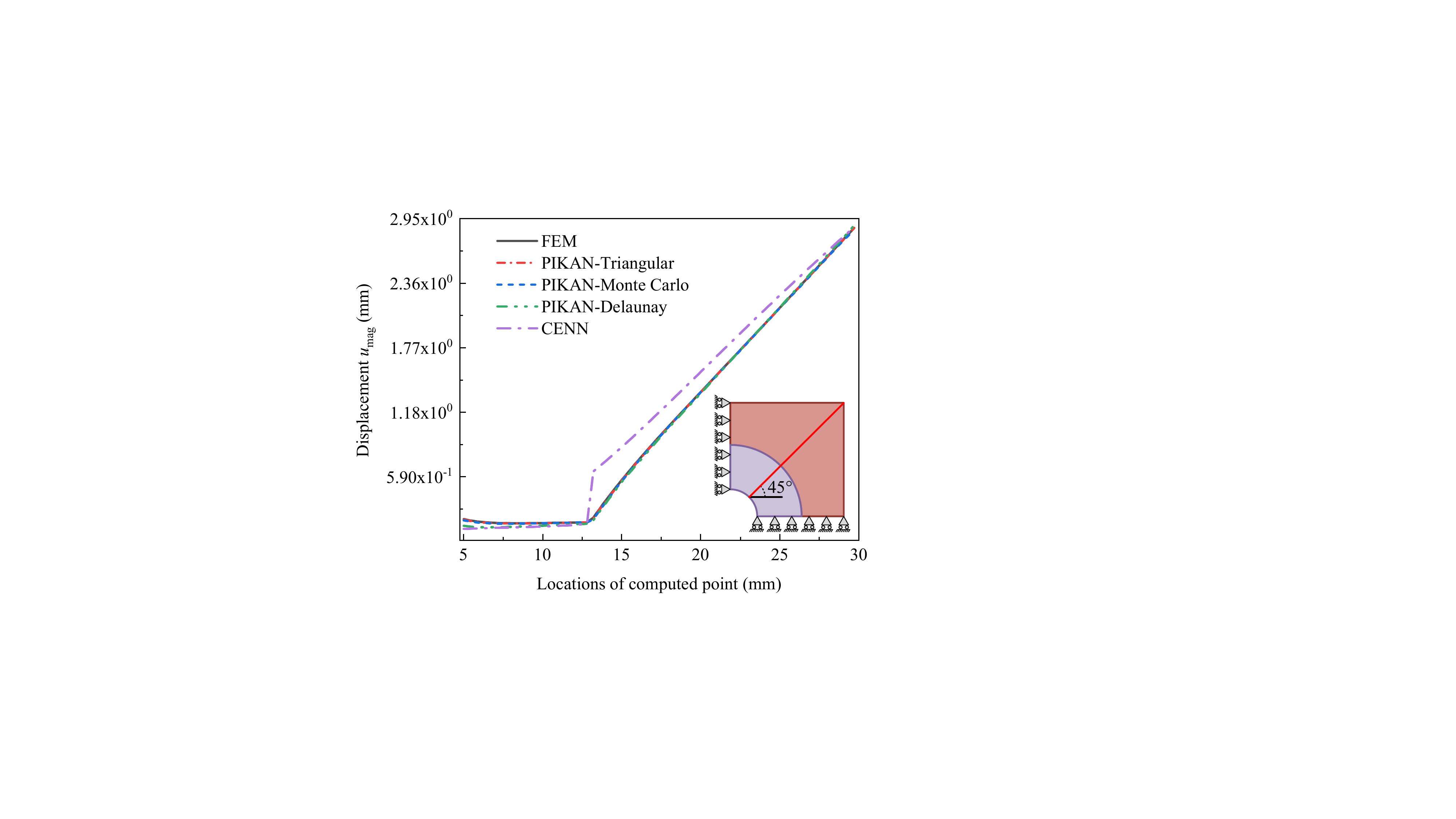}
        \caption{$u_{\text{mag}}$ along diagonal line.}
        \label{fig:diagonal_displacement}
    \end{minipage}
    \hspace{1em}
    \begin{minipage}{0.45\textwidth}
    \centering
        \includegraphics[width=0.93\textwidth]{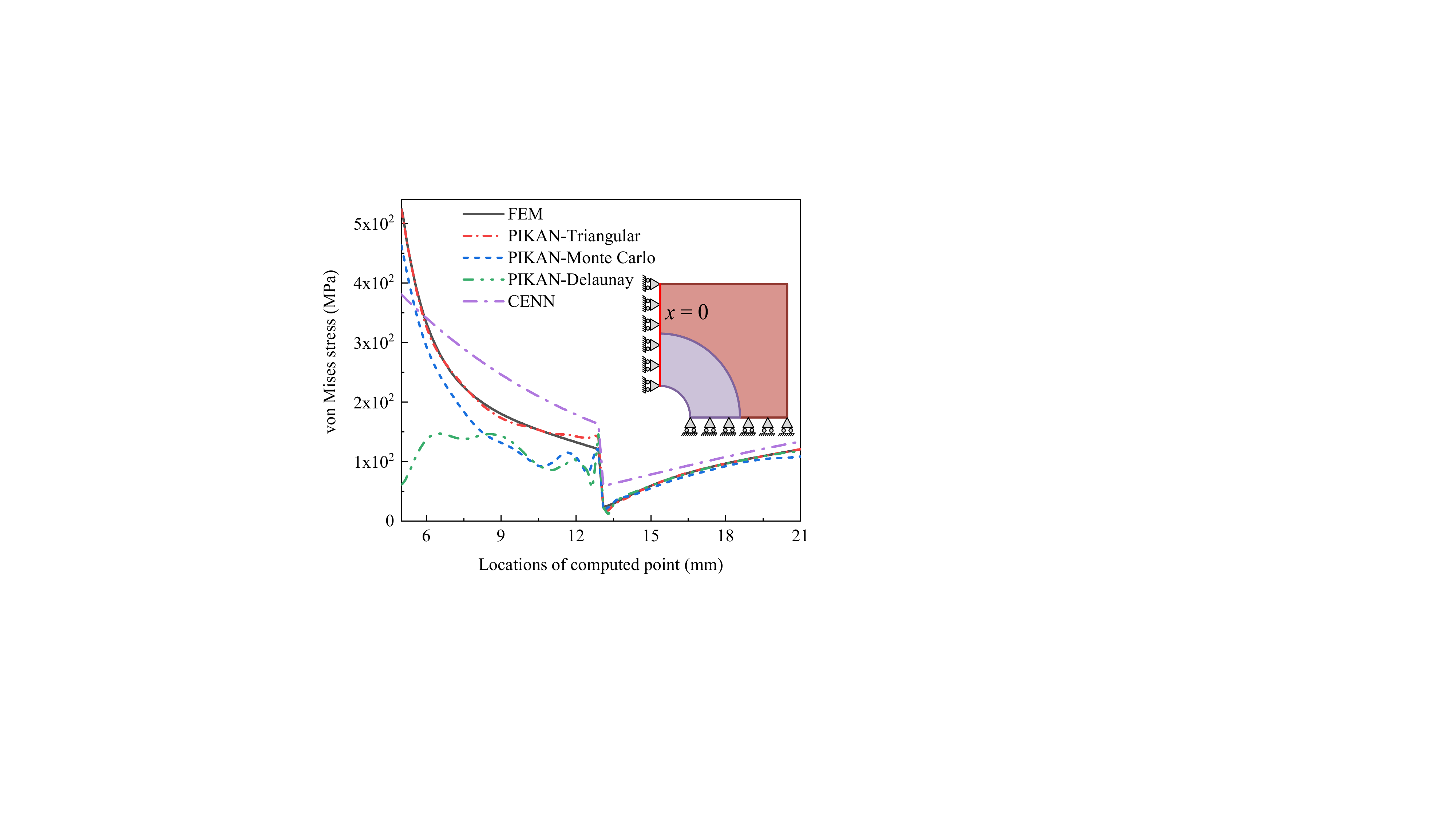}
        \caption{Stress along the line $x = 0$.}
        \label{fig:stress_analysis}
\end{minipage}
\end{figure}

To examine solution behavior across the entire domain, displacement profiles along the diagonal line $\theta = 45^\circ$ (from the hole edge to the outer boundary, indicated by the red line) are analyzed. Fig.~\ref{fig:diagonal_displacement} shows the comparison of displacement magnitude $u_{\text{mag}}$ along this diagonal line. The results demonstrate that triangular integration and Monte Carlo integration schemes produce solutions that agree well with the FEM reference, maintaining accuracy across both material regions and the interface transition. However, the Delaunay integration scheme exhibits significantly inferior accuracy. This scheme assigns point weights as one-third of the associated triangle area sum, and this discrete integral representation using uniformly sampled points introduces integration errors that prevent the loss function from accurately approximating the true potential energy. Consequently, the network trained on this basis produces relatively high solution errors. This result confirms a critical principle: in variational principle-based physics-informed neural networks, the accuracy of the numerical integration scheme used in the loss function directly affects the solver's final output accuracy. In contrast, CENN predictions show significant deviations from the reference solution throughout the domain, failing to accurately capture the displacement field distribution patterns. This demonstrates the inherent challenges of domain decomposition approaches that rely on penalty terms for interface continuity, particularly for problems with circular interfaces where geometric complexity compounds convergence difficulties.

Next, we evaluate PIKAN's performance by comparing the stress results from PIKAN, FEM and CENN along the vertical line $x = 0$ (from $y = 5$ to $y = 21$). This analysis line extends from the hole boundary through both material regions to the outer boundary, providing insight into stress prediction accuracy across material interfaces and high-gradient regions. The results are presented in Fig.~\ref{fig:stress_analysis}.
The analysis reveals that triangular integration achieves the highest accuracy, closely converging to the FEM reference solution throughout the domain. In contrast, Monte Carlo and Delaunay integration schemes show poor prediction performance in the Material 1 domain (high-stiffness inner region), where stress concentrations are more pronounced and strain energy density varies dramatically. The numerical approximation errors of these integration schemes are substantial in this region, preventing the loss function from providing sufficiently accurate gradients for network training, thereby degrading final solution accuracy. This differential performance suggests that stress calculations are particularly sensitive to integration scheme choice in regions with steep stress gradients and high material property contrasts.
CENN predictions follow the general trend of the reference solution but with substantially larger errors, highlighting the challenges of maintaining stress accuracy in domain decomposition approaches that rely on penalty terms for interface continuity.

\subsection{PIKAN application to DBC substrate structure}

PIKAN is applied to analyze DBC (Direct Bonding Copper) substrates featuring multi-material domains. DBC substrates are widely used in power electronics, high-power modules, and aerospace applications due to their excellent thermal conductivity and mechanical strength. Material property differences and thermal/mechanical stresses during manufacturing or operation cause these substrates to warp, affecting assembly performance and compromising equipment functionality and reliability~\cite{yu2021thermal}.

\begin{figure}[htbp]
  \centering
  % ---------- 一行三图 ----------
  \begin{subfigure}[t]{0.5\textwidth}   % 顶端对齐
    \centering
    \includegraphics[width=\linewidth]{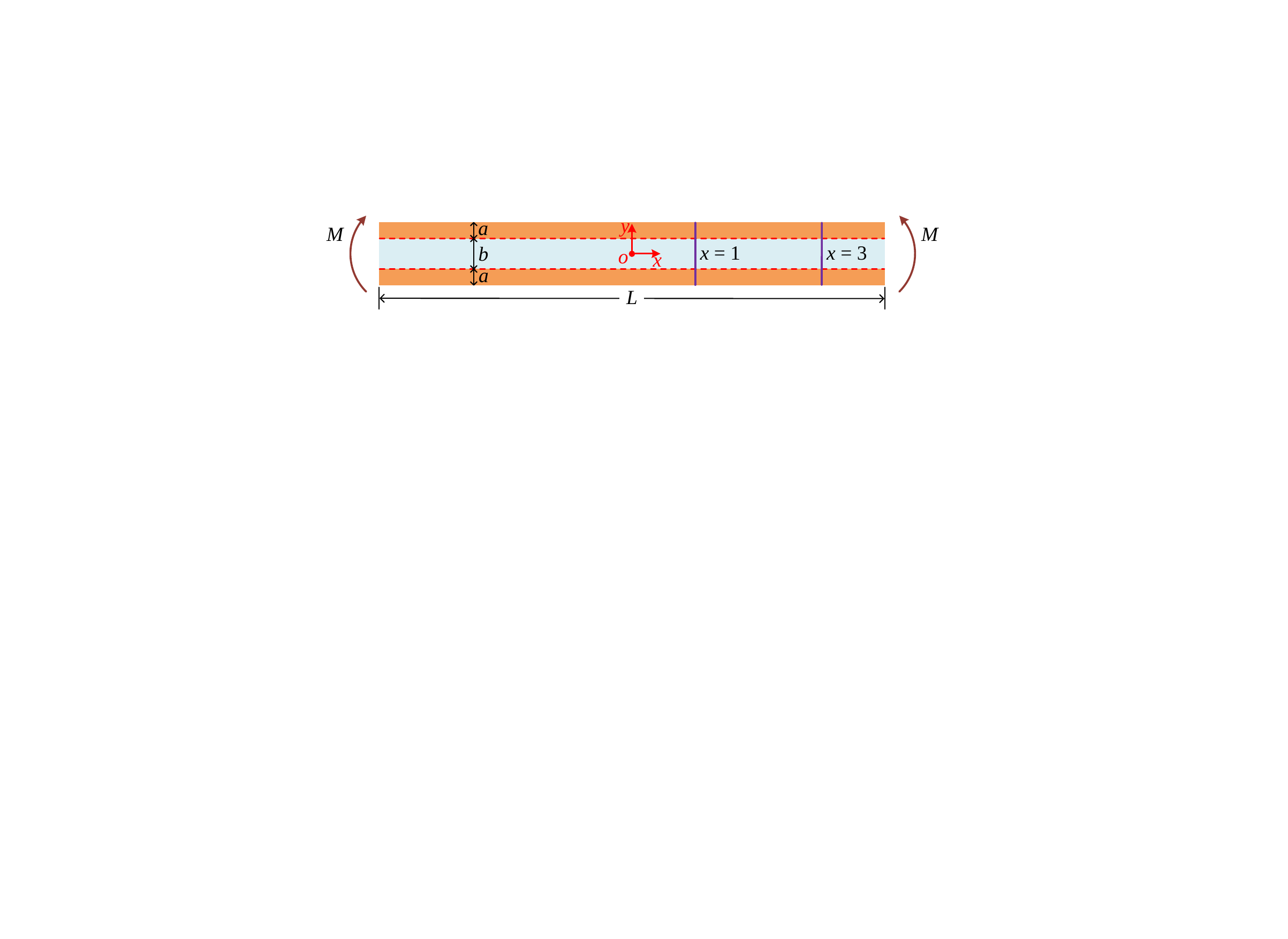}
    \caption{Full setup with applied moments}
    \label{fig:full_setup}
  \end{subfigure}\hfill
  \begin{subfigure}[t]{0.26\textwidth}
    \centering
    \includegraphics[width=\linewidth]{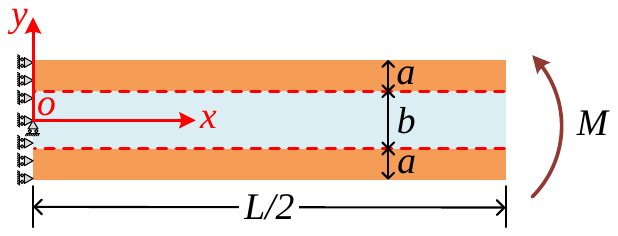}
    \caption{Half model}
    \label{fig:half_model}
  \end{subfigure}\hfill
  \begin{subfigure}[t]{0.2\textwidth}
    \centering
    \includegraphics[width=\linewidth]{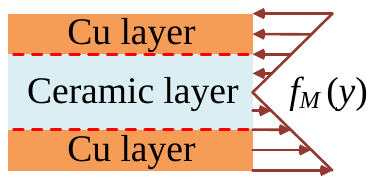}
    \caption{Equivalent loading}
    \label{fig:equivalent_loading}
  \end{subfigure}
  \caption{DBC substrate bending problem.}
  \label{fig:dbc_model}
\end{figure}

A simplified DBC substrate warping problem is modeled as shown in Fig.~\ref{fig:dbc_model}a. The substrate has total length $L = 8$ with an Al$_2$O$_3$ ceramic core of thickness $b = 0.48$ sandwiched between upper and lower copper layers, each with thickness $a = 0.26$. The substrate experiences equal and opposite bending moments $M = 100$ N$\cdot$mm at both ends, creating pure bending deformation.
Due to symmetry about the midplane, only the right half of the model requires analysis, as shown in Fig.~\ref{fig:dbc_model}b. The coordinate origin is located at the center of the substrate cross-section. Symmetry boundary conditions are applied at the left boundary, while the bending moment is applied as an equivalent distributed force at the right end. The equivalent force distribution, illustrated in Fig.~\ref{fig:dbc_model}c, is given by
\begin{equation}
    f_M(y) = \frac{12M}{h^3}y = 1200y \quad \text{for } y \in (-0.5, 0.5) 
    \label{eq:equivalent_force}
\end{equation}
where $h = 1.0$ is the total substrate thickness. This linear distribution ensures that the resultant force is zero and the moment equals the applied bending moment $M$. The DBC substrate consists of copper layers with Young's modulus $E_{\text{Cu}} = 128,000$ MPa and Poisson's ratio $\nu_{\text{Cu}} = 0.34$, and an alumina ceramic core with $E_{\text{Al}_2\text{O}_3} = 270,000$ MPa and $\nu_{\text{Al}_2\text{O}_3} = 0.28$.

\begin{figure}[htbp]
    \centering
    \begin{subfigure}{0.4\textwidth}
        \centering
        \includegraphics[width=\textwidth]{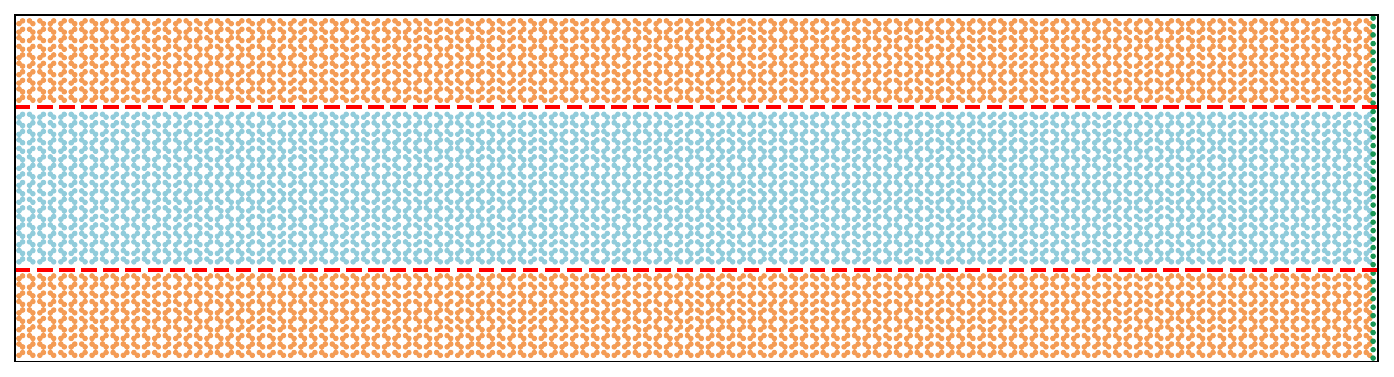}
        \caption{Triangular distribution scheme}
        \label{fig:triangular_distribution}
    \end{subfigure}
     \hspace{2em}
    \begin{subfigure}{0.4\textwidth}
        \centering
        \includegraphics[width=\textwidth]{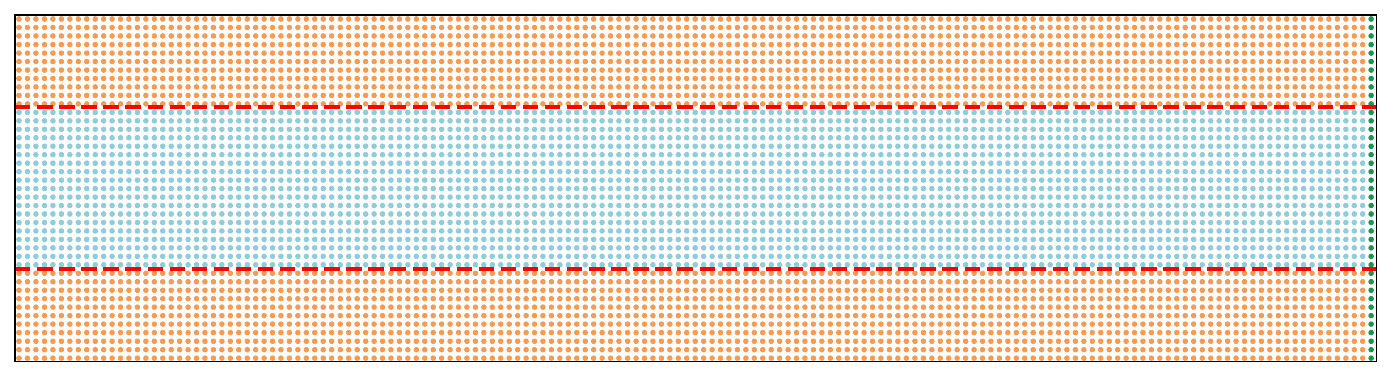}
        \caption{Uniform distribution scheme}
        \label{fig:uniform_distribution}
    \end{subfigure}
    \caption{Training point distribution schemes for DBC substrate model.}
    \label{fig:dbc_sample_distribution}
\end{figure}

A KAN with architecture $[2,5,5,5,2]$ is employed to predict the displacement field $(u_x, u_y)$ with grid size 10, B-spline order $k = 2$, and grid range $[-1,1]$. The neural network contains 980 trainable parameters.
Fig.~\ref{fig:dbc_sample_distribution} shows the training point distribution schemes for the DBC substrate model. Due to the layered structure of the DBC substrate, sample points are distributed across three material regions: the upper copper layer, the ceramic core, and the lower copper layer. In the centroid sampling scheme, the two copper layers contain a total of 166,400 sample points (83,200 points per layer) and the ceramic core contains 150,400 sample points. In the uniform sampling strategy, the copper layers contain a total of 84,906 sample points (42,453 points per layer) and the ceramic core contains 76,095 sample points. 
For both distribution schemes, the right-end natural boundary condition is discretized using 201 uniformly distributed sample points to accurately represent the linear force distribution given by Eq.~\eqref{eq:equivalent_force}.

\begin{figure}[htbp]
    \centering
    \includegraphics[width=0.45\textwidth]{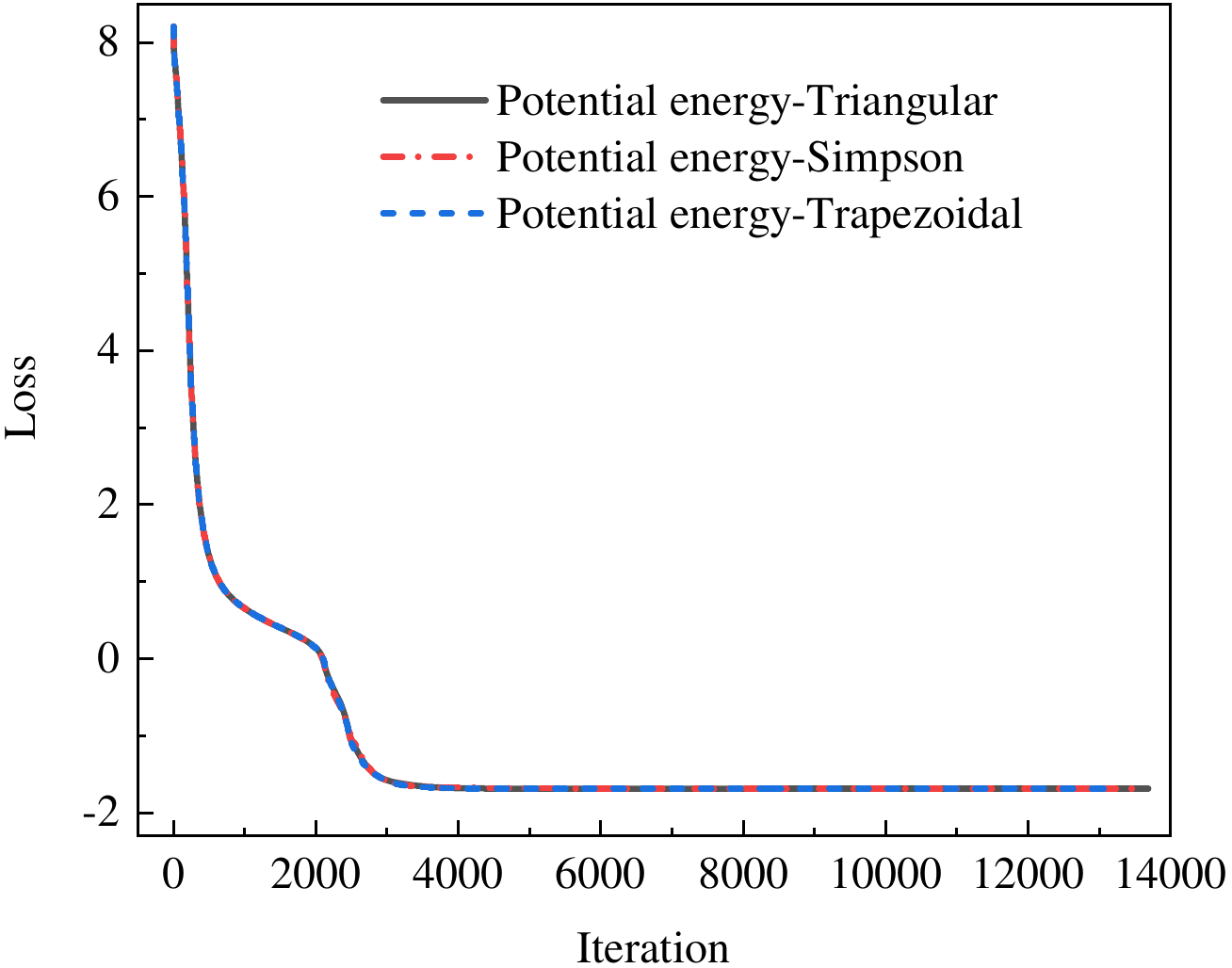}
    \caption{Evolution of loss functions during training for the DBC substrate model.}
    \label{fig:dbc_energy_convergence}
\end{figure}

\begin{figure}[htbp]
    \centering
    \begin{subfigure}{0.32\textwidth}
        \centering
        \includegraphics[width=\textwidth]{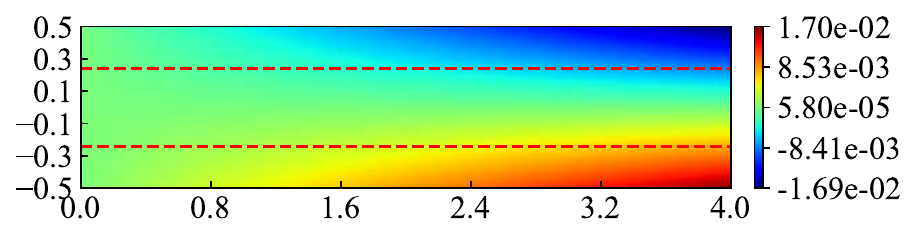}
        \caption{PIKAN: $u_x$}
        \label{fig:pikan_ux}
    \end{subfigure}
    \hfill
    \begin{subfigure}{0.32\textwidth}
        \centering
        \includegraphics[width=\textwidth]{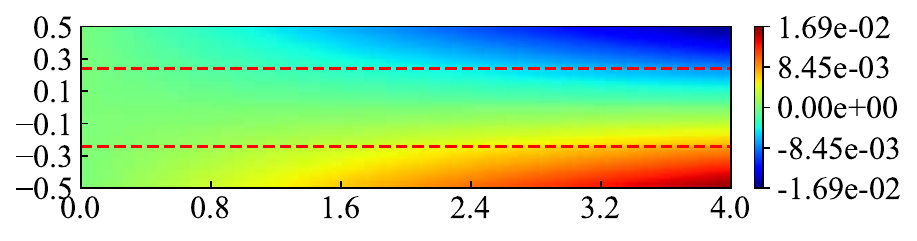}
        \caption{FEM: $u_x$}
        \label{fig:fem_ux}
    \end{subfigure}
    \hfill
    \begin{subfigure}{0.32\textwidth}
        \centering
        \includegraphics[width=\textwidth]{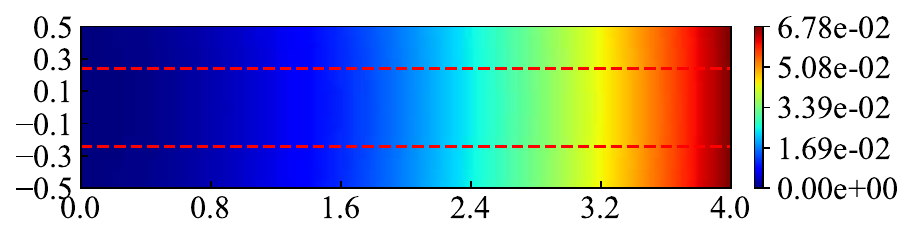}
        \caption{PIKAN: $u_y$}
        \label{fig:pikan_uy}
    \end{subfigure}
    
    \begin{subfigure}{0.32\textwidth}
        \centering
        \includegraphics[width=\textwidth]{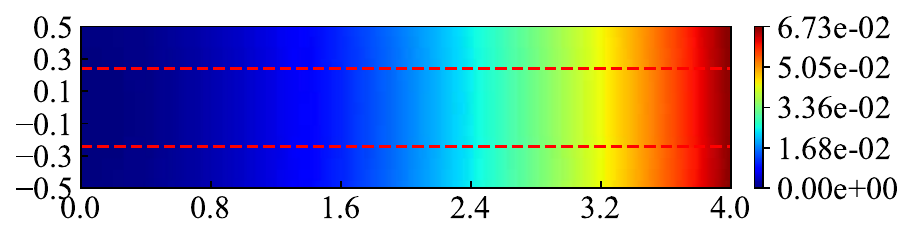}
        \caption{FEM: $u_y$}
        \label{fig:fem_uy}
    \end{subfigure}
    \hfill
    \begin{subfigure}{0.32\textwidth}
        \centering
        \includegraphics[width=\textwidth]{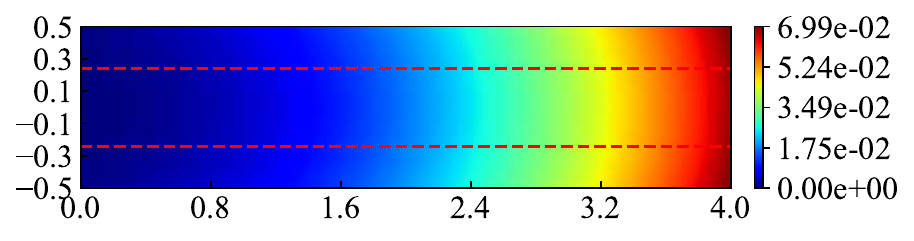}
        \caption{PIKAN: $u_{\text{mag}}$}
        \label{fig:pikan_umag}
    \end{subfigure}
    \hfill
    \begin{subfigure}{0.32\textwidth}
        \centering
        \includegraphics[width=\textwidth]{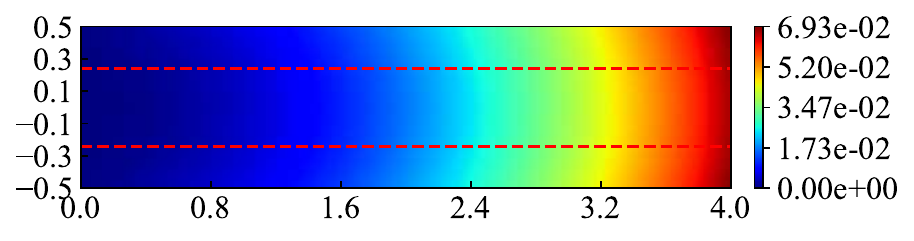}
        \caption{FEM: $u_{\text{mag}}$}
        \label{fig:fem_umag}
    \end{subfigure}
    \caption{Displacement field solutions for DBC substrate model: PIKAN vs. FEM.}
    \label{fig:dbc_displacement_comparison}
\end{figure}

\begin{figure}[htbp]
    \centering
    \begin{subfigure}{0.45\textwidth}
        \centering
        \includegraphics[width=0.9\textwidth]{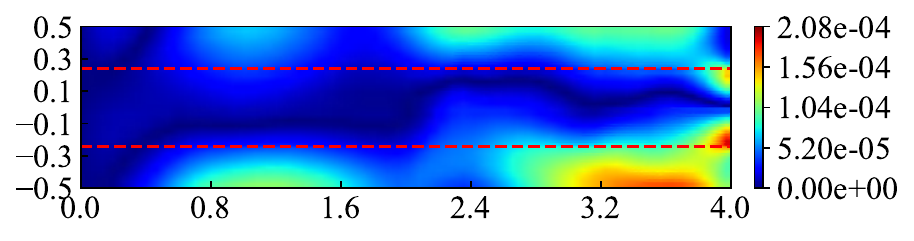}
        \caption{$u_x$ absolute error}
        \label{fig:error_ux}
    \end{subfigure}
    % \hfill
    \begin{subfigure}{0.45\textwidth}
        \centering
        \includegraphics[width=0.9\textwidth]{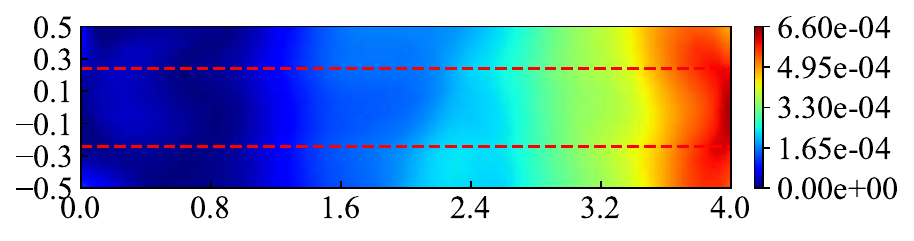}
        \caption{$u_y$ absolute error}
        \label{fig:error_uy}
    \end{subfigure}
    \caption{Absolute error contours for displacement fields in DBC substrate model.}
    \label{fig:dbc_error_analysis}
\end{figure}

Fig.~\ref{fig:dbc_energy_convergence} illustrates the evolution of the potential energy loss function for the DBC substrate warpage model.
Fig.~\ref{fig:dbc_displacement_comparison} presents a comprehensive comparison of displacement field solutions between PIKAN predictions (using simpson integration) and FEM reference solutions for the DBC substrate bending problem. The results demonstrate excellent agreement between PIKAN and FEM in displacement solutions. The horizontal displacement $u_x$ (Figs.~\ref{fig:pikan_ux} and \ref{fig:fem_ux}) shows the expected symmetric distribution about the neutral surface, with the upper and lower surfaces deflecting in opposite directions. For the vertical displacement $u_y$ (Figs.~\ref{fig:pikan_uy} and \ref{fig:fem_uy}), both methods capture the characteristic bending deformation pattern with maximum displacements occurring at the loaded end. The displacement magnitude $u_{\text{mag}}$ (Figs.~\ref{fig:pikan_umag} and \ref{fig:fem_umag}) provides a comprehensive view of the overall deformation.
Fig.~\ref{fig:dbc_error_analysis} quantifies the prediction accuracy through absolute error contours for both displacement components. The error analysis reveals that PIKAN maintains high accuracy throughout the domain, with maximum absolute errors typically below $6.6\times 10^{-4}$ for both $u_x$ and $u_y$ components. The error distribution shows slightly higher values near the material interfaces and the loaded boundary, which is expected due to the discontinuities in material properties and the complex stress distributions in these regions. Overall, the error magnitudes are several orders smaller than the displacement values, confirming the excellent accuracy of the PIKAN approach for multi-material structural analysis.

\begin{figure}[htbp]
    \centering
    \begin{subfigure}{0.32\textwidth}
        \centering
        \includegraphics[width=\textwidth]{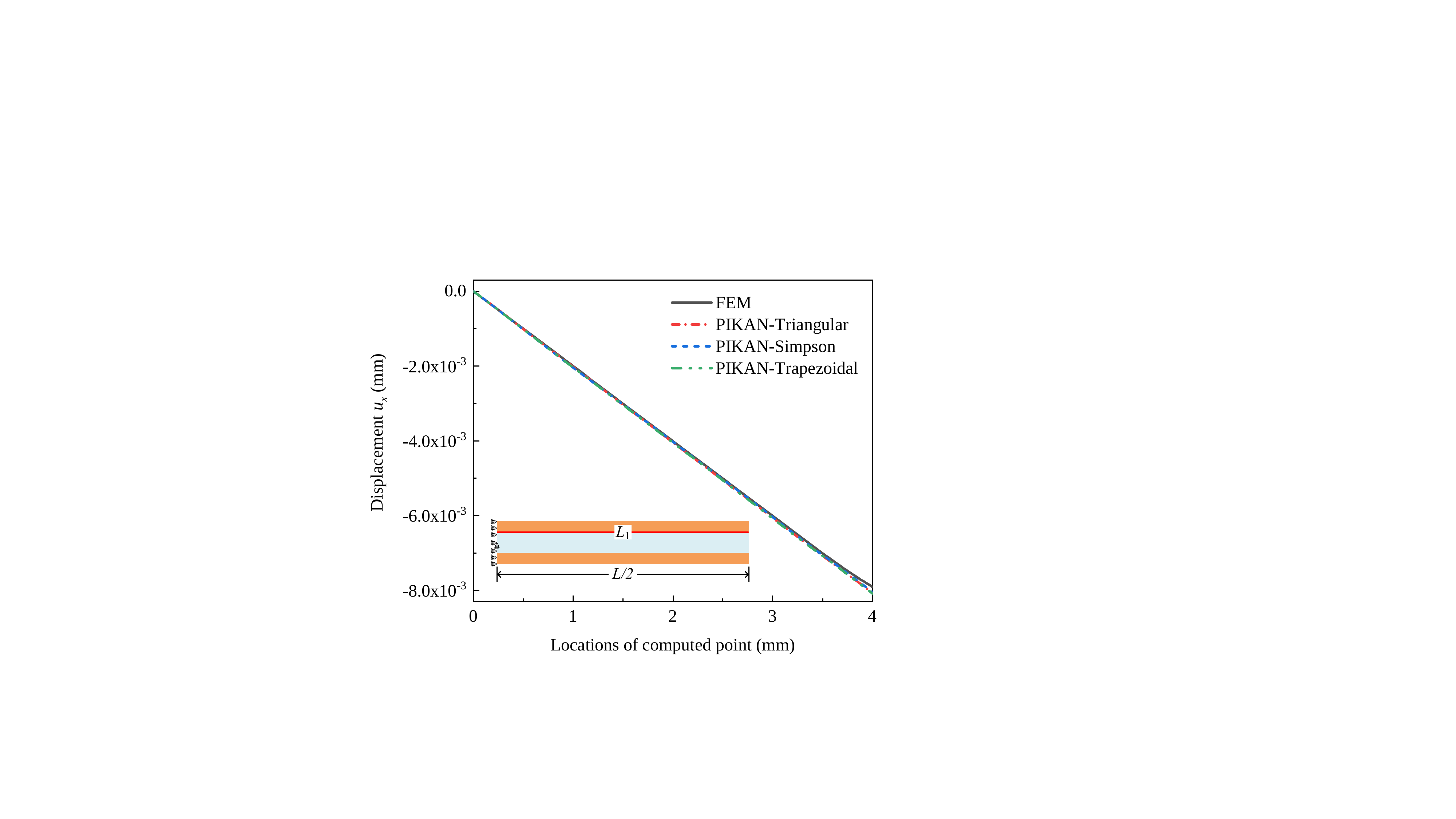}
        \caption{$u_x$ along upper interface $L_1$}
        \label{fig:interface_l1}
    \end{subfigure}
    \hfill
    \begin{subfigure}{0.32\textwidth}
        \centering
        \includegraphics[width=\textwidth]{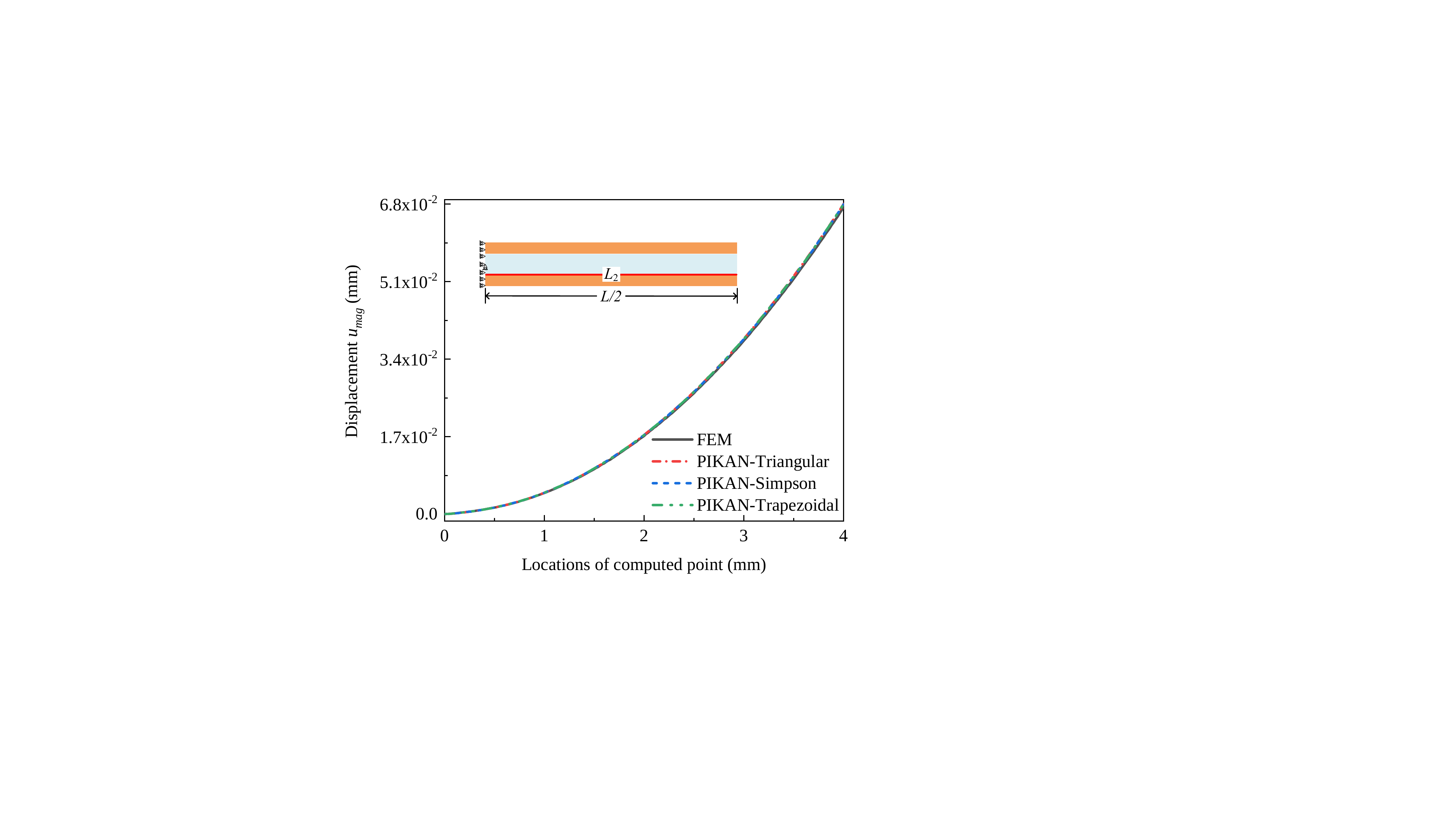}
        \caption{$u_{\text{mag}}$ along lower interface $L_2$}
        \label{fig:interface_l2}
    \end{subfigure}
    \hfill
    \begin{subfigure}{0.32\textwidth}
        \centering
        \includegraphics[width=\textwidth]{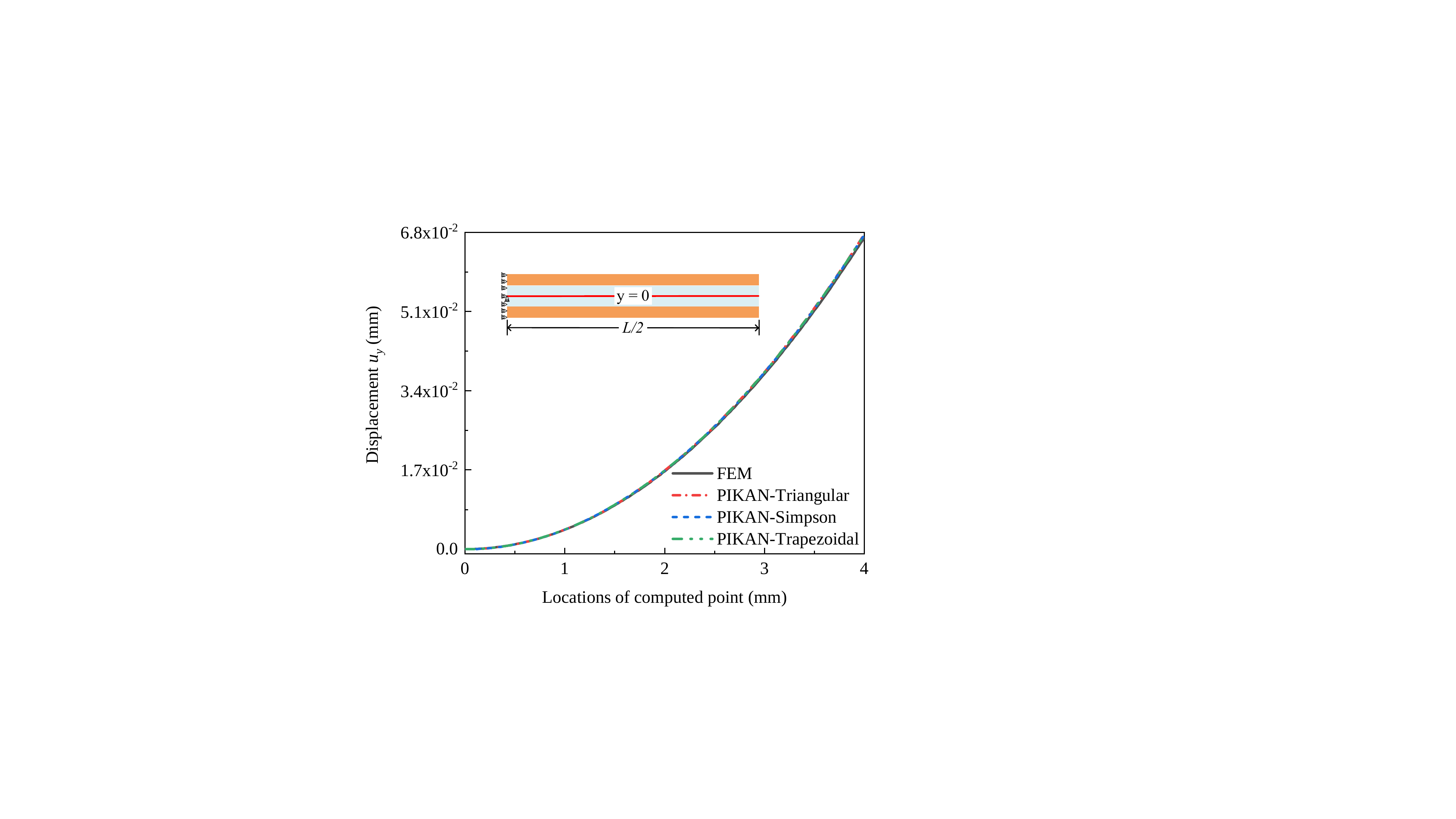}
        \caption{$u_y$ along neutral axis ($y = 0$) }
        \label{fig:line_y0}
    \end{subfigure}
    \caption{Displacement distributions at critical locations in DBC substrate.}
    \label{fig:dbc_interface_displacement}
\end{figure}

To further demonstrate the effectiveness of PIKAN for multi-material problems, we analyze displacement distributions at critical locations within the DBC substrate. The material interfaces $L_1$ and $L_2$ represent the upper and lower boundaries between the copper layers and the ceramic core, respectively. These interfaces are critical for understanding stress transfer and potential delamination sites in DBC substrates. 
Fig.~\ref{fig:dbc_interface_displacement} presents displacement distributions at different locations. Fig.~\ref{fig:dbc_interface_displacement}a shows the horizontal displacement $u_x$ along the upper material interface $L_1$ (copper-ceramic boundary), where the displacement varies linearly due to the pure bending deformation. Fig.~\ref{fig:dbc_interface_displacement}b presents the displacement magnitude $u_{\text{mag}}$ along the lower material interface $L_2$, Fig.~\ref{fig:dbc_interface_displacement}c shows the vertical displacement $u_y$ along the neutral surface (line $y = 0$, where $0 \leq x \leq 4$) within the ceramic layer.
The results demonstrate that PIKAN solutions achieve excellent agreement with FEM reference solutions across all analyzed locations. All numerical integration schemes exhibit comparable accuracy levels, indicating the robustness of the PIKAN approach. The smooth displacement distributions at material interfaces confirm that PIKAN effectively handles multi-material problems without requiring explicit interface continuity constraints.

\begin{figure}[htbp]
    \centering
    \begin{subfigure}{0.43\textwidth}
        \centering
        \includegraphics[width=\textwidth]{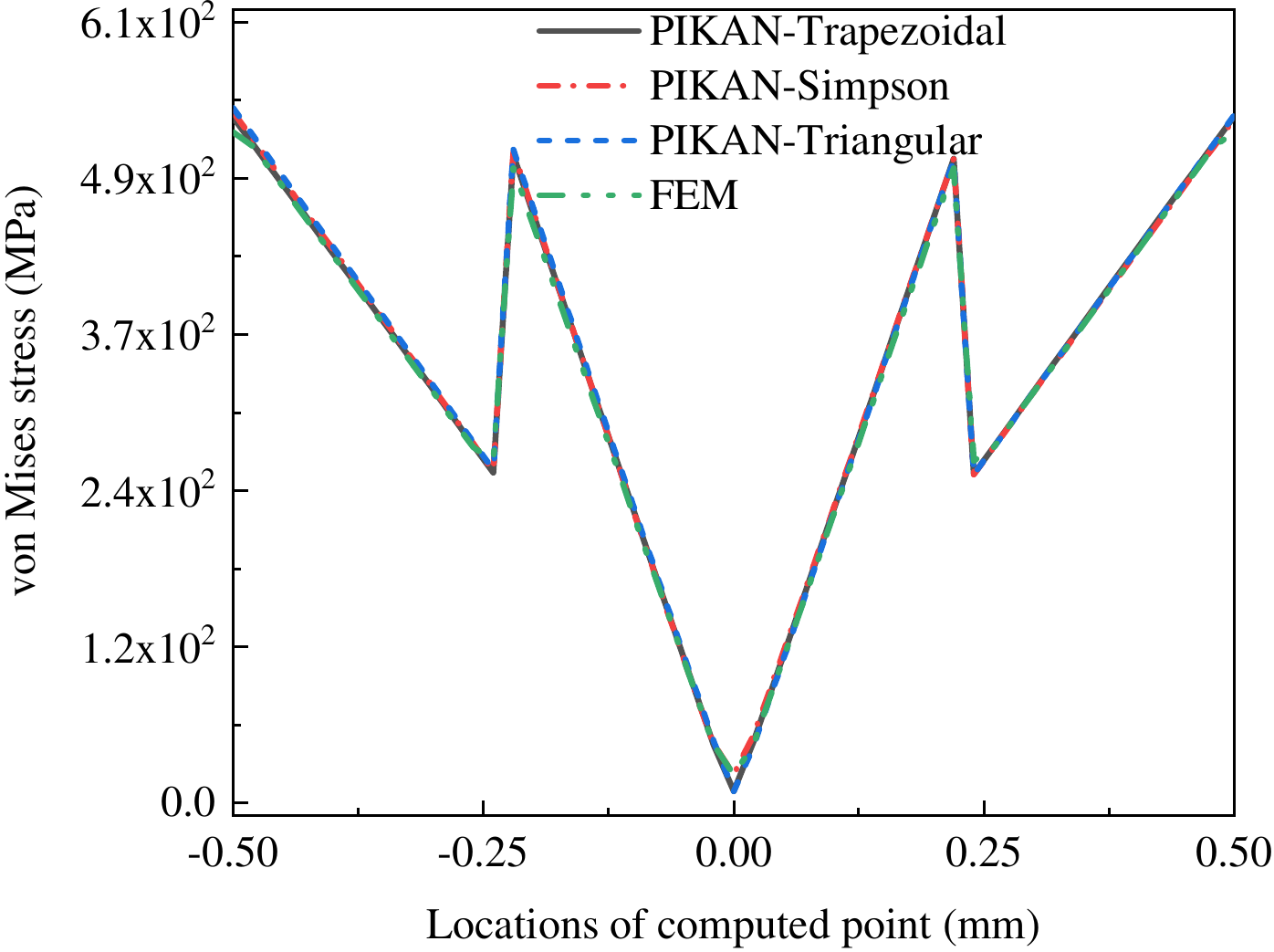}
        \caption{Stress distribution at $x = 1$.}
        \label{fig:stress_x1}
    \end{subfigure}
    \hspace{1.5em}
    \begin{subfigure}{0.43\textwidth}
        \centering
        \includegraphics[width=\textwidth]{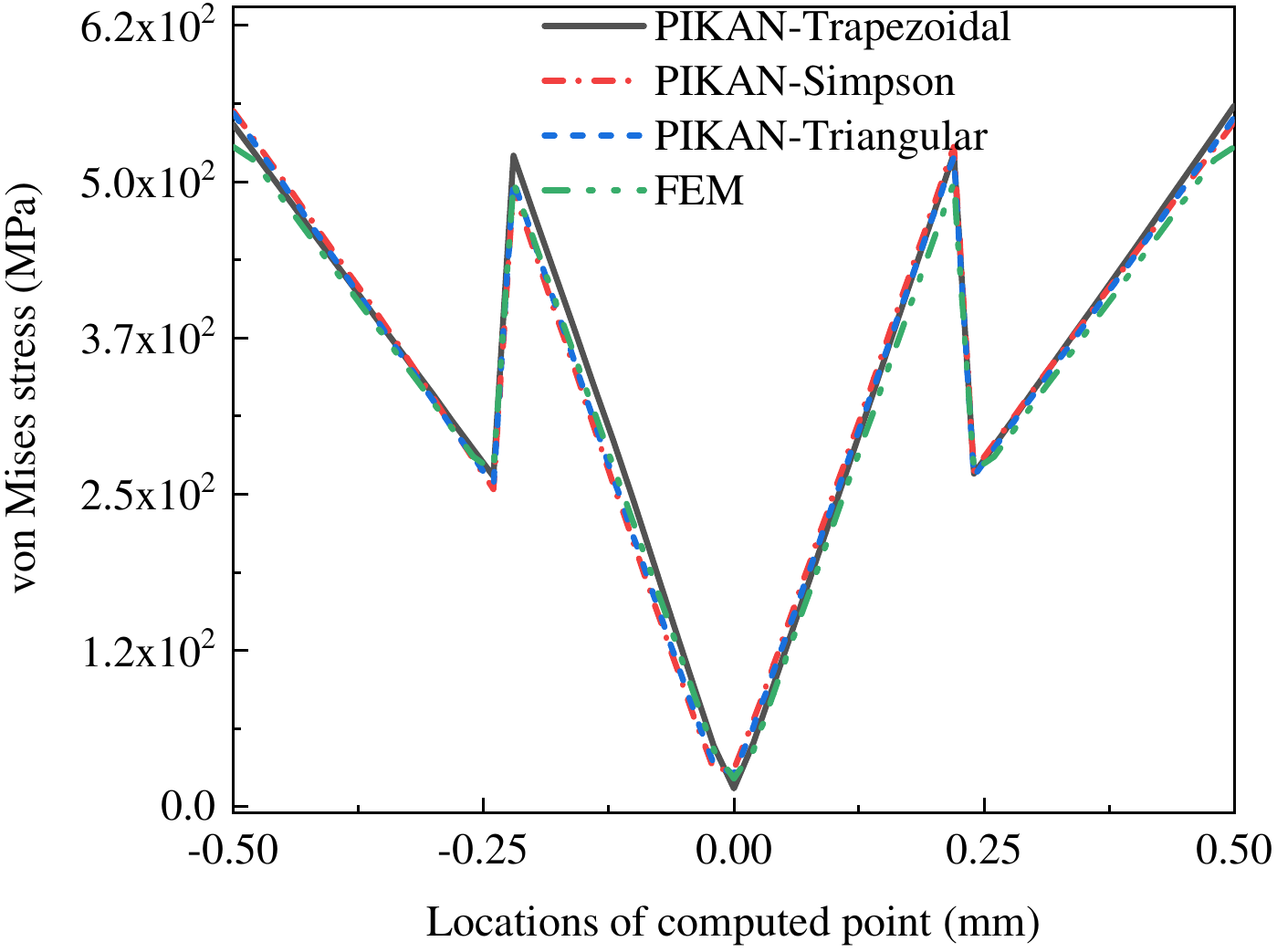}
        \caption{Stress distribution at $x = 3$.}
        \label{fig:stress_x3}
    \end{subfigure}
    \caption{Stress distributions along vertical lines in the DBC substrate.}
    \label{fig:dbc_stress_distribution}
\end{figure}

Figs.~\ref{fig:dbc_stress_distribution}a and b show stress distributions along lines $x = 1$ and $x = 3$ (where $-0.5 \leq y \leq 0.5$, as shown in Fig. \ref{fig:dbc_model}) in the DBC substrate. These lines capture stress variations in the interior region and near the loaded end.
PIKAN stress results agree well with FEM reference solutions at both locations. The stress distributions show typical multi-material behavior with stress jumps at material interfaces ($y = \pm 0.24$) due to different elastic moduli. Significant stress variations occur in the stiffer ceramic core, while copper layers show lower stress levels. Among these integration schemes, all computational results demonstrate similar accuracy and converge well to the FEM solutions. This indicates that PIKAN can effectively capture stress variations near material interfaces even without employing higher-order polynomial approximation methods (such as Simpson's rule), confirming the effectiveness and robustness of PIKAN for solving multi-material problems.

\subsection{PIKAN application to TGV-Cu structure}

\begin{figure}[htbp]
    \centering
    \begin{subfigure}{0.43\textwidth}
        \centering
        \includegraphics[width=\textwidth]{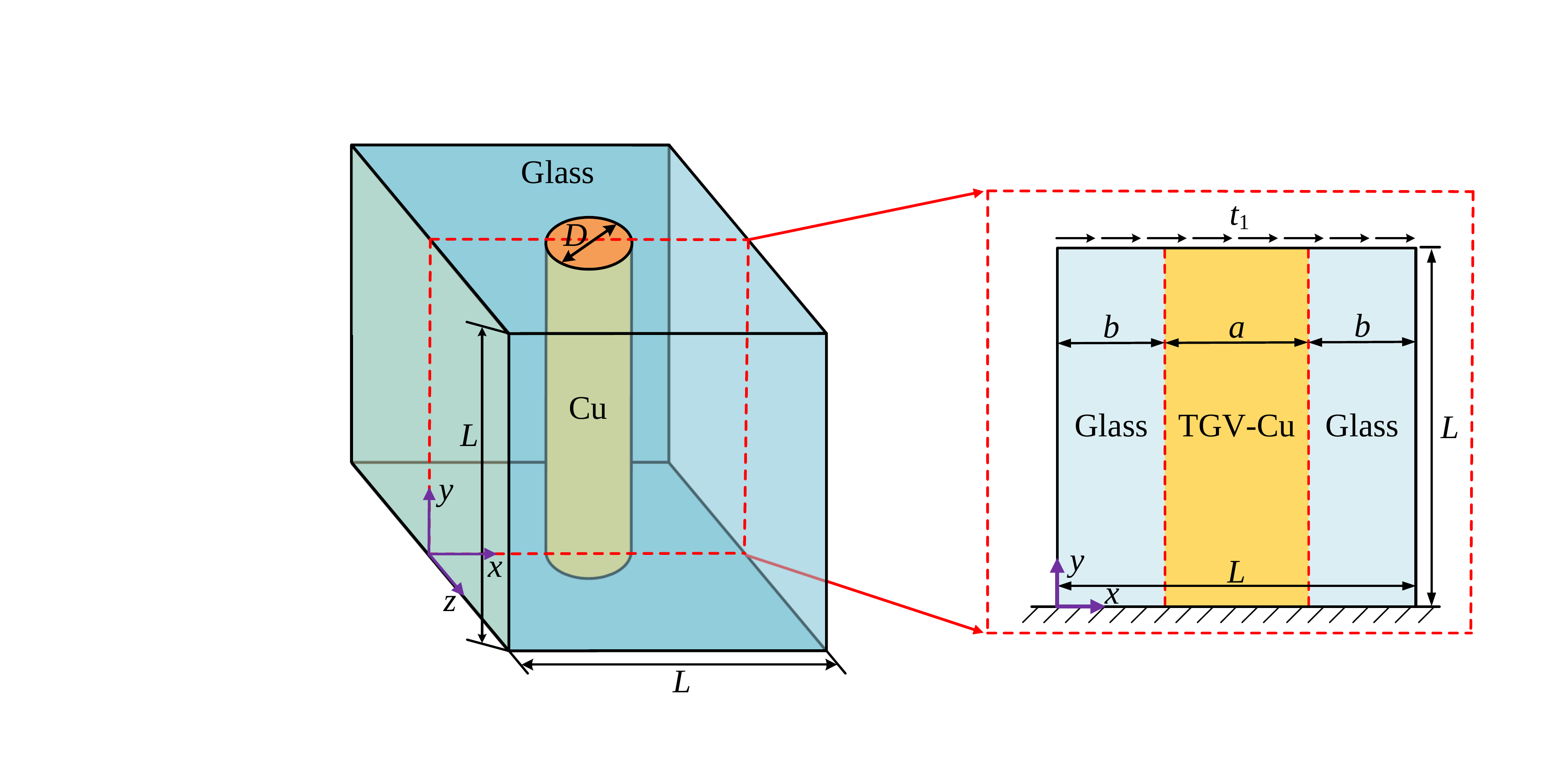}
        \caption{Single TGV-Cu structure}
        \label{fig:single_tgv}
    \end{subfigure}
    % \hfill
    \begin{subfigure}{0.4\textwidth}
        \centering
        \includegraphics[width=0.5\textwidth]{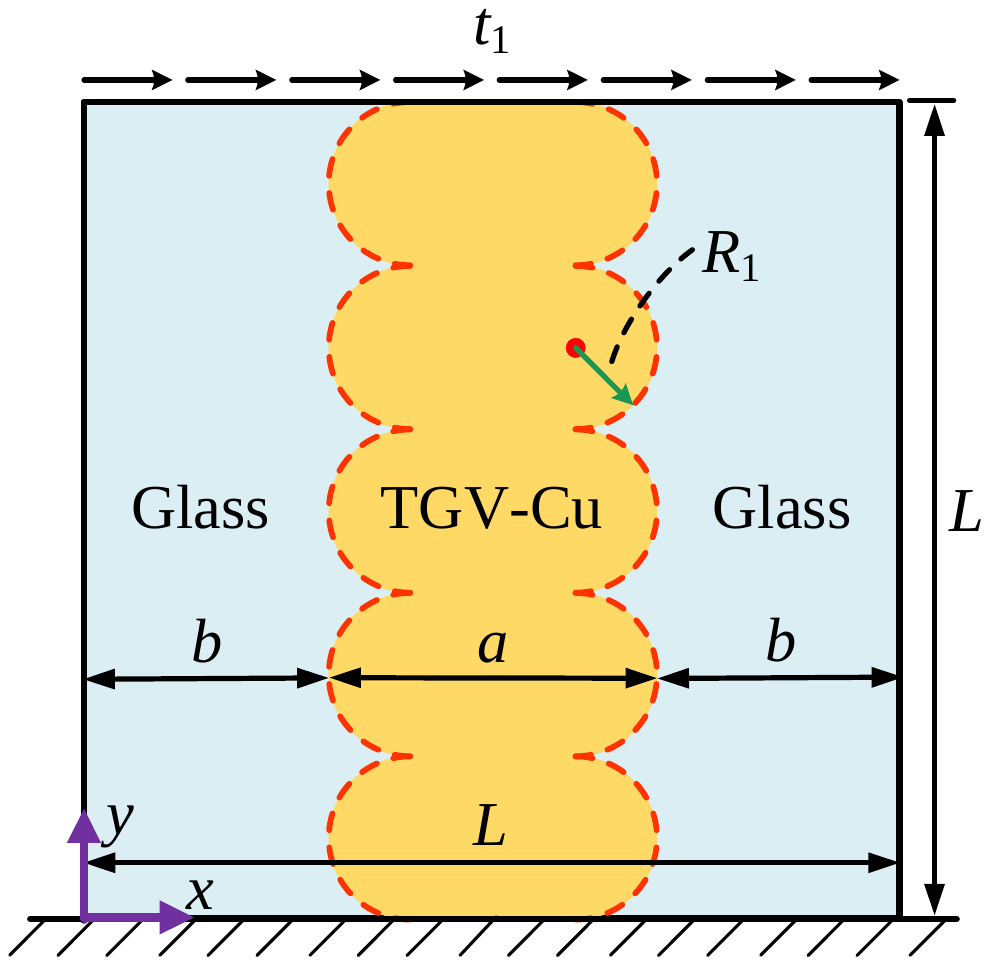}
        \caption{Interface roughness type 1 ($R_1 = 0.1$)}
        \label{fig:roughness_type1}
    \end{subfigure}
    
    \begin{subfigure}{0.4\textwidth}
        \centering
        \includegraphics[width=0.5\textwidth]{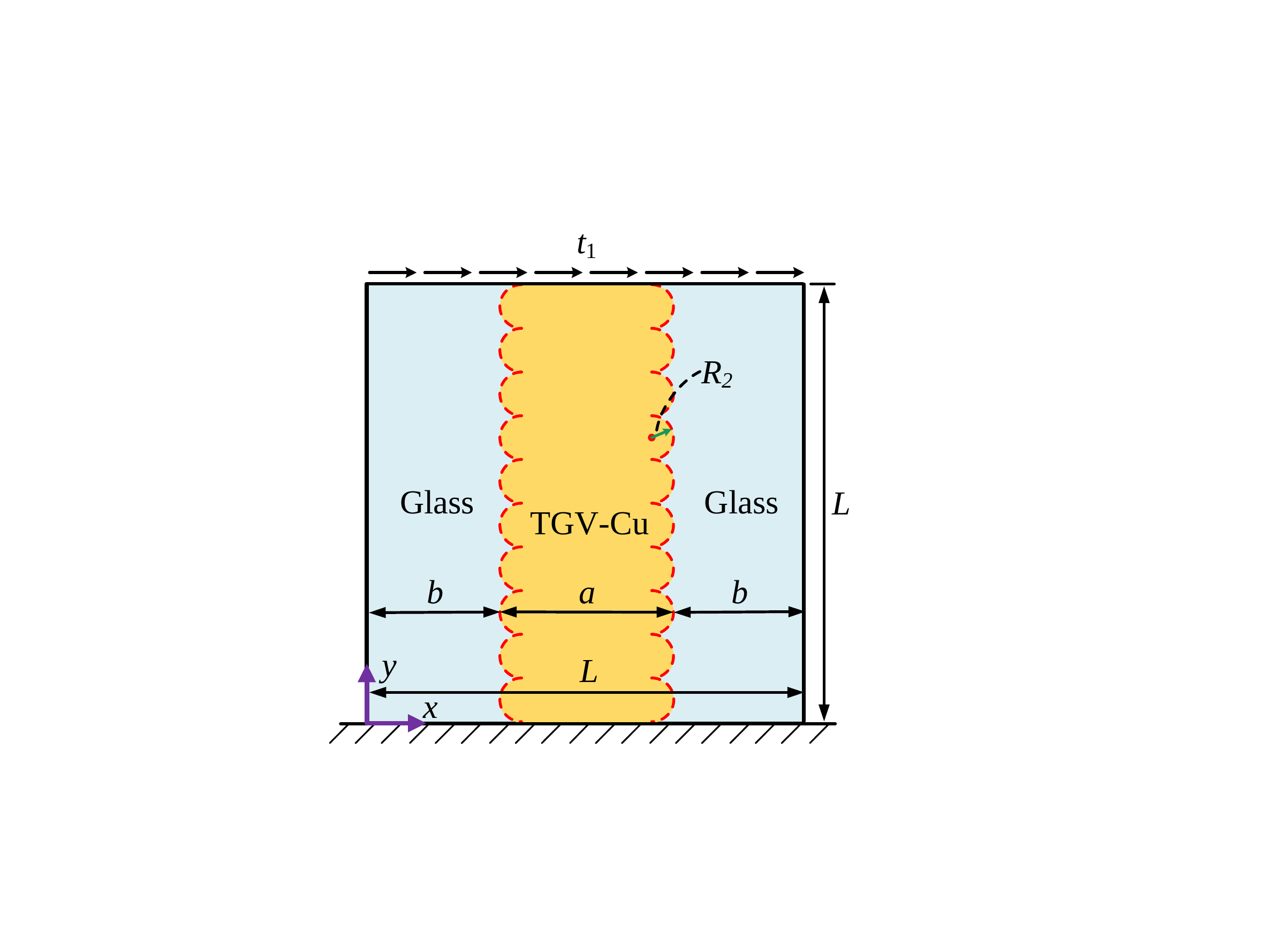}
        \caption{Interface roughness type 2 ($R_2 = 0.05$)}
        \label{fig:roughness_type2}
    \end{subfigure}
    % \hfill
    \begin{subfigure}{0.4\textwidth}
        \centering
        \includegraphics[width=0.9\textwidth]{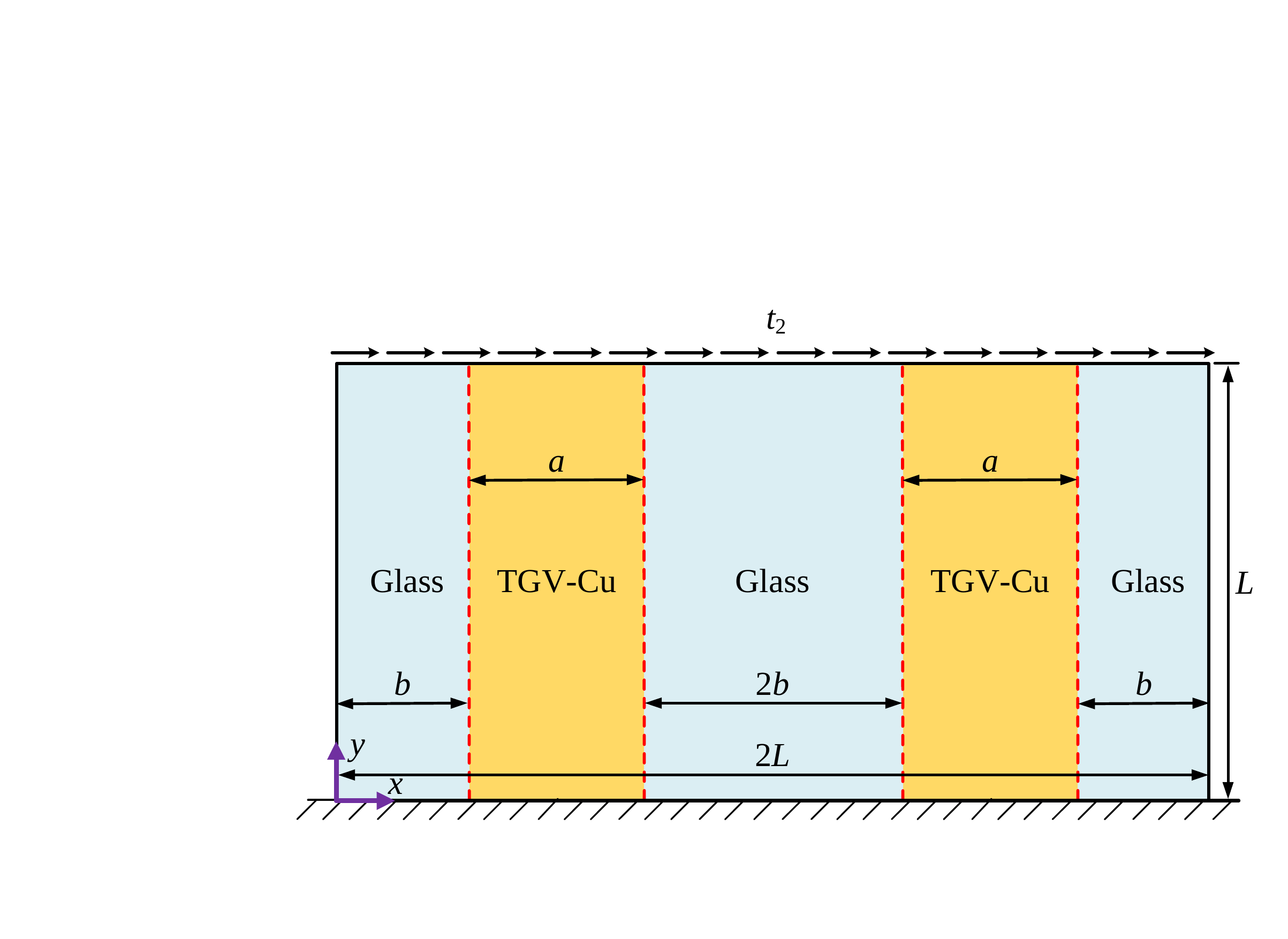}
        \caption{Dual-via TGV-Cu structure}
        \label{fig:dual_via}
    \end{subfigure}
    \caption{TGV-Cu structure models for CMP analysis: (a) baseline single via; (b,c) single via with interface roughness variations; (d) dual-via configuration.}
    \label{fig:tgv_models}
\end{figure}

Glass-substrate interposer packaging enables high-density interconnection of multi-chip systems through Through Glass Via (TGV) technology. TGV structures provide mechanical support and electrical interconnection via electroplated copper pillars fabricated through electroplating processes. After electroplating, residual copper layers on the glass surface must be removed through chemical mechanical polishing (CMP) to achieve smooth interfaces and proper electrical isolation~\cite{Jeffrey2015tgv}.
PIKAN is employed to analyze TGV-Cu structures during the CMP process, as shown in Fig.~\ref{fig:tgv_models}. Four model configurations are investigated: (1) single TGV-Cu structure (Fig.~\ref{fig:tgv_models}a); (2) and (3) TGV-Cu structures with different interface roughness patterns (Figs.~\ref{fig:tgv_models}b and c); and (4) dual-via TGV-Cu structure (Fig.~\ref{fig:tgv_models}d).

The geometric dimensions are $L = 1$, $a = 0.4$, and $b = 0.3$. Interface roughness is modeled as semicircular protrusions with radii $R_1 = 0.1$ and $R_2 = 0.05$ at the glass-copper interface. The dual-via configuration features two copper vias with center-to-center spacing of 0.6.
Boundary conditions simulate the CMP process: the bottom surface is fixed, while horizontal tractions is applied at the top surface. The applied tractions are $t_1 = 200$ N/mm for single-via structures and $t_2 = 250$ N/mm for the dual-via structure. Material properties are listed in Table~\ref{tab:tgv_materials}.

\begin{table}[htbp]
    \scriptsize
    \centering
    \caption{Material parameters for TGV-Cu models}
    \label{tab:tgv_materials}
    \begin{tabular}{@{}ccc@{}}
    \toprule
    Material & Young's Modulus $E$ (MPa) & Poisson's Ratio $\nu$ \\
    \midrule
    Glass substrate& 77,000 & 0.24 \\
    Cu (electroplated) & 150,000 & 0.3 \\
    \bottomrule
    \end{tabular}
\end{table}

\begin{figure}[htbp]
  \centering
  % ---------- 第一行：4 张子图 ----------
  \begin{subfigure}[t]{0.2\textwidth}
    \centering
    \includegraphics[width=0.7\linewidth]{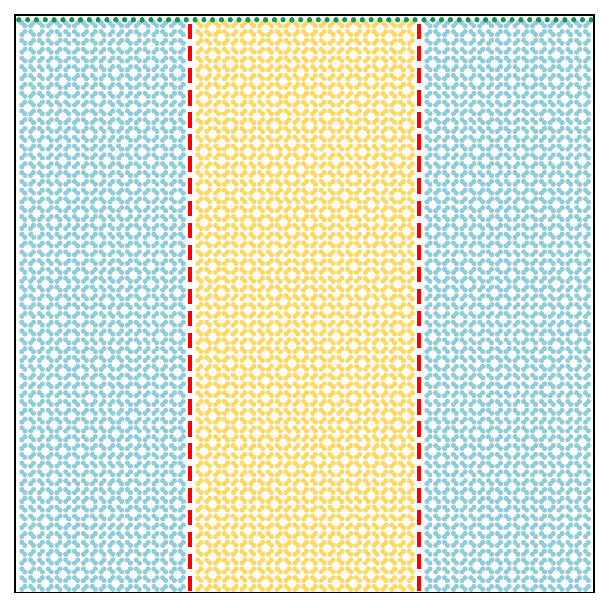}
    \caption{Single via: triangular}
    \label{fig:single_triangular}
  \end{subfigure}\hspace{0.01em}
  \begin{subfigure}[t]{0.2\textwidth}
    \centering
    \includegraphics[width=0.7\linewidth]{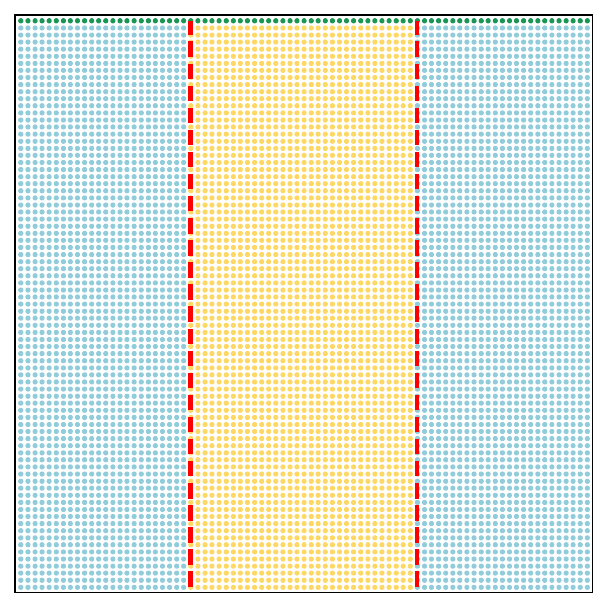}
    \caption{Single via: uniform}
    \label{fig:single_uniform}
  \end{subfigure}\hspace{0.01em}
  \begin{subfigure}[t]{0.2\textwidth}
    \centering
    \includegraphics[width=0.7\linewidth]{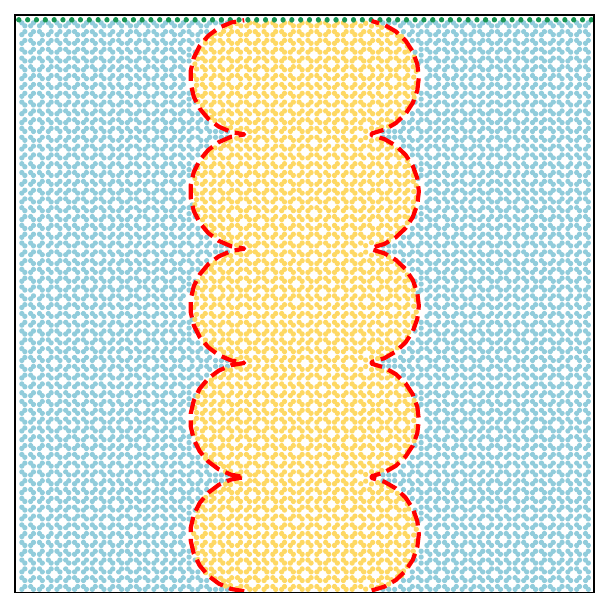}
    \caption{Rough-1: triangular}
    \label{fig:roughness1_triangular}
  \end{subfigure}\hspace{0.01em}
  \begin{subfigure}[t]{0.2\textwidth}
    \centering
    \includegraphics[width=0.7\linewidth]{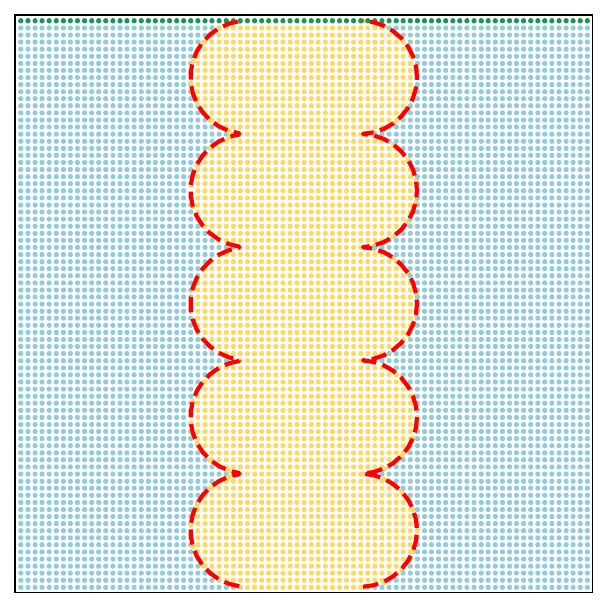}
    \caption{Rough-1: uniform}
    \label{fig:roughness1_uniform}
  \end{subfigure}

  \par\medskip   % 强制换行 + 垂直间距
  % ---------- 第二行：4 张子图 ----------
  \begin{subfigure}[t]{0.2\textwidth}
    \centering
    \includegraphics[width=0.7\linewidth]{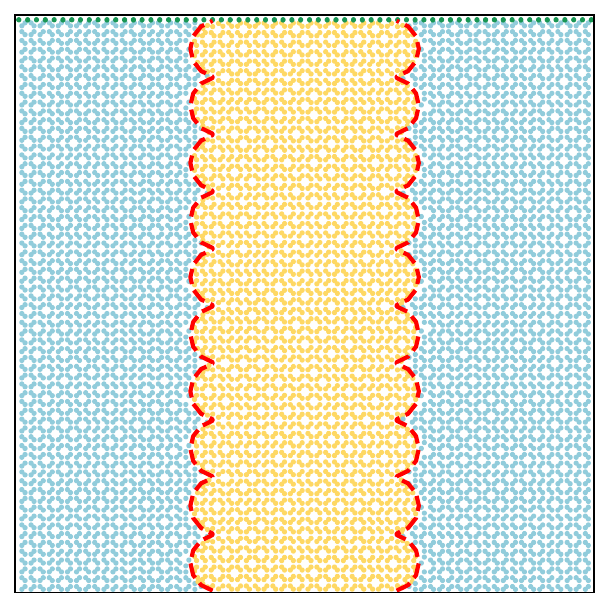}
    \caption{Rough-2: triangular}
    \label{fig:roughness2_triangular}
  \end{subfigure}\hspace{0.01em}
  \begin{subfigure}[t]{0.2\textwidth}
    \centering
    \includegraphics[width=0.7\linewidth]{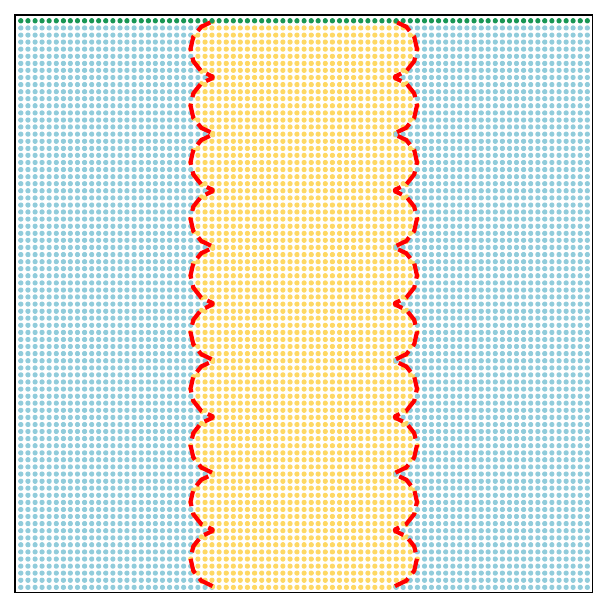}
    \caption{Rough-2: uniform}
    \label{fig:roughness2_uniform}
  \end{subfigure}\hspace{0.01em}
  \begin{subfigure}[t]{0.29\textwidth}
    \centering
    \includegraphics[width=0.7\linewidth]{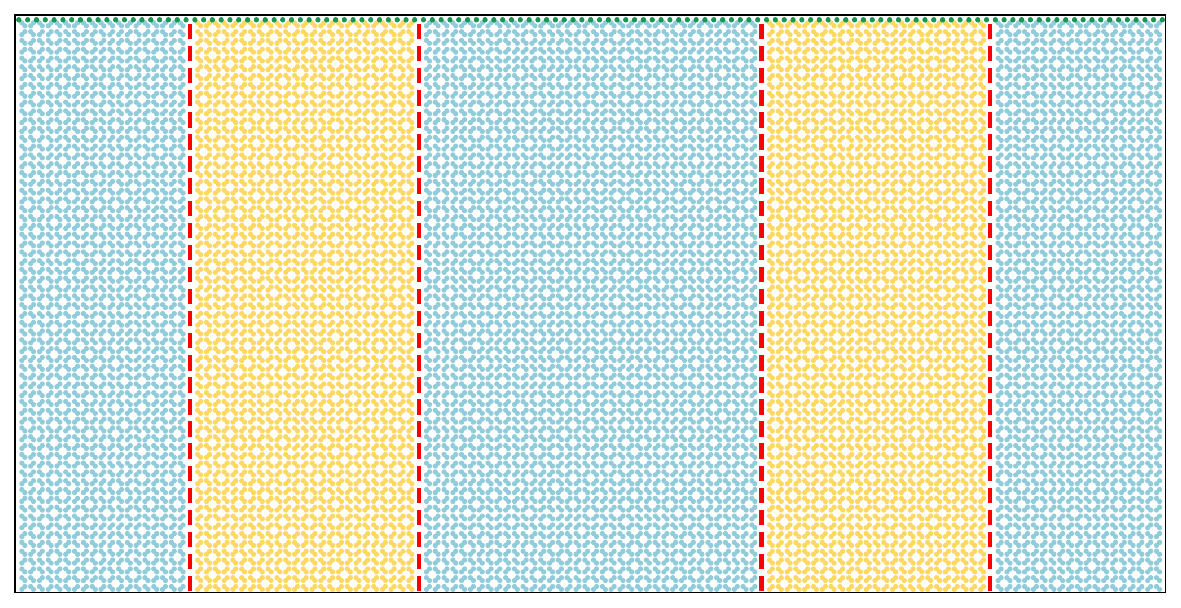}
    \caption{Dual-via: triangular}
    \label{fig:dual_triangular}
  \end{subfigure}\hspace{0.01em}
  \begin{subfigure}[t]{0.29\textwidth}
    \centering
    \includegraphics[width=0.7\linewidth]{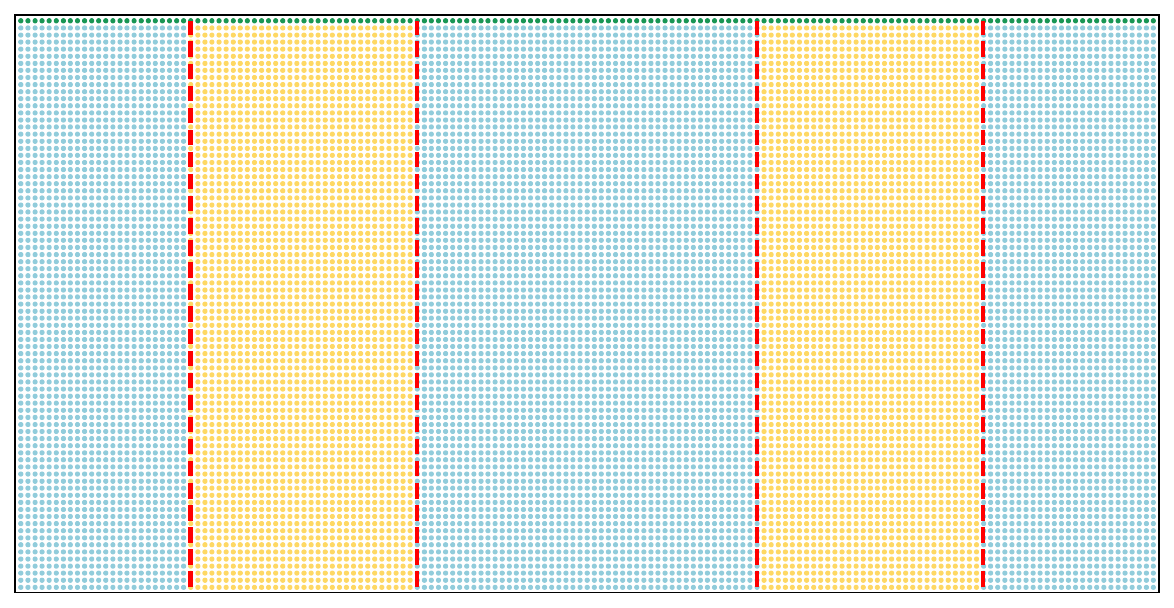}
    \caption{Dual-via: uniform}
    \label{fig:dual_uniform}
  \end{subfigure}

  \caption{Training point distribution schemes for TGV-Cu structure models: different sampling strategies applied to single and dual via configurations (Rough-1 and Rough-2 denote roughness 1 and roughness 2, respectively).}
  \label{fig:tgv_sample_distribution}
\end{figure}

The KAN architecture for this problem maintains consistency with previous examples: $[2,5,5,5,2]$ structure with grid size 20, B-spline order 3, and 1,750 trainable parameters. Fig.~\ref{fig:tgv_sample_distribution} shows the training sample point distribution schemes for all TGV-Cu structure models.
For centroid point sampling, the domain point counts are: single structure - 79,200 points; roughness type 1 - 80,094 points; roughness type 2 - 80,064 points; and dual-via structure - 158,400 points. For uniform distribution sampling patterns, the distributions are: single structure - 40,401 points; both roughness models - 40,761 points each; and dual-via structure - 80,601 points.
Natural boundary condition training points are distributed as 201 points for single TGV-Cu structures and 401 points for the dual-via structure, corresponding to the different boundary lengths requiring load application.

\begin{figure}[htbp]
  \centering
  \begin{subfigure}[t]{0.23\textwidth}
    \centering
    \includegraphics[width=\linewidth]{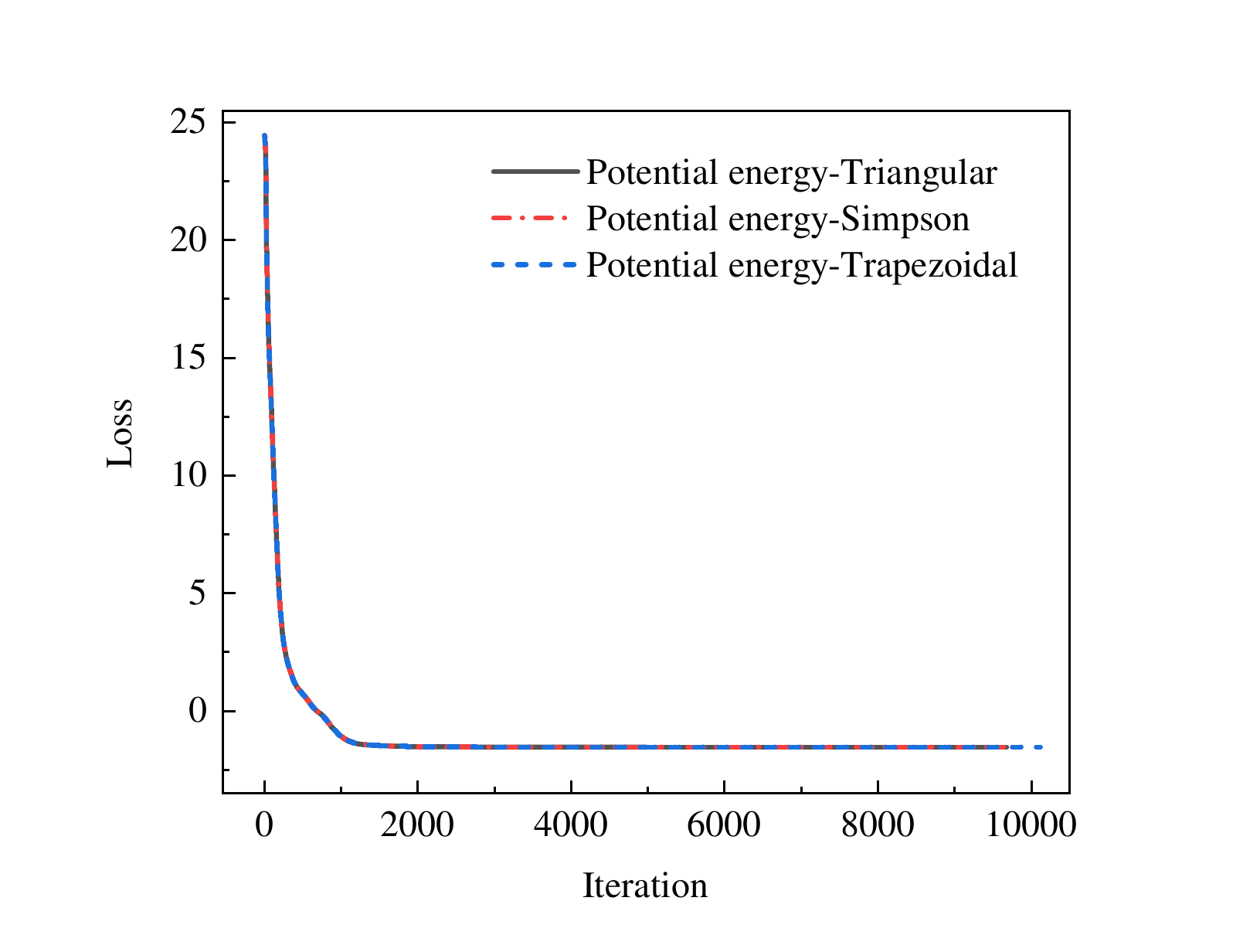}
    \caption{Single TGV-Cu}
    \label{fig:single_energy}
  \end{subfigure}\hfill
  \begin{subfigure}[t]{0.24\textwidth}
    \centering
    \includegraphics[width=\linewidth]{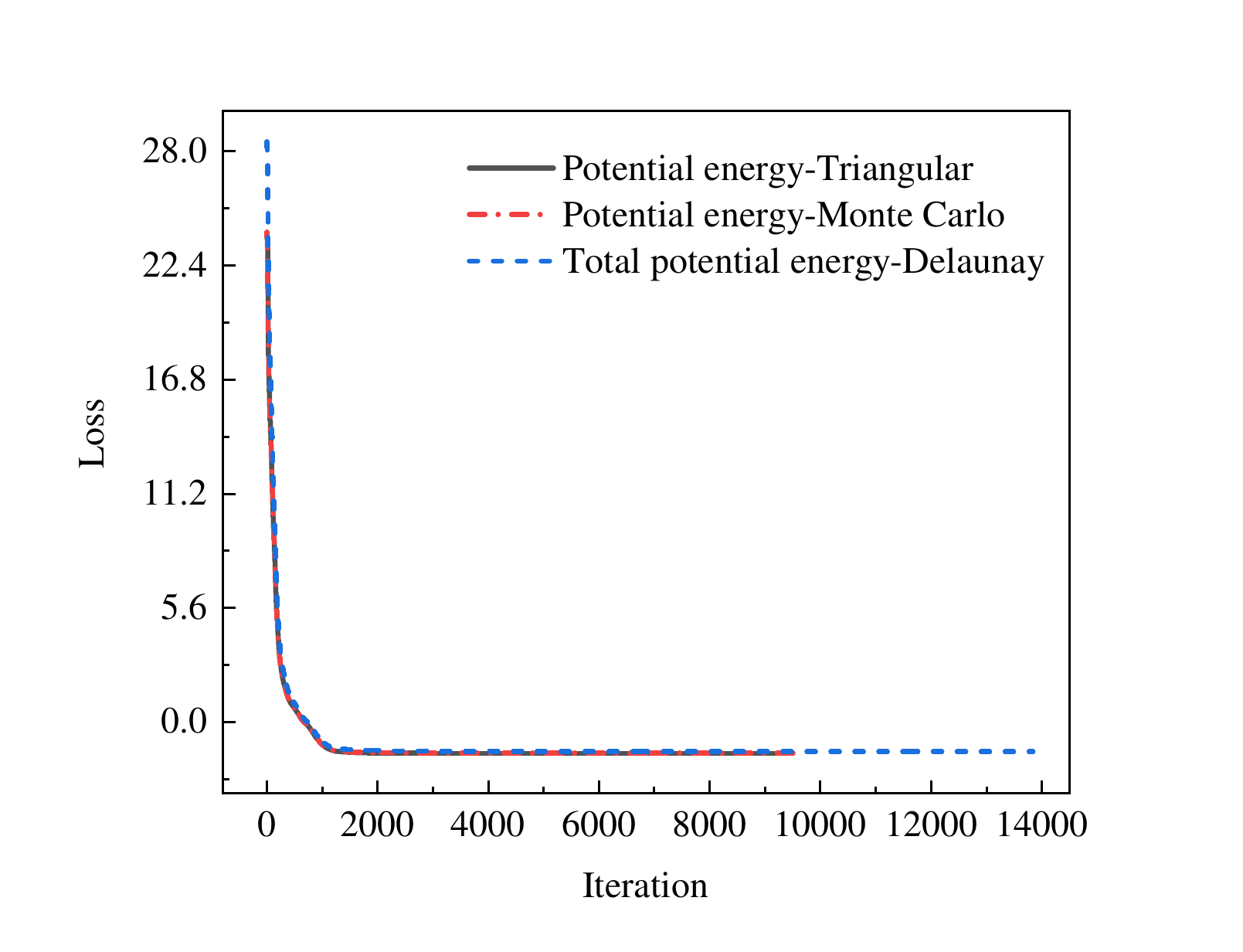}
    \caption{Interface rough-1}
    \label{fig:roughness1_energy}
  \end{subfigure}\hfill
  \begin{subfigure}[t]{0.235\textwidth}
    \centering
    \includegraphics[width=\linewidth]{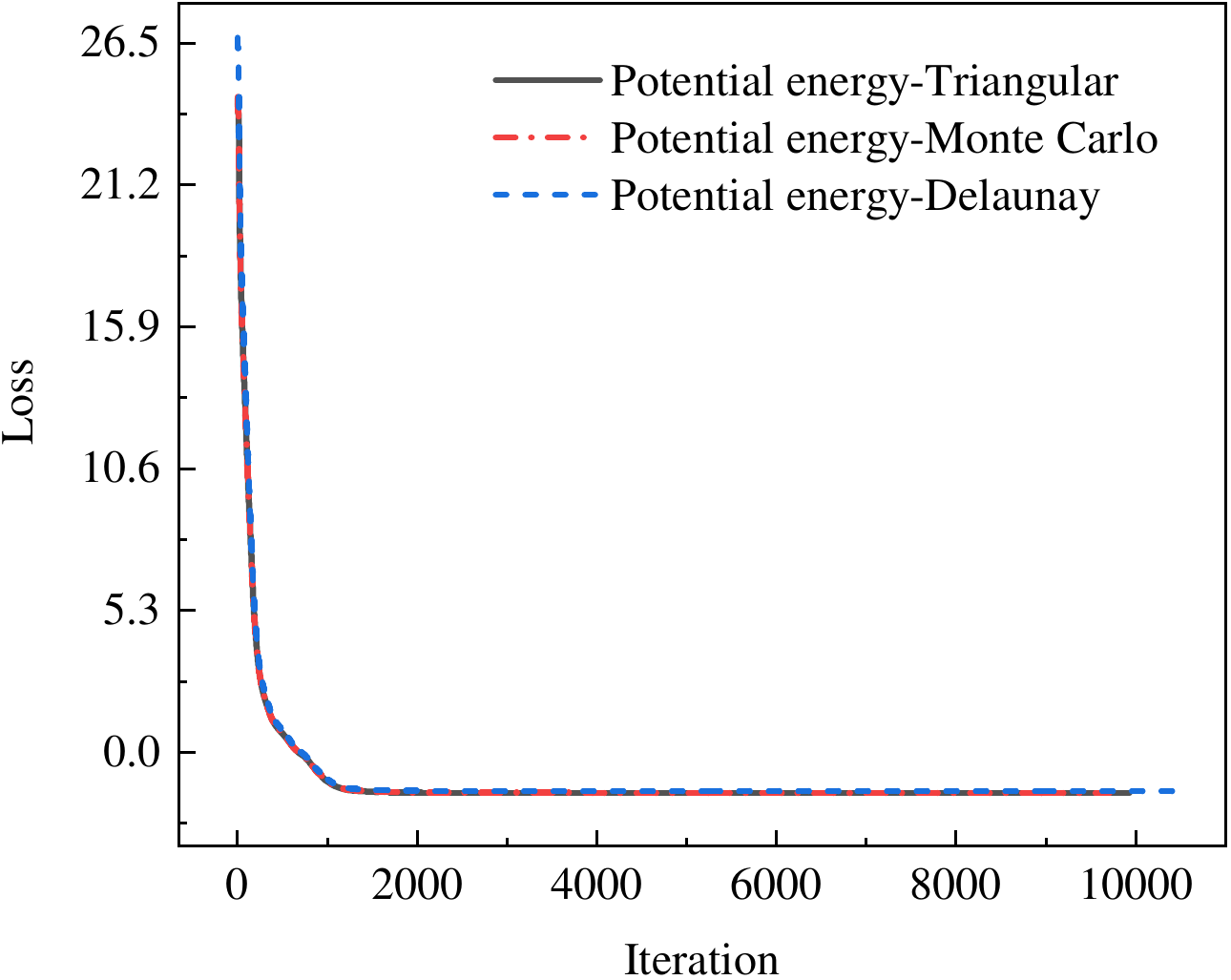}
    \caption{Interface rough-2}
    \label{fig:roughness2_energy}
  \end{subfigure}\hfill
  \begin{subfigure}[t]{0.245\textwidth}
    \centering
    \includegraphics[width=\linewidth]{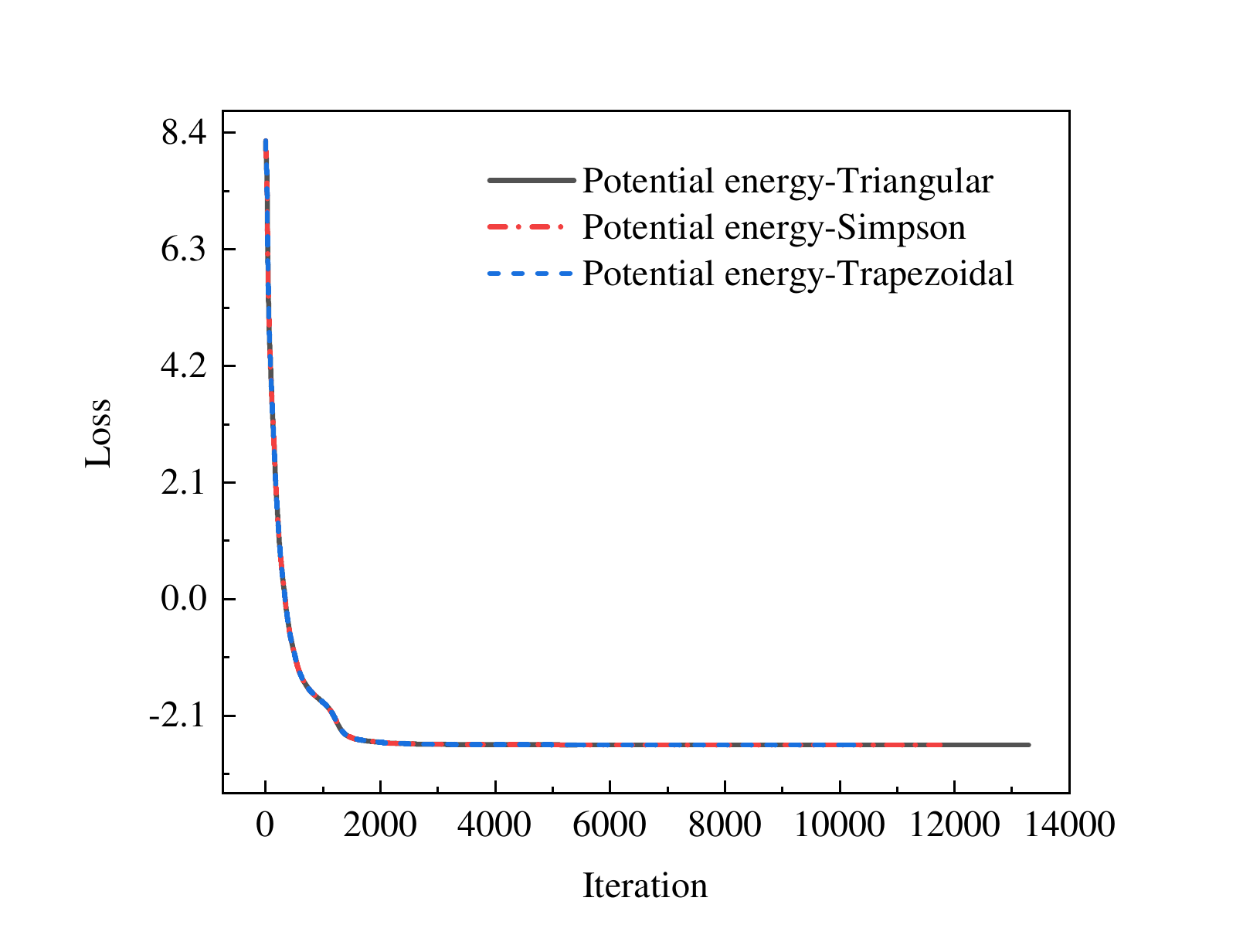}
    \caption{Dual-via}
    \label{fig:dual_energy}
  \end{subfigure}
  \caption{Loss function evolution for different TGV-Cu structure configurations.}
  \label{fig:tgv_energy_convergence}
\end{figure}

Fig.~\ref{fig:tgv_energy_convergence} shows the evolution of loss functions for different TGV-Cu structure models during neural network training. Despite varying geometric complexities and structural arrangements, all models exhibit similar convergence behavior with consistent energy reduction trends.
The convergence patterns demonstrate PIKAN's robustness across different geometric configurations. All integration schemes achieve stable convergence, with triangular integration typically showing slightly faster convergence. The decreasing loss function indicates the system approaching an equilibrium state. When the loss function stabilizes after a certain number of iterations (approximately 2000 iterations in this case), it indicates the system is approaching a stationary point of the potential energy functional. In practical applications, an early stopping criterion can be adopted: training can be terminated when the relative change in loss function over consecutive iterations falls below a preset tolerance $\varepsilon$, thereby balancing accuracy and computational efficiency. In this work, to demonstrate the complete convergence process and ensure result stability, we chose to continue iterations to a sufficiently conservative number.

    % Most optimized presentation - separate figures for each comparison
    \begin{figure}[htbp]
    \centering
    \begin{subfigure}{0.48\textwidth}
        \centering
        \includegraphics[width=\textwidth]{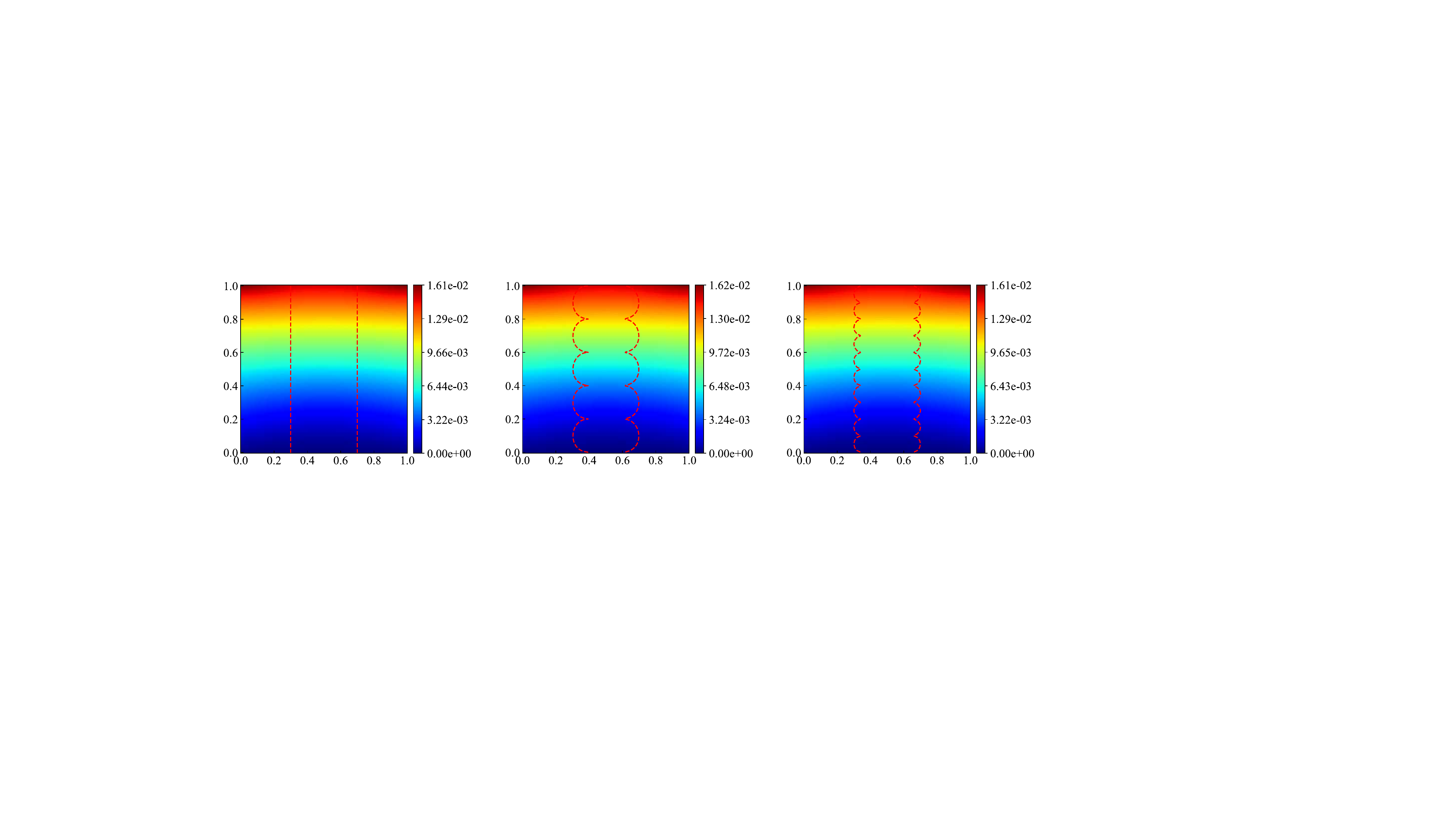}
        \caption{PIKAN: $u_x$}
    \end{subfigure}
    \hfill
    \begin{subfigure}{0.48\textwidth}
        \centering
        \includegraphics[width=\textwidth]{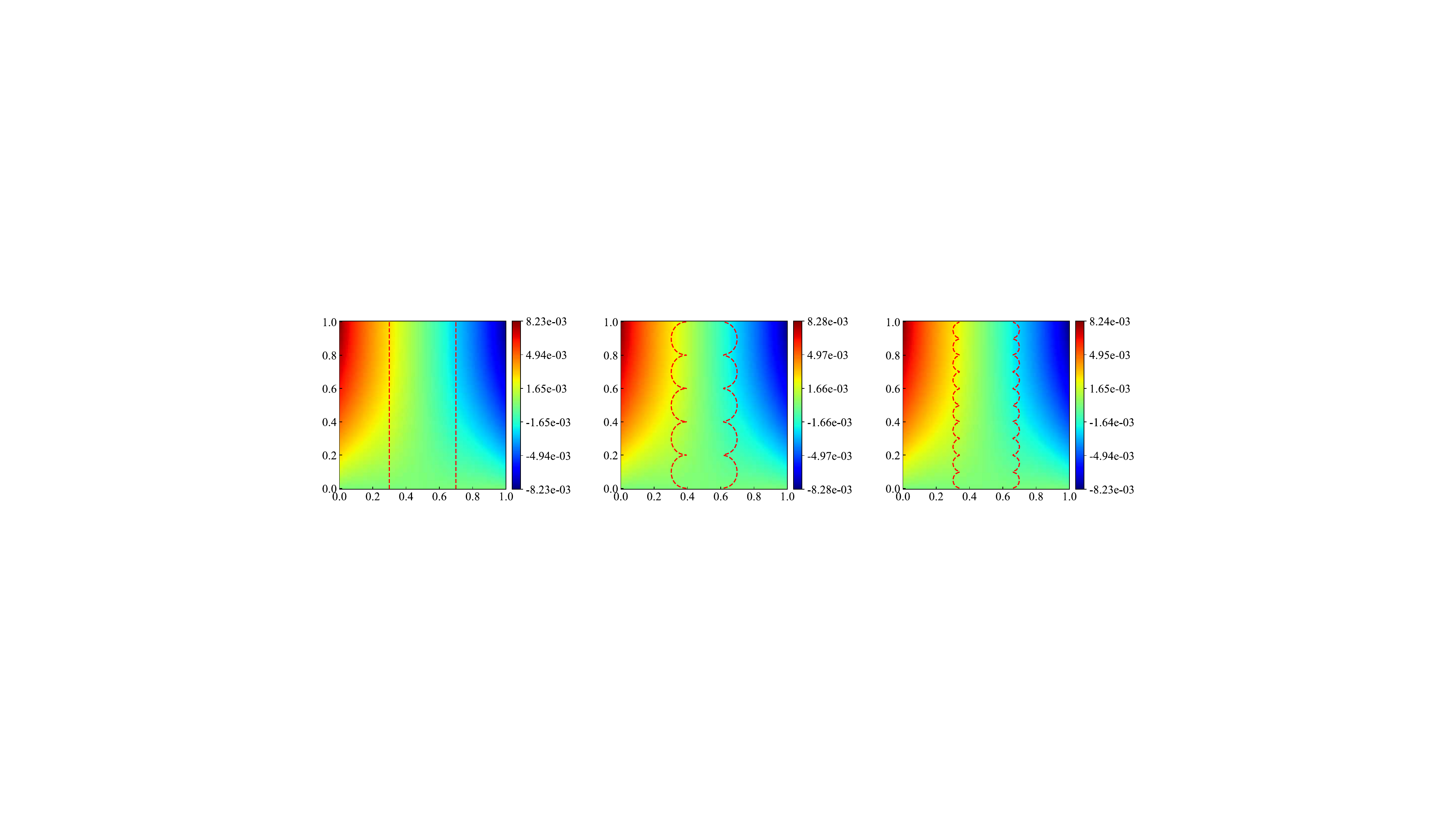}
        \caption{PIKAN: $u_y$}
    \end{subfigure}

    \begin{subfigure}{0.48\textwidth}
        \centering
        \includegraphics[width=\textwidth]{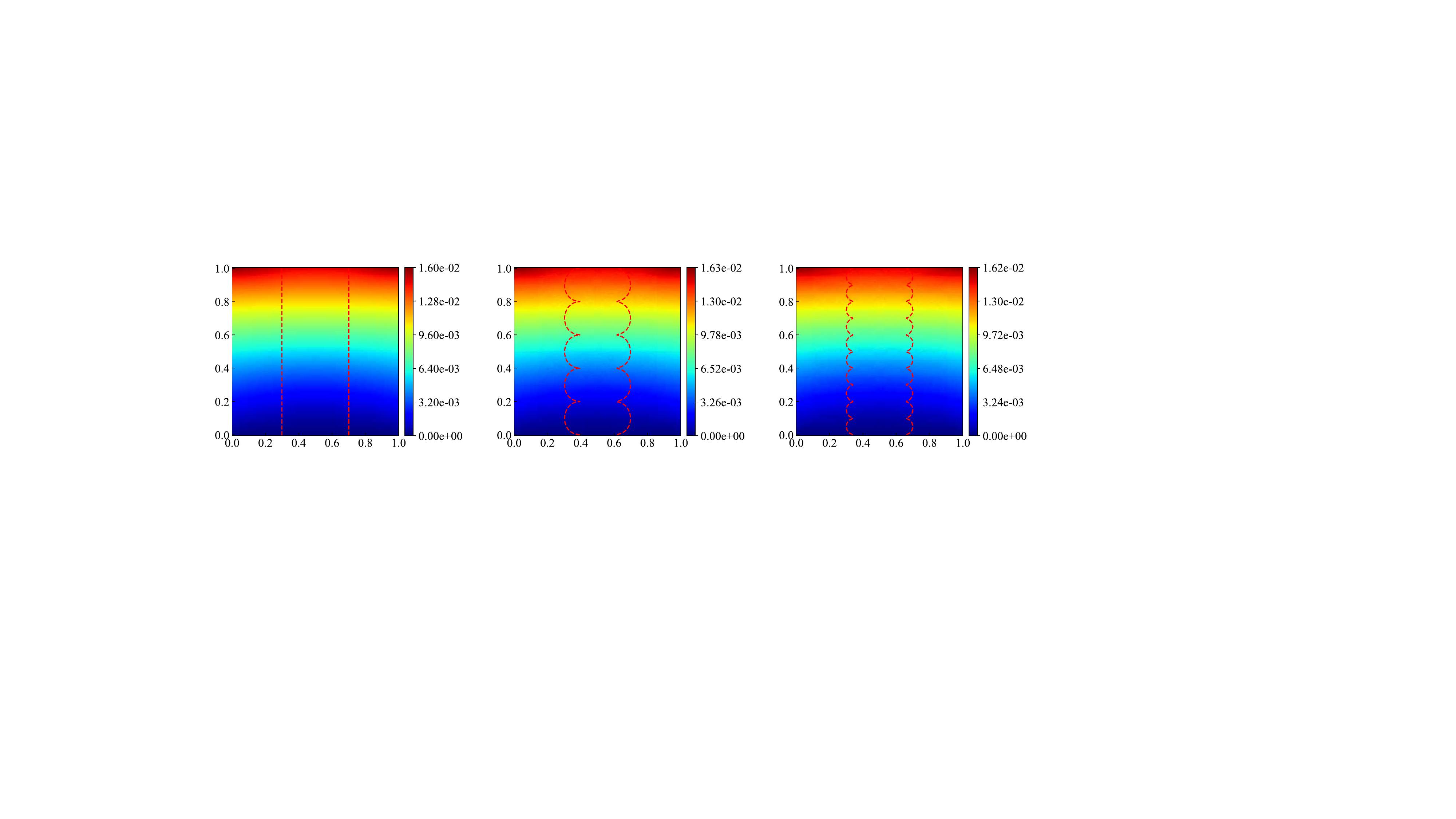}
        \caption{FEM: $u_x$}
    \end{subfigure}
     \hfill
    \begin{subfigure}{0.48\textwidth}
        \centering
        \includegraphics[width=\textwidth]{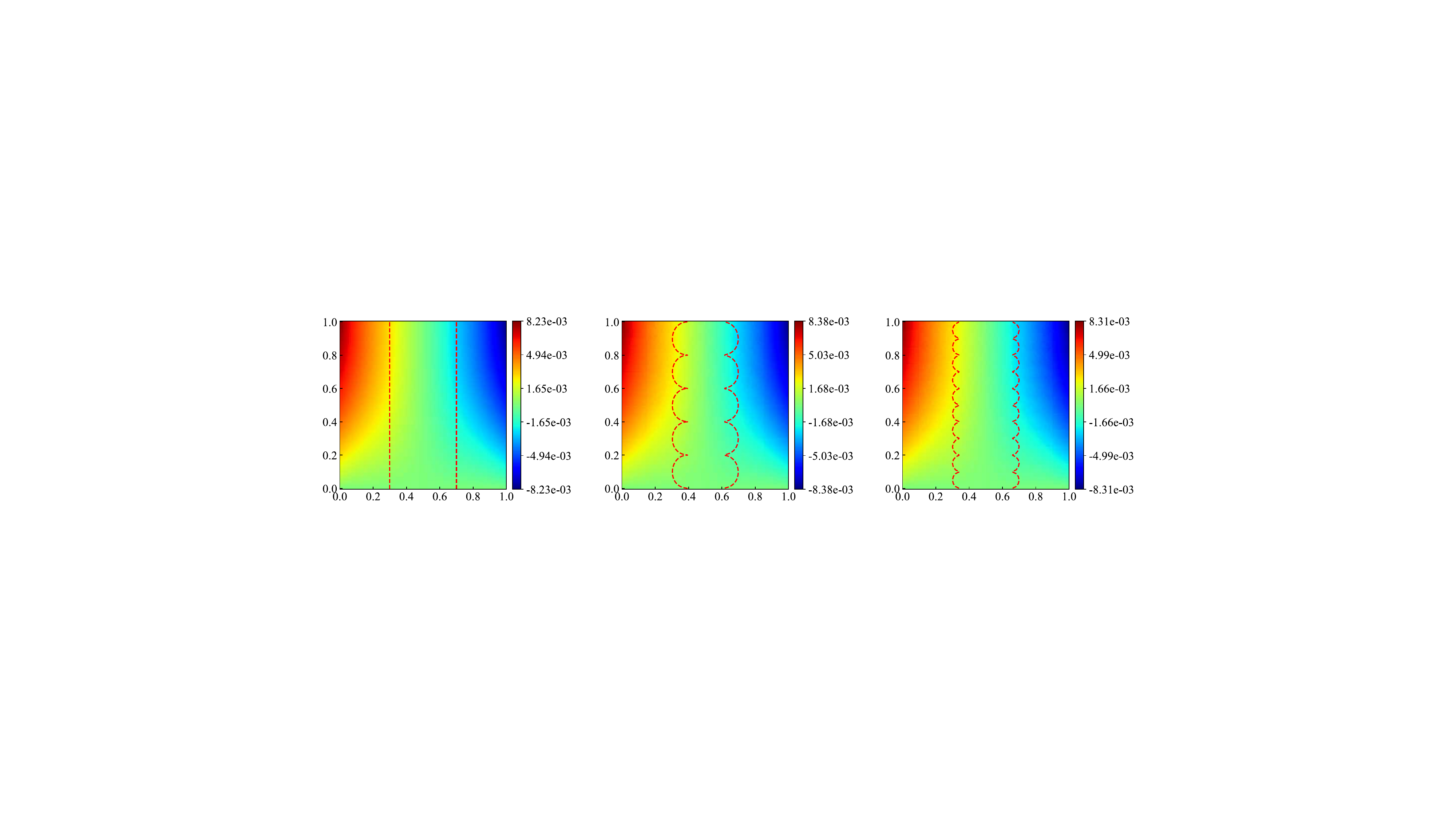}
        \caption{FEM: $u_y$}
    \end{subfigure}
    
    \begin{subfigure}{0.48\textwidth}
        \centering
        \includegraphics[width=\textwidth]{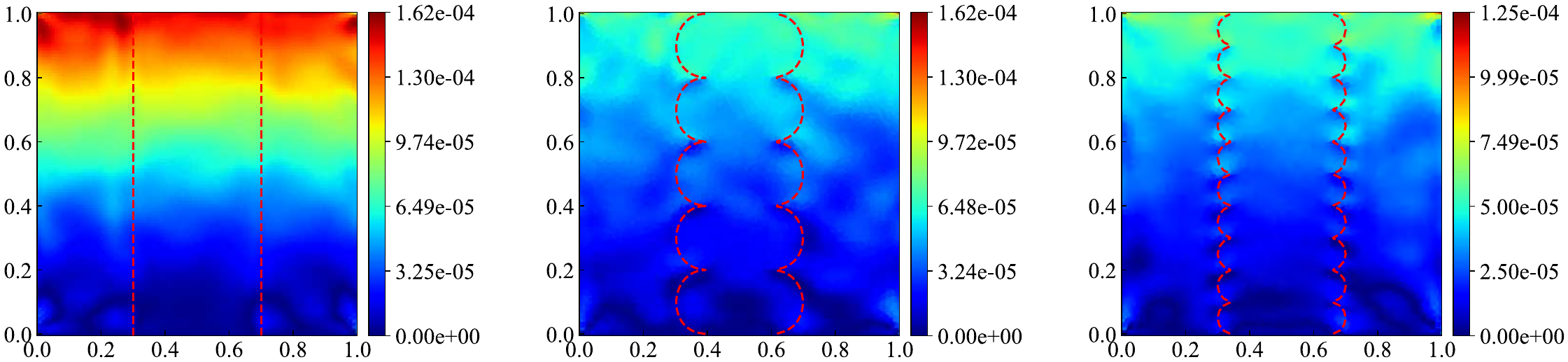}
        \caption{Absolute error: $u_x$}
    \end{subfigure}
      \hfill
    \begin{subfigure}{0.48\textwidth}
        \centering
        \includegraphics[width=\textwidth]{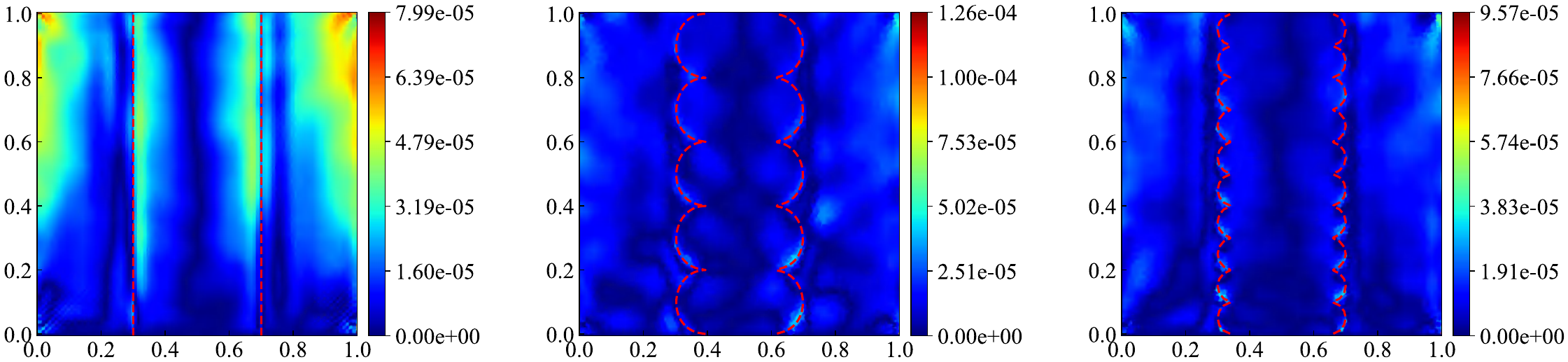}
        \caption{Absolute error: $u_y
        $}
    \end{subfigure}
    \caption{Comparison of PIKAN and FEM displacement fields and the corresponding  absolute error in single-via TGV-Cu structures.}
    \label{fig:tgv_smooth_comparison}
    \end{figure}

\begin{figure}[htbp]
  \centering
  % ---------- 第一行 ----------
  \begin{subfigure}{0.3\textwidth}
    \centering
    \includegraphics[width=0.95\linewidth]{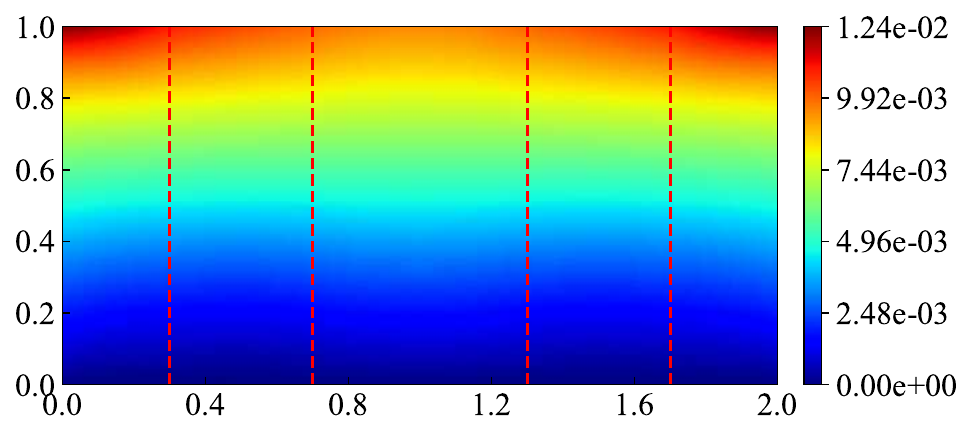}
    \caption{PIKAN: $u_x$}
  \end{subfigure}\hspace{0.1em}
  \begin{subfigure}{0.3\textwidth}
    \centering
    \includegraphics[width=0.95\linewidth]{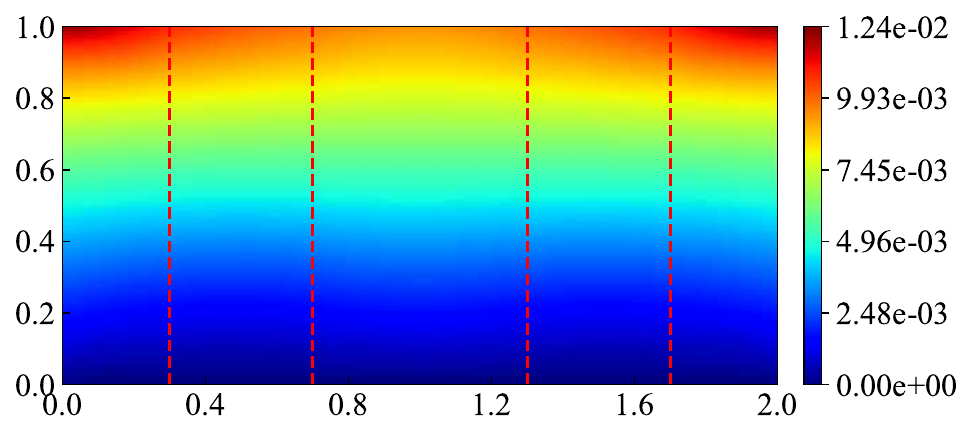}
    \caption{FEM: $u_x$}
  \end{subfigure}\hspace{0.1em}
  \begin{subfigure}{0.3\textwidth}
    \centering
    \includegraphics[width=0.95\linewidth]{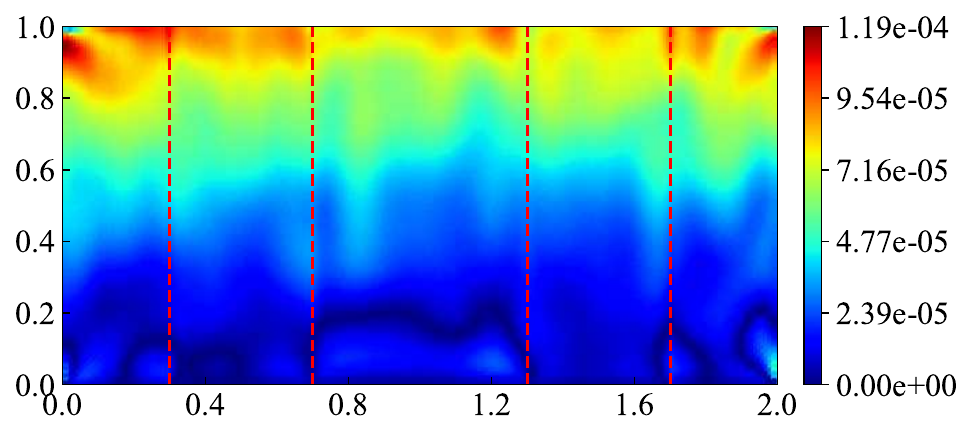}
    \caption{Absolute error: $u_x$}
  \end{subfigure}

  \par\medskip          % ← 强制换行 + 垂直间距（可改 \smallskip 或 \vspace{...}）
  % ---------- 第二行 ----------
  \begin{subfigure}{0.3\textwidth}
    \centering
    \includegraphics[width=0.95\linewidth]{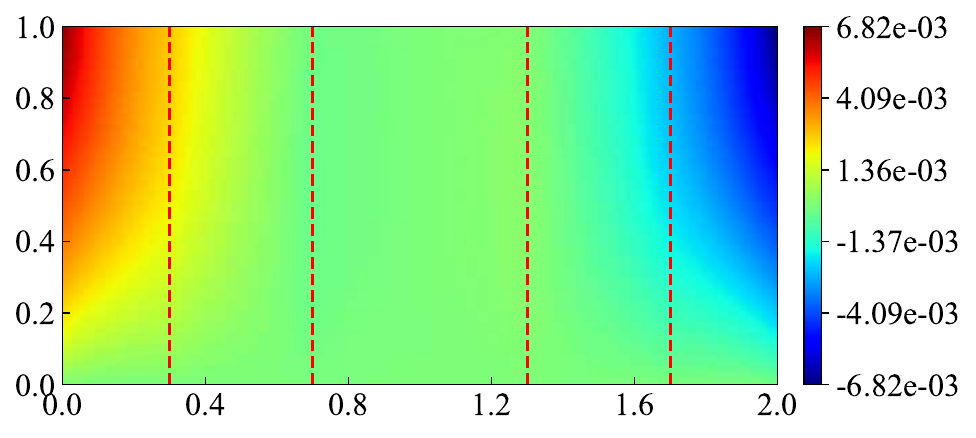}
    \caption{PIKAN: $u_y$}
  \end{subfigure}\hspace{0.1em}
  \begin{subfigure}{0.3\textwidth}
    \centering
    \includegraphics[width=0.95\linewidth]{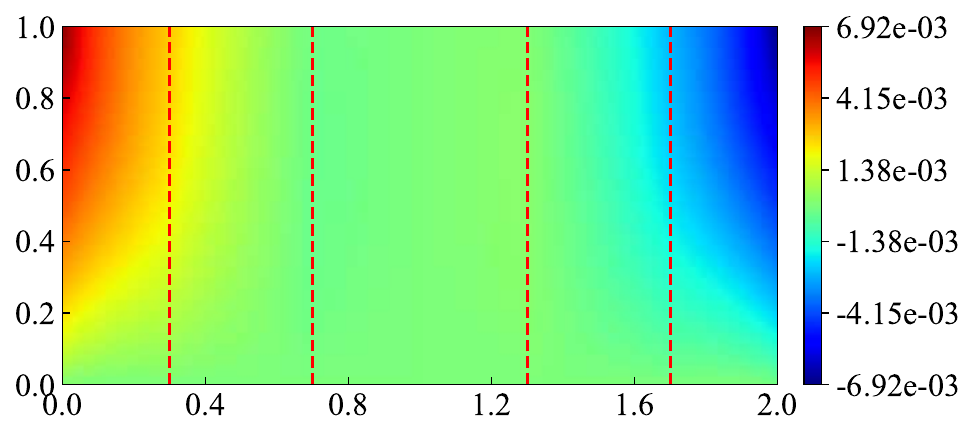}
    \caption{FEM: $u_y$}
  \end{subfigure}\hspace{0.1em}
  \begin{subfigure}{0.3\textwidth}
    \centering
    \includegraphics[width=0.95\linewidth]{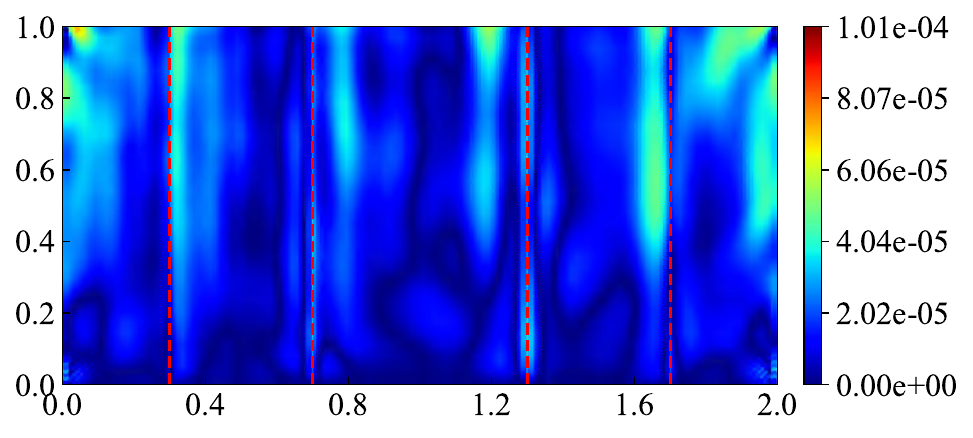}
    \caption{Absolute error: $u_y$}
  \end{subfigure}

  \caption{Comparison of PIKAN and FEM displacement fields and the corresponding absolute error in dual-via TGV-Cu structure.}
  \label{fig:tgv_dual_displacement}
\end{figure}

Fig.~\ref{fig:tgv_smooth_comparison}a - d compares PIKAN results (using triangular integration) with FEM solutions for single-via TGV-Cu structures with smooth and various roughness interfaces. Similarly, Fig.~\ref{fig:tgv_smooth_comparison}e, f shows pointwise absolute error plots for the single-via structure to quantify prediction accuracy. The displacement field distribution patterns predicted by PIKAN show excellent agreement with FEM solutions across $u_x$ and $u_y$. PIKAN successfully captures the deformation patterns in both glass and copper materials.
Fig.~\ref{fig:tgv_dual_displacement} presents the comparison between PIKAN results (using trapezoidal integration) and FEM solutions for the dual-via TGV-Cu structure. The pointwise absolute errors for the dual-via TGV-Cu structure are shown in Fig.~\ref{fig:tgv_dual_displacement}c, f. Despite the increased geometric complexity with multiple material interfaces and via interactions, PIKAN maintains high accuracy in predicting displacement fields. The method effectively captures the coupling effects between adjacent vias and the displacement distribution in the multi-via configuration.

    \begin{figure}[htbp]
    \centering
    \begin{subfigure}{0.39\textwidth}
        \centering
        \includegraphics[width=0.95\textwidth]{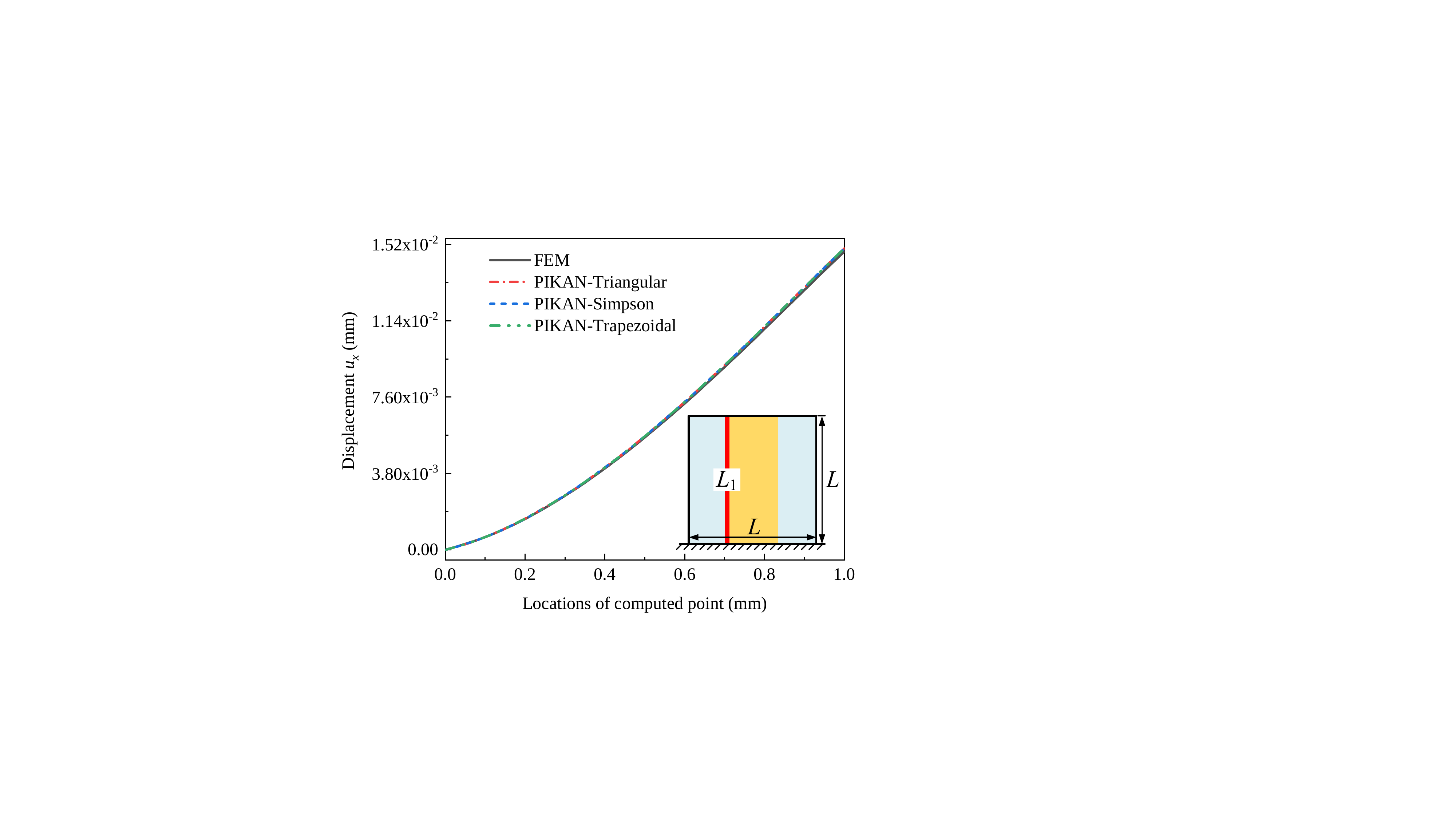}
        \caption{Smooth interface}
        \label{fig:smooth_interface_ux}
    \end{subfigure}
    %\hfill
    \begin{subfigure}{0.39\textwidth}
        \centering
        \includegraphics[width=0.95\textwidth]{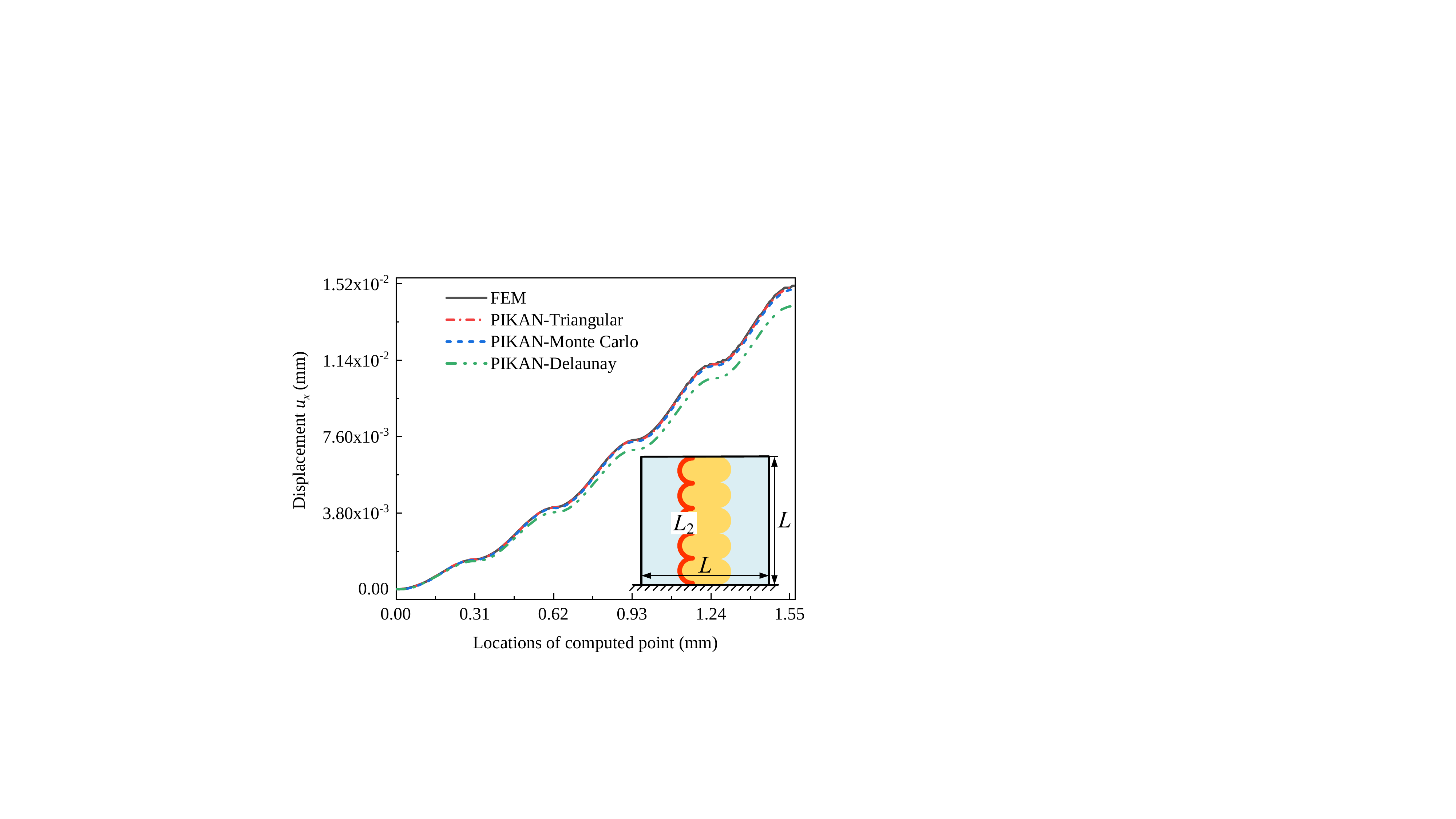}
        \caption{Interface roughness type 1}
        \label{fig:roughness1_interface_ux}
    \end{subfigure}
    % \hfill
    \begin{subfigure}{0.39\textwidth}
        \centering
        \includegraphics[width=0.95\textwidth]{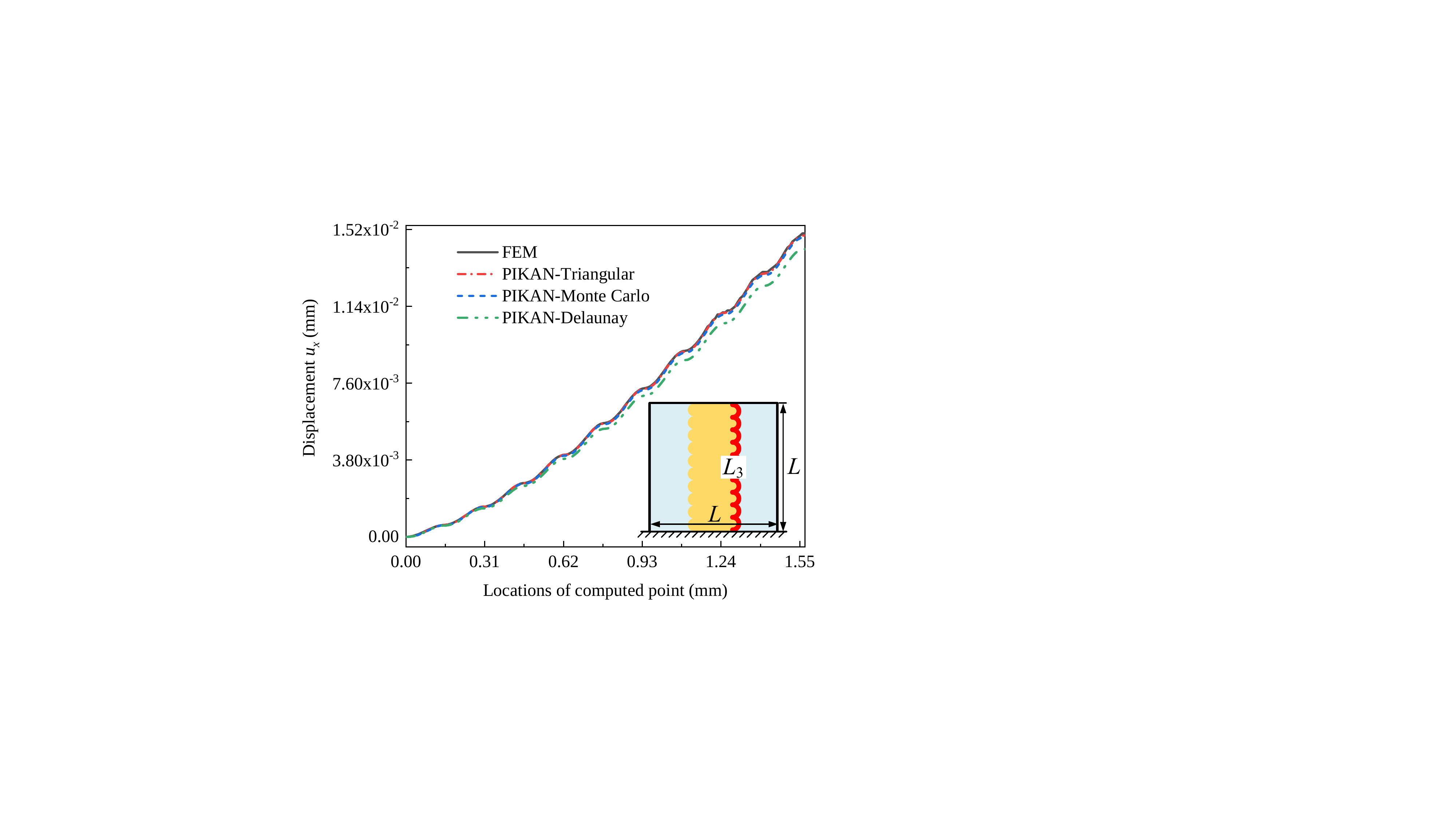}
        \caption{Interface roughness type 2}
        \label{fig:roughness2_interface_ux}
    \end{subfigure}
    % \hfill
    \begin{subfigure}{0.39\textwidth}
        \centering
        \includegraphics[width=0.95\textwidth]{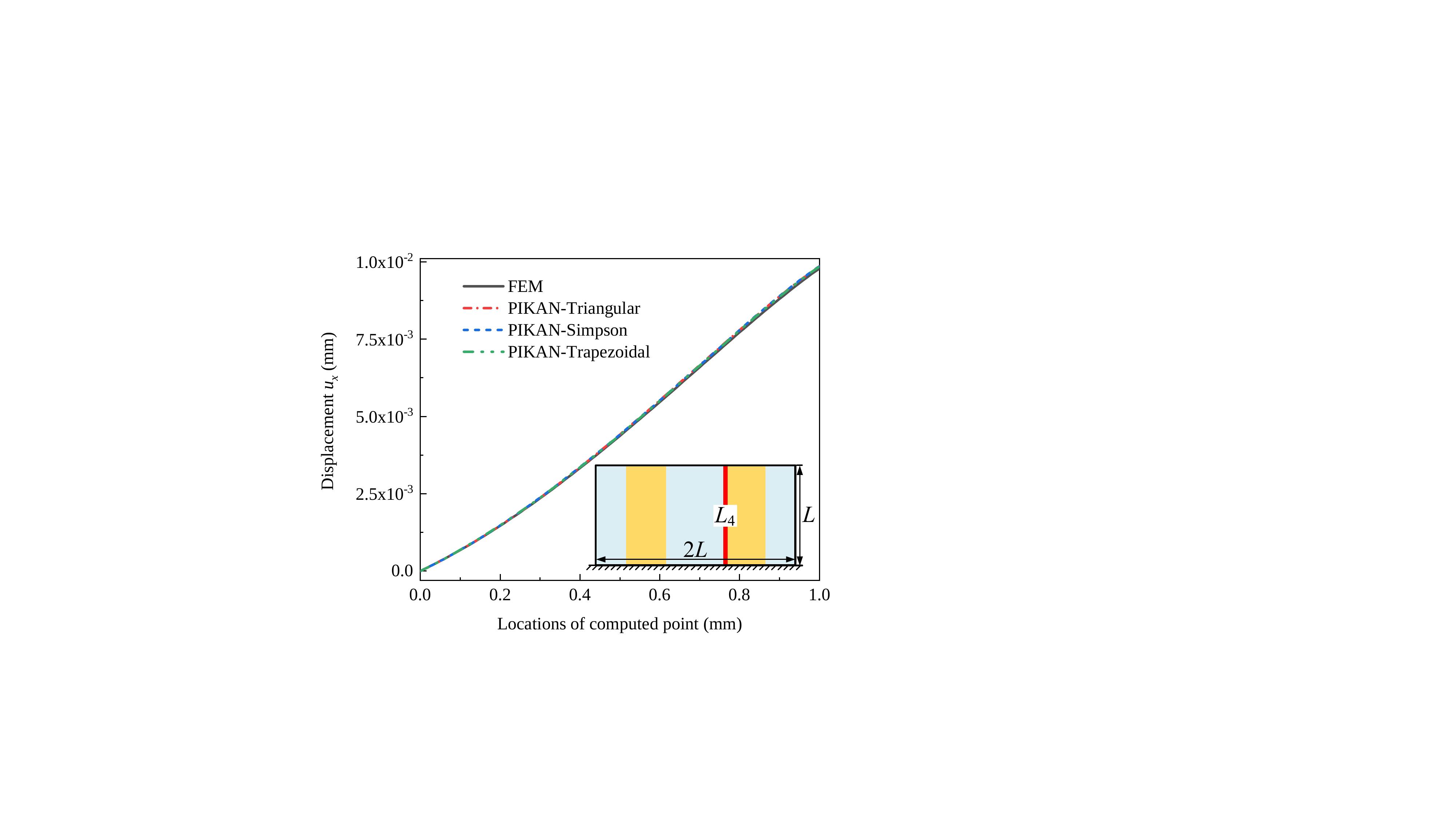}
        \caption{Dual-via configuration}
        \label{fig:dual_via_interface_ux}
    \end{subfigure}
    
    \caption{Displacement $u_x$ along interfaces for different TGV-Cu configurations (comparison of PIKAN predictions using different numerical integration schemes).}
    \label{fig:tgv_interface_comparison}
    \end{figure}

To evaluate PIKAN's performance at critical material interfaces, displacement field $u_x$ values are extracted along the glass-copper interfaces for all TGV-Cu configurations and compared with FEM reference solutions. To evaluate the method’s robustness with respect to numerical integration schemes, we systematically compare results obtained using different integration approaches in Fig.~\ref{fig:tgv_interface_comparison}. These interface analyses are particularly important for understanding stress transfer mechanisms and potential failure locations in TGV structures.
The comparison results demonstrate that PIKAN maintains excellent agreement with reference solutions across all geometric configurations, from simple smooth interfaces to complex roughness patterns and multi-via arrangements. Triangular, Monte Carlo, Simpson, and trapezoidal integration schemes all achieve high accuracy and converge well to reference solutions. However, Delaunay integration exhibits noticeable deviations from reference solutions across all curved interface geometries. This indicates that this integration scheme has inherent limitations when computing domain-wide strain energy integrals for complex geometries, with the introduced errors directly causing observable deviations in the trained network solution. 
The robust performance across different interface geometries confirms PIKAN's effectiveness for complex multi-material structural analysis without requiring specialized interface treatments or domain decomposition approaches.
    
\begin{figure}[htbp]
    \centering
    \begin{subfigure}{0.39\textwidth}
        \centering
        \includegraphics[width=0.95\textwidth]{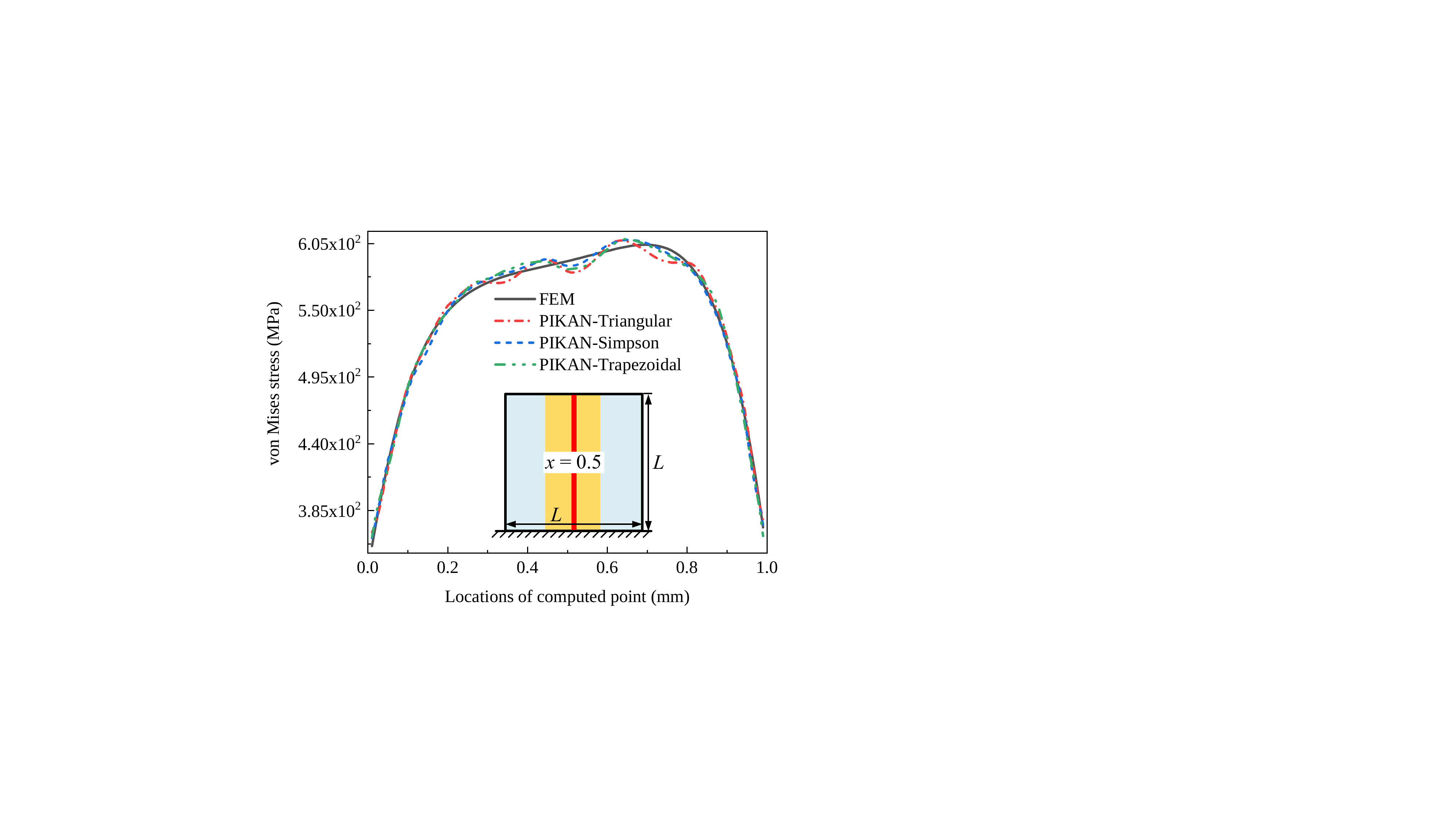}
        \caption{Copper via region: $x = 0.5$}
        \label{fig:stress_cu_domain}
    \end{subfigure}
    % \hfill
    \begin{subfigure}{0.39\textwidth}
        \centering
        \includegraphics[width=0.95\textwidth]{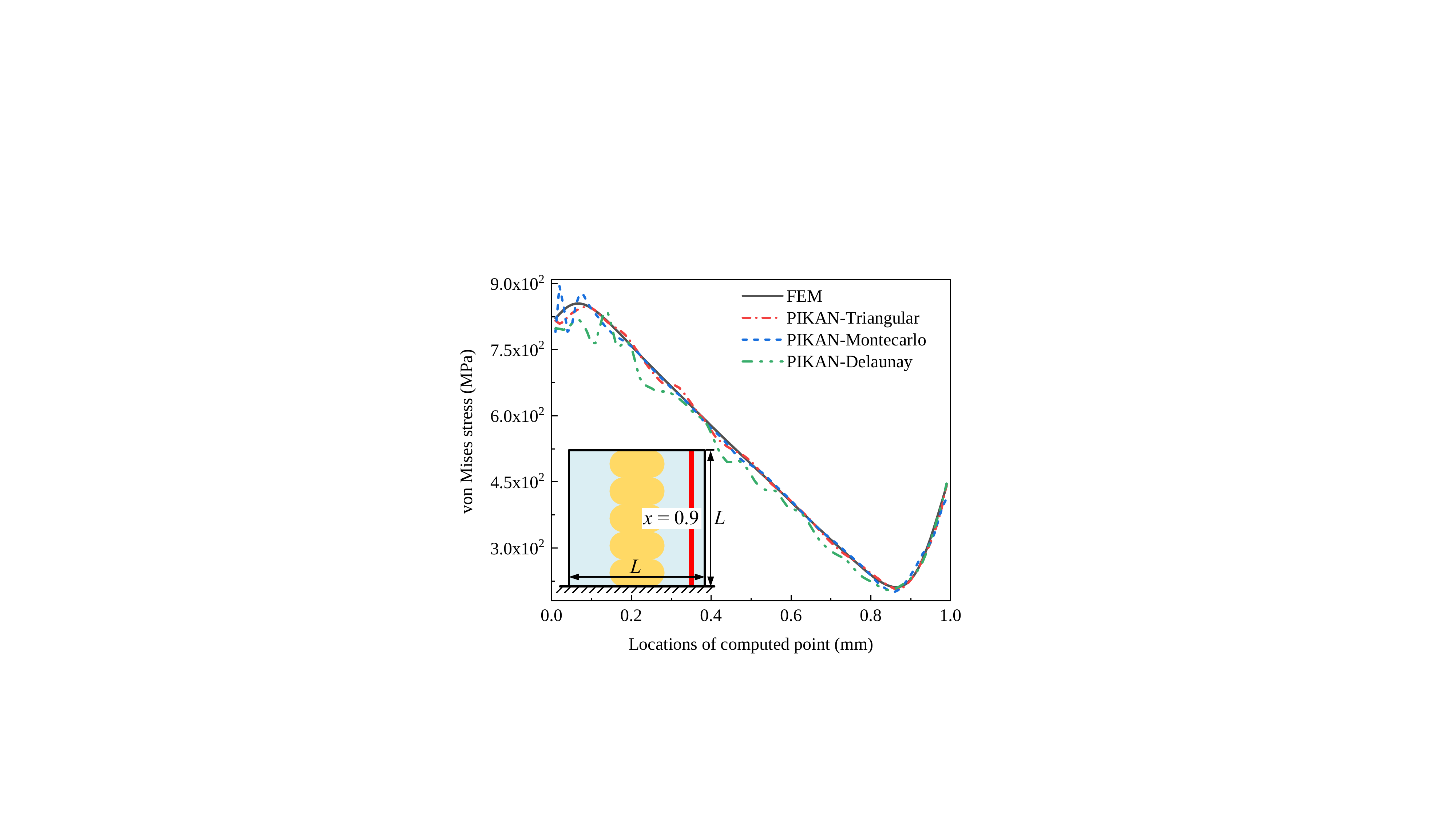}
        \caption{Glass substrate region: $x = 0.9$}
        \label{fig:stress_glass_domain}
    \end{subfigure}
    
    \begin{subfigure}{0.39\textwidth}
        \centering
        \includegraphics[width=0.95\textwidth]{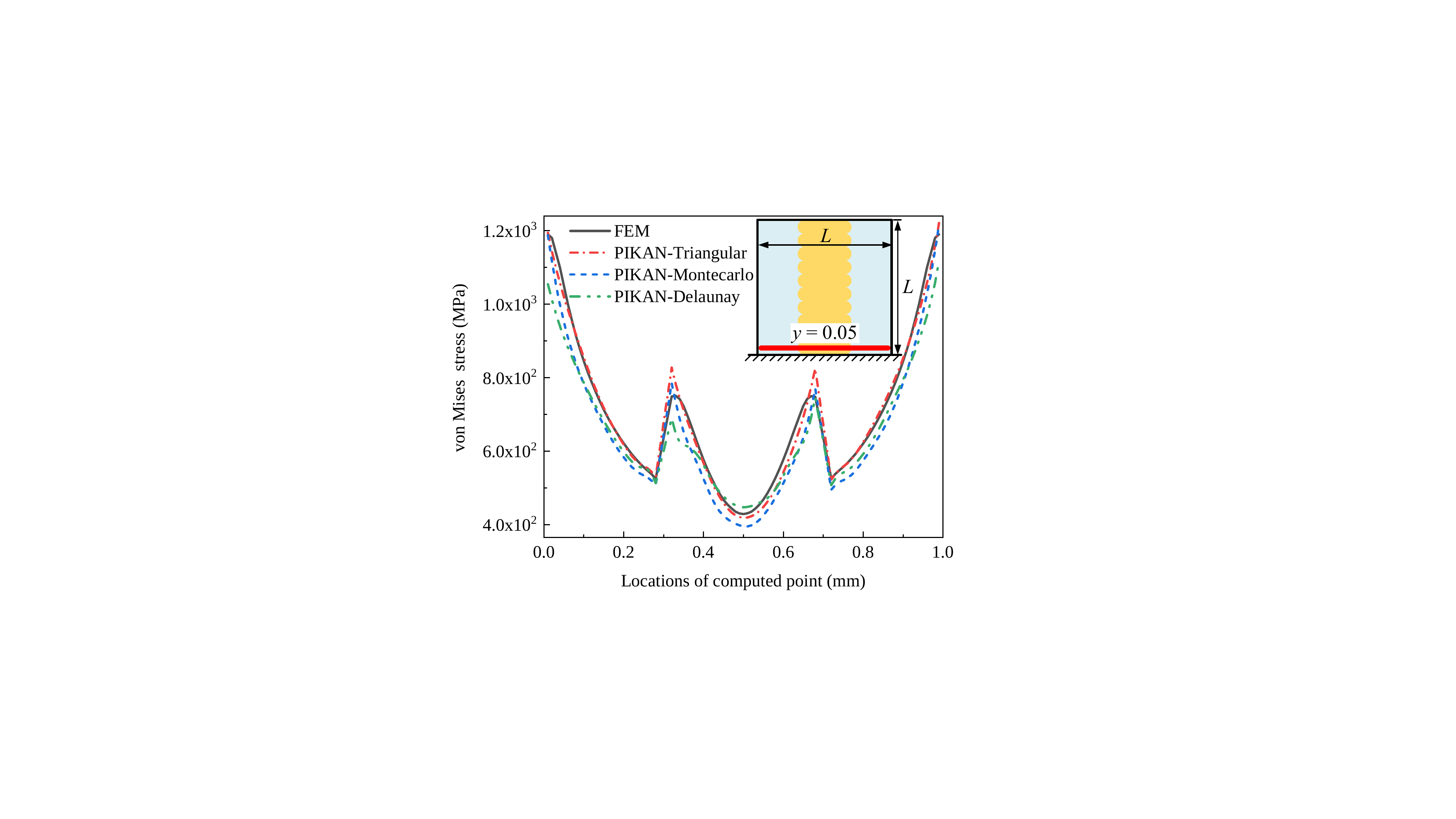}
        \caption{Single-via interface: $y = 0.05$}
        \label{fig:stress_interface_single}
    \end{subfigure}
    % \hfill
    \begin{subfigure}{0.39\textwidth}
        \centering
        \includegraphics[width=0.95\textwidth]{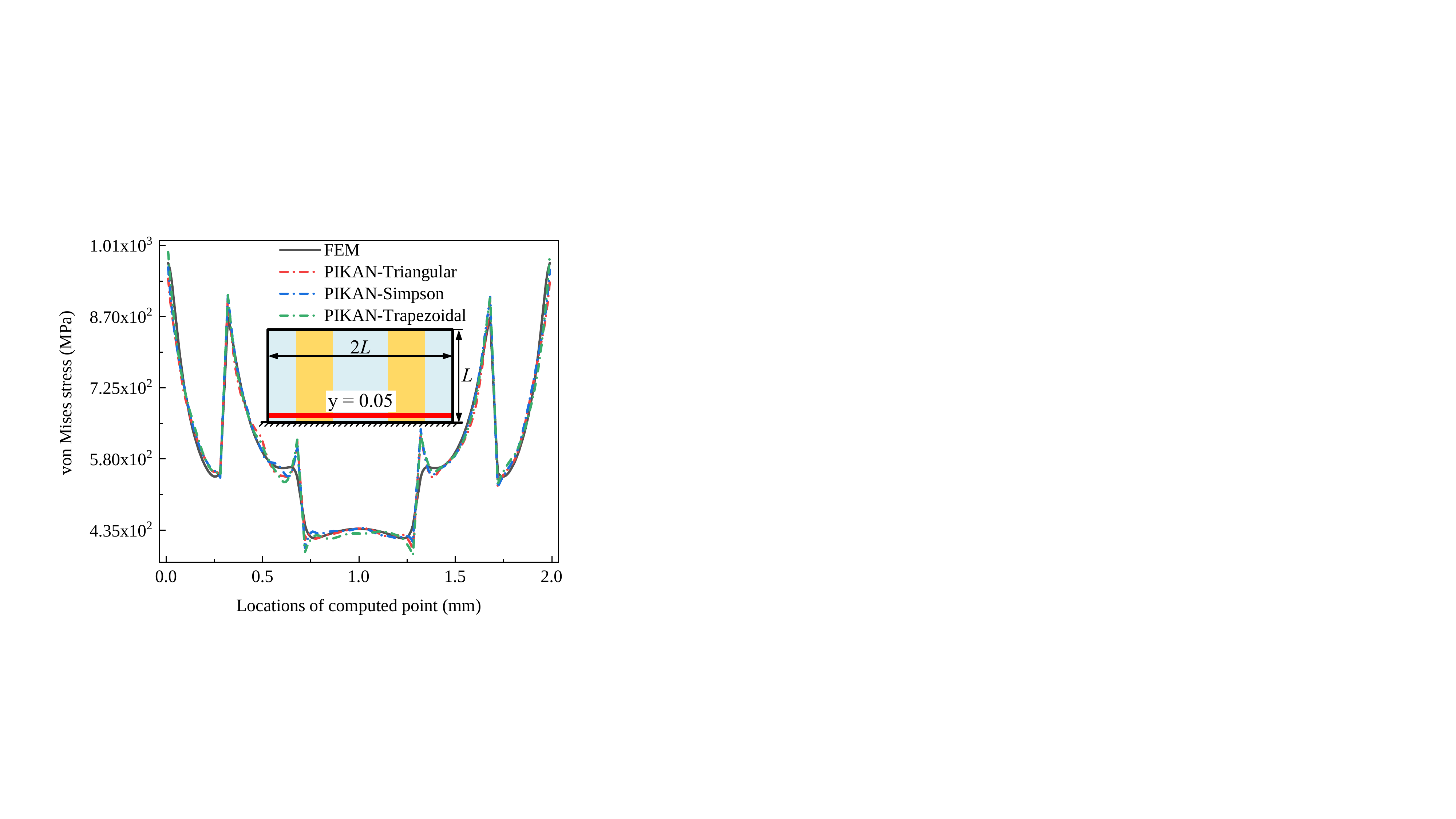}
        \caption{Dual-via interface: $y = 0.05$}
        \label{fig:stress_interface_dual}
    \end{subfigure}
    
    \caption{von Mises stress distributions at critical locations in TGV-Cu structures.}
    \label{fig:tgv_stress_comparison}
    \end{figure}

Fig.~\ref{fig:tgv_stress_comparison} presents stress distributions along strategically selected analysis lines to evaluate PIKAN's accuracy in capturing stress concentrations and material interface effects.
Figs.~\ref{fig:tgv_stress_comparison}a and b show stress distributions along vertical lines at $x = 0.5$ (within the copper via) and $x = 0.9$ (within the glass substrate). These locations capture stress variations in regions of high stress concentration (near the via) and far-field behavior (in the glass matrix). Figs.~\ref{fig:tgv_stress_comparison}c and d present stress comparisons along the horizontal line $y = 0.05$, which crosses the glass-copper interface near the bottom boundary for both single-via and dual-via configurations. This analysis line is critical for understanding stress transfer across material interfaces and potential stress concentration effects.

    The results demonstrate clear performance differences among integration schemes. Triangular integration achieves the highest accuracy due to its adaptive mesh density that better resolves stress gradients near geometric discontinuities. Although Simpson's integration method has higher theoretical precision than the trapezoidal rule, PIKAN results using both methods achieve similar accuracy levels when the sample point density is sufficiently high. Monte Carlo integration performs moderately well but shows some scatter typical of stochastic methods. Delaunay integration exhibits the poorest numerical stability, especially at interface discontinuities. Its suboptimal control area allocation introduces non-negligible discretization errors into the loss function, adversely affecting training process robustness and convergence.

\section{Conclusion}
\label{sec:conclusion}

This study proposes PIKAN, a novel method combining KAN with the DEM for analyzing multi-material elasticity problems in electronic packaging structures. KAN networks inherently possess piecewise function characteristics that naturally accommodate material property discontinuities, while DEM provides superior computational efficiency and accuracy compared to strong-form PINNs. The PIKAN method effectively integrates these advantages to achieve accurate analysis of complex multi-material structures. Unlike traditional approaches that introduce penalty terms to enforce essential boundary conditions, PIKAN ensures solution validity by constructing admissible displacement fields that automatically satisfy essential boundary conditions. The method requires only a single KAN network to approximate displacement fields across the entire computational domain, eliminating the subdomain partitioning required in conventional multi-material analysis. This approach avoids hyperparameter tuning for interface continuity conditions, reducing training complexity while achieving superior accuracy compared to domain decomposition methods.

Numerical validation demonstrates PIKAN's accuracy and robustness across diverse multi-material configurations. Comparative analysis with traditional CENN approaches confirms PIKAN's superior performance, particularly for complex interface geometries where penalty-based methods struggle with convergence. The study reveals that numerical integration strategy affects solution accuracy. Additionally, PIKAN requires fewer trainable parameters compared to MLP-based methods. Future work will systematically compare the performance of KANs with different basis functions (such as Chebyshev polynomials) for multi-material problems~\cite{SHUKLA2024117290}. Additionally, exploring hybrid MLP-KAN architectures represents a promising direction~\cite{kashefi2025}. For instance, employing MLPs for global approximation while integrating KANs in the output layer to capture local discontinuities with high precision could achieve an optimal balance between efficiency and accuracy. Such hybrid approaches may leverage the complementary strengths of both architectures: MLPs’ efficiency in smooth regions and KANs’ superior ability to handle sharp discontinuities. At the theoretical level, future research will focus on establishing rigorous mathematical foundations for the proposed PIKAN method, with emphasis on rigorous convergence and stability analysis of PIKAN for solving solid mechanics problems. Furthermore, the energy-based framework has been proven applicable to multiphysics and nonlinear problems~\cite{WOS:001115578900001,WOS:001314528200001}, while the KAN architecture has demonstrated strong potential in solving complex differential equations~\cite{SHUKLA2024117290,toscano2025pinns,patra2024physics}. Therefore, the proposed method establishes a solid foundation for addressing broader nonlinear and multiphysics problems, which represents an important future research direction. The combination of DEM’s variational formulation with KAN’s superior nonlinear representation capability positions PIKAN as a promising approach for tackling increasingly complex engineering problems.

The current study focuses on 2D validation examples. Extension to 3D electronic packaging structures presents significant challenges beyond simply increasing problem dimensionality. To begin with, for analyzing large-scale electronic packaging structures with multi-scale features, further enhancing computational efficiency and scalability will be critical. Future research needs to develop adaptive sampling and training strategies to efficiently capture physical field details in critical regions, and to design specialized optimization techniques for large-scale KAN networks to improve training stability and convergence speed. These advances will enable PIKAN to tackle increasingly complex, real-world electronic packaging problems while maintaining its core advantage of seamless multi-material interface handling. Furthermore, recent research has revealed that standard KANs face theoretical bottlenecks in effectively capturing E(3) transformations, which is essential for representing physical laws in 3D space~\cite{alesiani2025geome}, and their parameter efficiency may degrade poorly with increasing dimensionality~\cite{WOS:001542855900001}. These limitations are particularly concerning for 3D solid mechanics problems where computational efficiency and physical consistency are critical. Despite these challenges, recent studies have begun addressing high-dimensional and symmetry-related issues through improved KAN architectures~\cite{hu2025incorp} or by integrating KANs into 3D physical modeling frameworks~\cite{huang2025dreamphysics}, providing valuable references for future research. However, developing efficient, physically-meaningful PIKAN solvers suitable for 3D solid mechanics problems remains an open research question and represents a key focus for our future work.

In conclusion, this research successfully extends PINNs to heterogeneous material problems, demonstrating significant potential for electronic packaging structure reliability analysis and broader multi-material engineering applications. Source codes are available at https://github.com/yanpeng-gong/PIKAN-MultiMateria.

\section*{CRediT authorship contribution statement}
\textbf{Yanpeng Gong}: Conceptualization, Methodology, Formal analysis, Writing – original draft, Writing – review \& editing, Supervision, Project administration, Funding acquisition. \textbf{Yida He}: Software, Validation, Formal analysis, Data curation, Writing – original draft, Visualization. \textbf{Yue Mei}: Supervision, Writing – review \& editing, Supervision. \textbf{Fei Qin}: Supervision, Writing – review \& editing, Supervision. \textbf{Xiaoying Zhuang}: Writing – review \& editing, Supervision. \textbf{Timon Rabczuk}: Writing – review \& editing, Supervision.

\section*{Acknowledgments}
This research was supported by the National Natural Science Foundation of China (Nos.12002009, 12472198, W2431038), Beijing Natural Science Foundation (L233001) and the Alexander von Humboldt Foundation, Germany. The authors would like to thank Mohammad Sadegh Eshaghi for his valuable suggestions and guidance during the algorithm implementation process of this work.

\hypertarget{AppendixA}{\section*{Appendix A}}
\label{AppendixA}

\begin{algorithm}[htbp]
    \scriptsize
    \caption{PIKAN algorithm for multi-material problems}
    \label{alg:pikan_multimaterial}
    \begin{algorithmic}[1]
    \STATE \textbf{Input:} 
    \STATE \quad Physical domains of different materials $\Omega_1, \Omega_2, \ldots, \Omega_n$
    \STATE \quad Boundary conditions: essential boundary $\Gamma_u$ and natural boundary $\Gamma_\sigma$
    \STATE \quad Material parameters: $E_1, \nu_1, E_2, \nu_2, \ldots, E_n, \nu_n$
    \STATE \quad Domain sample points $\mathbf{x}$ from all material domains
    \STATE \quad Natural boundary condition points $\mathbf{x}_\sigma$ from $\Gamma_\sigma$
    \STATE \quad KAN architecture and hyperparameters
    \STATE \quad Neural network optimizer
    \STATE \textbf{Output:} Optimized KAN parameters $\boldsymbol{\theta}^*$
    
    \STATE \textbf{Initialize:} KAN parameters $\boldsymbol{\theta} = \{c_p^{(j,i)}, m_{j,i}, w_{j,i}\}$
    \STATE Normalize all sample point coordinates
    
    \WHILE{Loss function not converged}
        \STATE Obtain displacement predictions $\mathbf{F}(\mathbf{x})$ from single KAN for all sample points
        \STATE Construct admissible displacement field $\mathbf{u}^{\text{pred}}(\mathbf{x})$
        \STATE Compute displacement gradients $\nabla \mathbf{u}^{\text{pred}}$ using automatic differentiation
        \FOR{each material domain $i = 1, 2, \ldots, n$}
            \STATE Compute strain energy $\Psi_{\text{in}}^{m_i}$ using material parameters $(E_i, \nu_i)$
        \ENDFOR
        \STATE Compute external work potential energy: $\Psi_{\text{ex}}$
        \STATE Calculate loss function: $\mathcal{L}_{\text{PIKAN}} = \sum_{i=1}^{n} \Psi_{\text{in}}^{m_i} - \Psi_{\text{ex}}$
        \STATE Update KAN parameters $\boldsymbol{\theta}$ using optimizer and backpropagation
    \ENDWHILE
    
    \STATE \textbf{Return:} Optimized parameters $\boldsymbol{\theta}^*$
    \end{algorithmic}
\end{algorithm}

\hypertarget{AppendixB}{\section*{Appendix B: Hyperparameter sensitivity analysis of PIKAN}}
\label{AppendixB}

\begin{figure}[htbp]
\centering
\begin{subfigure}{0.45\textwidth}
    \centering
    \includegraphics[width=\textwidth]{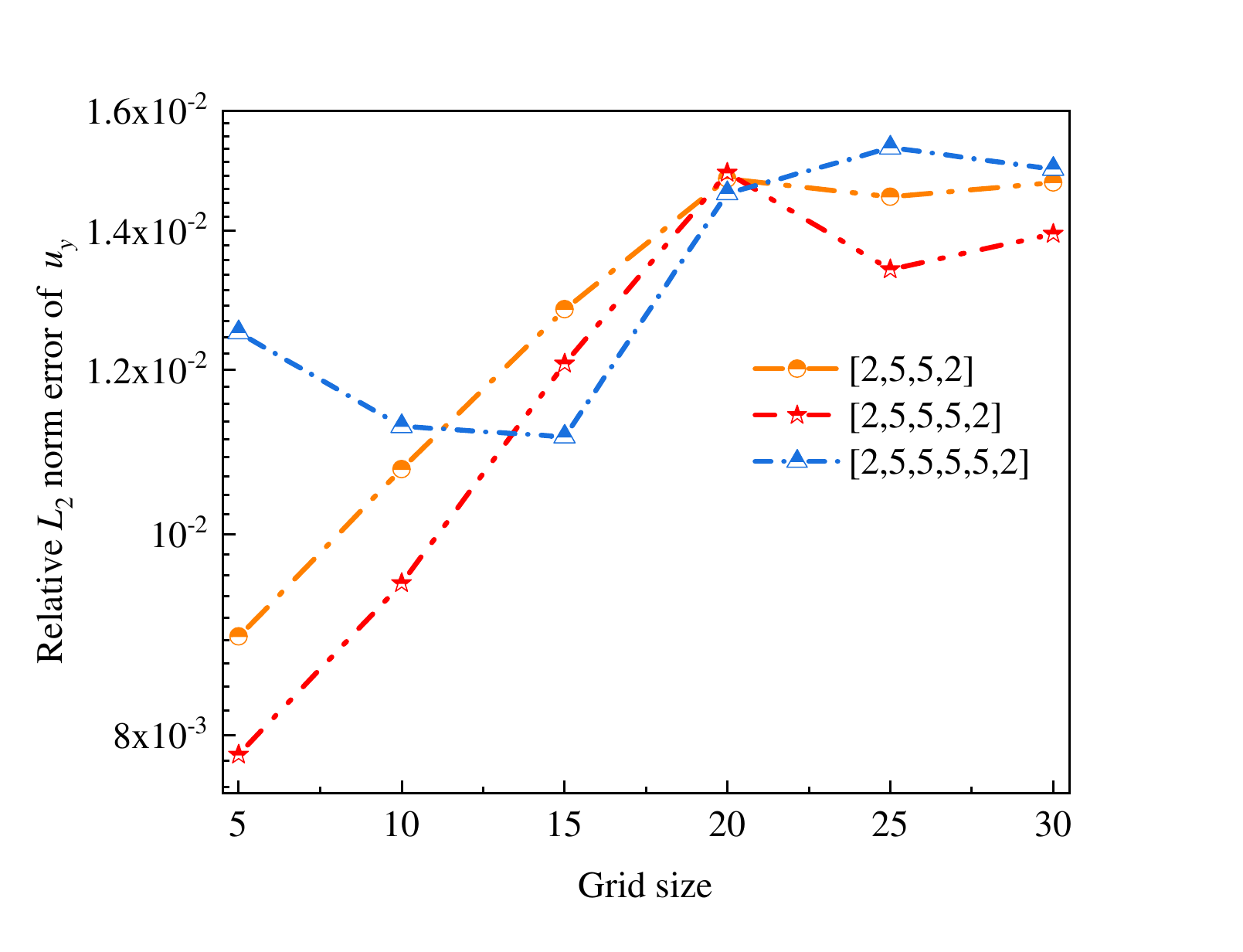}
    \caption{Grid size effect (B-spline order = 3)}
    \label{fig:grid_size_effect}
\end{subfigure}
\hspace{1em}
\begin{subfigure}{0.45\textwidth}
    \centering
    \includegraphics[width=\textwidth]{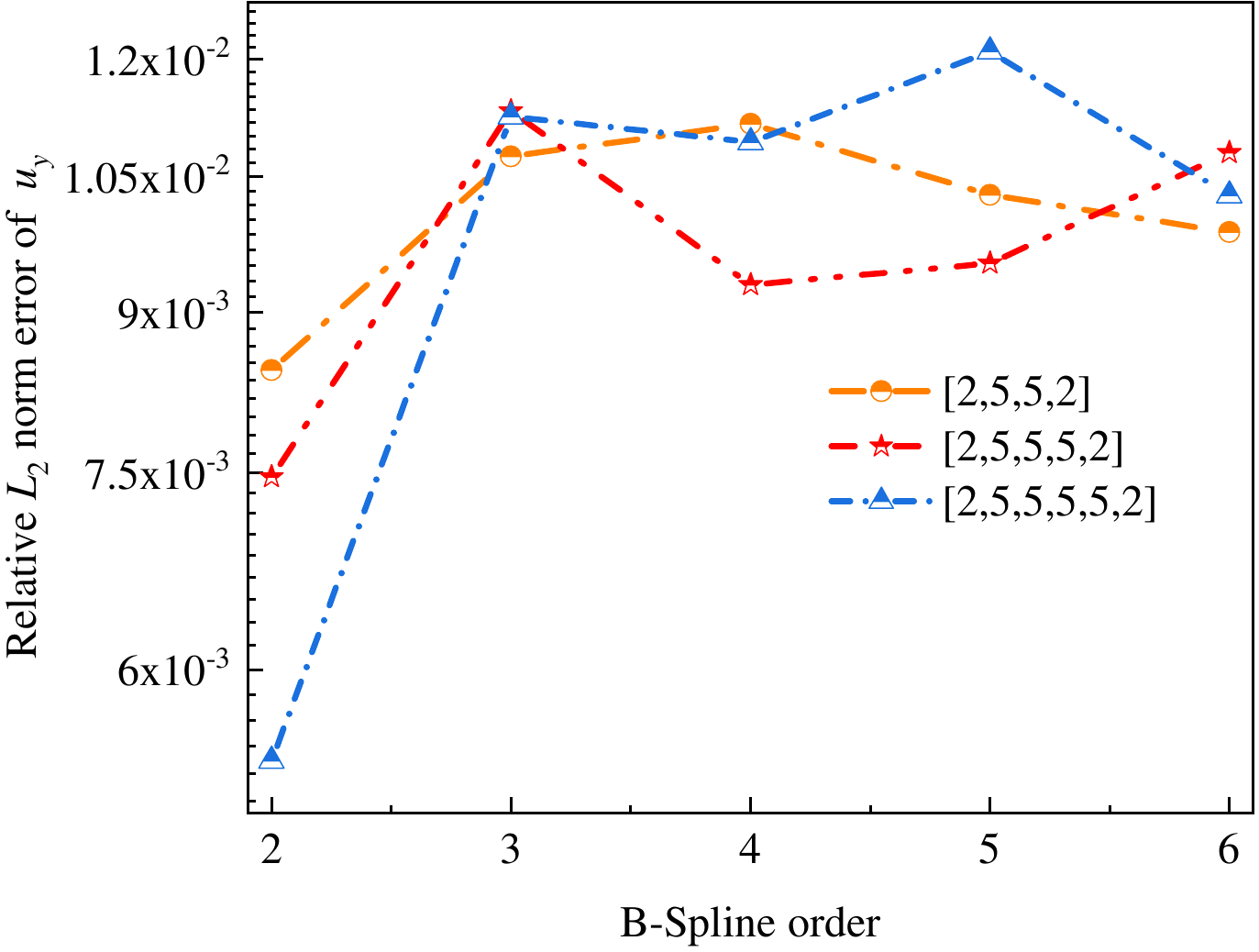}
    \caption{B-spline order effect (grid size = 10)}
    \label{fig:bspline_order_effect}
\end{subfigure}
\caption{Relative $L_2$ norm error of $u_y$ for different KAN hyperparameter combinations.}
\label{fig:hyperparameter_analysis}
\end{figure}

Hyperparameter selection critically affects PINN computational results. This section investigates how different hyperparameter combinations impact PIKAN prediction accuracy using the straight interface cantilever beam model from Fig.~\ref{fig:cantilever_straight}. We analyze the relative $L_2$ norm error variation for displacement component $u_y$ under different KAN configurations, where each configuration is defined by network depth and width, grid size, and B-spline order.
Fig.~\ref{fig:hyperparameter_analysis} presents the hyperparameter sensitivity analysis results. Fig.~\ref{fig:grid_size_effect} shows the $L_2$ norm error of $u_y$ with B-spline order fixed at 3 for three different KAN architectures: [2,5,5,2], [2,5,5,5,2], and [2,5,5,5,5,2]. The results indicate that increasing grid size does not produce consistent accuracy improvements across different network depths, suggesting that optimal grid size depends on the specific network architecture.
Fig.~\ref{fig:bspline_order_effect} demonstrates the effect of varying B-spline order on computational accuracy when grid size is fixed at 10. The results show that the optimal order varies with network architecture. These findings highlight the importance of joint optimization of KAN hyperparameters rather than individual parameter tuning.

The analysis reveals complex relationships between KAN hyperparameters and computational accuracy. While smaller grid sizes and B-spline orders often yield better performance in these specific examples, the optimal hyperparameter combination varies with network architecture and problem complexity. This sensitivity occurs because non-optimal hyperparameter selections can lead to either underfitting (insufficient model capacity to capture solution complexity) or overfitting (excessive model complexity leading to poor generalization), consistent with findings reported in machine learning literature~\cite{liu2024kan}. 
The non-monotonic relationship between hyperparameters and accuracy highlights a fundamental challenge in KAN optimization: larger grid sizes or higher B-spline orders do not automatically improve performance. Instead, these parameters must be carefully balanced with network depth and width to achieve optimal results. Currently, effective methods for global optimization of KAN hyperparameter combinations are lacking due to the vast search space and extensive computational requirements for systematic exploration. Future work should focus on developing efficient hyperparameter optimization strategies specifically tailored for physics-informed KAN applications.




%\begin{thebibliography}{00}

%% For numbered reference style

%   \textit{\LaTeX: a document preparation system},
%   Addison Wesley, Massachusetts,

% \end{thebibliography}

% 在文档末尾，\end{document} 之前添加：
\bibliographystyle{elsarticle-num}  % 使用elsarticle数字风格
\bibliography{PINN_ref_lib.bib}  % 替换为您的.bib文件名
\end{document}